\RequirePackage{snapshot}
\documentclass[final]{scrartcl}
\usepackage[english]{babel}
\usepackage[utf8]{inputenc}
\usepackage[margin=2cm]{geometry}

\usepackage{mathtools}

\usepackage{amssymb}  

\usepackage{bm}
\usepackage{amsthm}
\usepackage{commath}

\usepackage[square,numbers]{natbib}
\usepackage[final]{graphicx}
\usepackage[dvipsnames]{xcolor} 

\usepackage{booktabs}

\usepackage{showlabels}
\usepackage{rotating}

\usepackage{placeins} 
\usepackage{subcaption}
\usepackage{framed}
\usepackage[normalem]{ulem}
\usepackage{enumerate}
\usepackage{siunitx}
\usepackage[ruled,vlined]{algorithm2e}
\SetKw{KwBy}{by}

\usepackage{textcomp,soul}

\usepackage[version=4]{mhchem}
\usepackage[T1]{fontenc}

\usepackage[final]{hyperref}
\hypersetup{
    colorlinks,
    linkcolor={red!50!black},
    citecolor={blue!50!black},
    urlcolor={blue!80!black}
}

\usepackage{url}

\usepackage{cleveref}
\graphicspath{{./images/}}

\newcommand{\partialderiv}[2]{\frac{\partial #1}{\partial #2}}

\newcommand{\pardiff}[2]{\frac{\partial #1}{\partial #2}}
\newcommand{\dive}[1]{\bm{\nabla}\cdot\left(#1\right)}
\newcommand{\lapl}[1]{\nabla^{2}#1}
\newcommand{\evalu}[3]{\bigg[#1\bigg]^{#2}_{#3}}

\newcommand{\secpardiff}[2]{\frac{\partial^2 #1}{\partial #2^2}}
\newcommand{\diff}[2]{\frac{d #1}{d #2}}
\newcommand{\secdiff}[2]{\frac{d^2 #1}{d #2^2}}

\newcommand{\of}{OpenFOAM\textsuperscript\textregistered\,}
\newcommand{\RE}{\textrm{Re}}
\newcommand{\PE}{\textrm{Pe}}
\newcommand{\DA}{\textrm{Da}}

\newcommand{\N}{\mathcal{N}}
\renewcommand{\L}{\mathcal{L}}
\renewcommand{\P}{\mathcal{P}}

\newcommand{\bx}{\mathbf{x}}
\newcommand{\flux}{\mathbf{j}}
\newcommand{\bu}{\mathbf{u}}

\title{Electrochemical transport modelling and open-source simulation of pore-scale solid-liquid systems}
\author{
Robert Barnett\footnote{School of Mathematical Sciences, University of Nottingham, University Park, NG7 2RD, Nottingham, UK},
Federico Municchi\footnote{Colorado School of Mines,
1500 Illinois St.,
Golden,
CO
80401,
USA},
John King$^*$, Matteo Icardi$^*$\footnote{Corresponding author, \url{matteo.icardi@nottingham.ac.uk}}
}

\date{\today}

\begin{document}

\maketitle

\section*{Abstract}

    The modelling of electrokinetic flows is a critical aspect spanning many industrial applications and research fields. This has introduced great demand in flexible numerical solvers to describe these flows. The underlying phenomena are microscopic, non-linear, and often involve multiple domains. Therefore often model assumptions and several numerical approximations are introduced to simplify the solution. In this work, we present a multi-domain multi-species electrokinetic flow model including complex interface and bulk reactions. After a dimensional analysis and an overview of some limiting regimes, we present a set of general purpose finite-volume solvers, based on \of, capable of describing an arbitrary number of electrochemical species over multiple interacting (solid or fluid) domains \cite{spnpfoam}. We provide verification of the computational approach for several cases involving electrokinetic flows, reactions between species, and complex geometries. We first present three one-dimensional verification test cases, and then show the capability of the solver to tackle two- and three-dimensional electrically driven flows and ionic transport in random porous structures.  The purpose of this work is to lay the foundation for a general-purpose open-source flexible modelling tool for problems in electrochemistry and electrokinetics at different scales.

\section*{Keywords}
Stokes Poisson Nernst Planck, dilute electrolyte, OpenFOAM, electrokinetic flow, ionic transport    


\section{Introduction}

    Electrokinetic flows are a highly active topic of discussion that reaches a multitude of scientific fields. Examples include chloride transport in reinforced concrete \cite{Chlo_Ingres_Mortar,Refin_PNP_Conc}, ion regulation in biological cells \cite{BPNP_Model_Q_Zheng,PNPF_Theo_Ion_Channels}, fuel cells \cite{Ceramic_Fuel_Cells} and electrochemical energy storage \cite{Charge_Trans_Modelling_Review_Richard,Naf_Membrane}, such as batteries and super-capacitors \cite{Compare_SM_PNP,Therm_Dyn_Cons_PNP_A_Latz,Newman_Revisit}. Consequently, the demand for the numerical modelling of electrokinetic flows, often due to complex geometries or multi-physics barring analytical solutions, spans a great deal of fields. While each real-world example given often comes with its own bespoke problems to consider, they are all describable as specific cases under a generalized model for electrokinetic flow problems. Therefore, the motivation for this paper is to lay the foundational work of a general-purpose open-source computational fluid dynamics (CFD) toolbox for the modelling of electrokinetic flows \cite{spnpfoam}. Built-in a modular and sequential fashion so additional physics, like steric effects \cite{Ji2012,Ma2020,Voukadinova2019} or chemical activity \cite{Samson1999}, can be subsequently 'grafted' onto the workflow and solved too. As such, we chose the finite-volume CFD package \of as our underlying CFD package for implementation. The reasons being: its open-source nature means development is unrestricted and the package is widely accessible; despite its age, \of still retains a large community of both industry and academic users that may find value in such solvers; while its numerical approach to coupled multi-physics system, handled in a sequential, i.e., segregated manner, makes inclusion of other electrokinetic phenomena not introduced here a simpler procedure -- albeit with caveats such as stability concerns -- so extensions are more permissible.
    
    With electrokinetic flows, the involvement of electrical forces leads to a number of interesting phenomena. One such phenomenon is electro-osmosis \cite{Alizadeh2021}, where an applied electric field induces fluid movement due to the formation of electric double layers, see \cref{fig::EDL}. Another interesting set of phenomena is known as induced-charge electrokinetics. Whilst similar to traditional electrokinetic phenomena, the difference comes from the double layer being induced by an applied electric field \cite{Squires2009}. The main difficulty of modelling electrokinetic flows stems from the microscopic scale of phenomena. For applications with large-scale domains, solving at the microscopic scale can be computationally taxing. For systems with complex geometries, such as porous media, this difficulty is furthered.
    
    To describe electrokinetic flows we require modelling the ion concentrations, electric potential and fluid velocity, with the Stokes-Nernst-Planck-Poisson (SPNP) model. In this work, for the velocity field we will, in fact, consider Stokes' flow as, at the micro- and meso-scale, the viscous forces far outweigh the inertial due to small velocities and length scale \cite{Mark_Schmu_PNP_Porous}. The Stokes equation is coupled with Poisson's equation for the electric potential, relating the electric field composition to the variation in ion charge density. Finally, these are both coupled with the Nernst-Planck equation to describe ion transport. In using Nernst-Planck, we neglect any ion-ion interaction by assuming the ionic solution is sufficiently dilute. 
    The Nernst-Planck ionic flux was first formulated in its steady form for a one-dimensional cylinder by Nernst \cite{Nernst1888}. It was later extended by Planck \cite{Planck1890} to a transient setting, furthermore introducing the continuity and Poisson's equations \cite{Maex2017}. In doing this the two paved the way in helping develop SPNP, providing a simple yet accurate description of electrokinetic flow for dilute solutions.     

    A common similarity in many applications containing electrokinetic flows is the involvement of multiple regions, such as fluid electrolytes and solid electrodes in batteries. These regions exchange ions with each other and in some instances contain chemical reactions exchanging mass between species. For batteries, these reactions are a necessary process for operational use. However, in examples such as chloride corrosion in reinforced concrete, it is a detriment \cite{Sim_Chl_Mig_Rein_Conc}. As such, modelling these reactions has just as much importance as the flow they reside in.

    Whilst used extensively in a wide range of physical settings by the research community, SPNP does come with its own caveats. For one, SPNP neglects any ion-ion interactions that may occur by assuming a dilute solution. This may not hold true for solutions with many ionic species \cite{Compare_SM_PNP}. Also, as SPNP is a continuum model, any steric effects are ignored. As such, many efforts have been made to extend SPNP to include other physical processes. To cover steric effects, free energy functionals using density functional theory (DFT) accounting for long-range Coulomb correlation and hard sphere (HS) interactions of ions \cite{Ji2012,Ma2020,Voukadinova2019} are formulated. Extensions to make Nernst-Planck more thermodynamically consistent under non-equilibrium thermodynamics \cite{Therm_Dyn_Cons_PNP_A_Latz,Dreyer2013} have also been proposed. To model non-ideal solutions, \cite{Samson1999} proposes an added term to the Nernst-Planck ionic flux considering varying chemical activities solved by the Debye-H\"uckel model.
    
    Whilst these extensions do further the physical realism of the original SPNP model, this often exacerbates other challenges of SPNP. One of the foremost is the non-linear coupling between fields. For example, in \cite{Samson1999} the Debye-H\"uckel model equates the chemical activity of a species to the solution's ionic strength. This in part creates explicit coupling between all ionic concentrations unlike in classic SPNP, resulting in more complex, often unviable,  computational approaches. 
    
    The development of numerical solvers for such equations within general PDE and CFD toolboxes is something that has been discussed for decades.
    Two common approaches can be taken to solve systems of coupled discretized equations. The first is the so-called block-coupled, with all equations solved at once in a large matrix. Whilst taking a large amount of memory, it upholds the coupling between fields and is numerically robust. The second is the segregated approach and consists in solving each equation separately and in sequence. Since this leads to a decoupling of the equations, appropriate iterative methods \cite{Patankar1972,Nuca2022}  must be used to ensure coupling between fields. The advantages of a segregated approach are the lower memory requirements, easy preconditioning of the equations, and their multi-stage structure that allows better control of the solution procedure. However, block-coupled approaches tend to scale better with the number of processors. When constructing our solvers to model electrokinetic flows, we chose a segregated iterative approach to couple the equations and the different domains.
    
    This work presents a multi-region multi-species SPNP model and discusses its implementation in finite-volumes segregated solvers, built with the \of library and released open-source \citep{spnpfoam}. We present the mathematical model, including a dimensional analysis, and consider multiple solid and fluid regions, with general reaction and interface models.    
    Whilst other  finite-volumes and finite-elements solvers have been developed \cite{Pimenta2019,Ying2021,Zheng2011,Basu2020,Li2020}, restrictions such as being designed for specific applications, single domains, dimensionality, absence of reactions, steady state or ignoring the fluid velocity are often made. 
    %
    %
    Another point of novelty also stems from a general non-linear reaction model, so that various reaction rate models such as Butler-Volmer \cite{Charge_Trans_Modelling_Review_Richard} or the \textit{rate law} \cite{RaymondChang2000} can be applied.
    
	
    
    This work is organised in the following way. Within \S\ref{sec::SPNP} we present the governing equations of Stokes-Poisson-Nernst-Planck and what fluid properties are assumed to give accurate flow description. In \S\ref{sec::dimensionalAnalysis} we perform dimensional analysis to understand the transport regimes possible and how this results in the often-used electro-neutrality approximation. For \S\ref{sec:multidomain} we outline what is required to capture reactions at a multi-region interface given different restrictions on ion movement. In \S\ref{sec:numericalImplementation} the implementation of our solvers for single and multi-region is discussed, as well as the iterative algorithm performed when introducing our reactive conditions. To verify the accuracy of our solvers and reactive conditions, we provide necessary numerical examples in \S\ref{sec:results}, with concluding remarks given in \S\ref{sec::Conclusion}.
    
\section{Stokes-Poisson-Nernst-Planck model}\label{sec::SPNP}

\begin{figure}
        \centering
        \includegraphics[width=.65\linewidth]{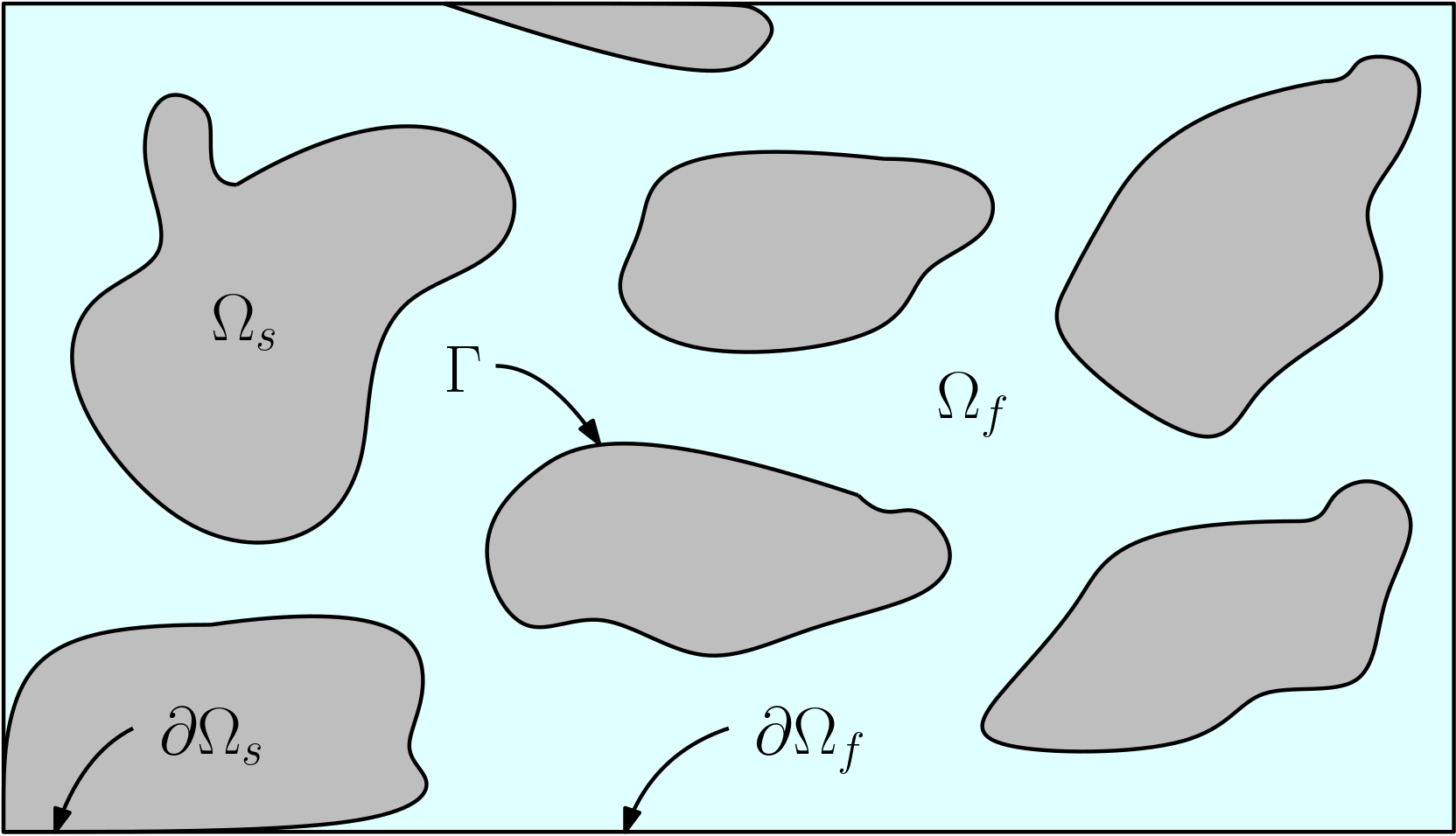}
        \caption{Graphical representation of the whole domain $\Omega = \Omega_{\text{s}} \cup \Omega_{\text{f}}$ considering two sub-domains: $\Omega_{\text{s}}$ a solid region (not necessarily connected), with external boundary $\partial\Omega_{\text{s}}$; and, $\Omega_{\text{f}}$ a fluid region with $\partial\Omega_{\text{f}}$ external boundary such that $\partial \Omega = \partial\Omega_{\text{s}}\cup\partial\Omega_{\text{f}}$. Solid-fluid interface denoted as $\Gamma$.}
        \label{fig::generalSolidFluidGeometry}
    \end{figure}

    Here we discuss the equations that makeup Stokes-Poisson-Nernst-Planck, modelling the advective, diffusive and electrostatic forces of an ionic solution. As many real-world applications of electrochemistry involve interacting solids and fluids we consider a  multi-domain scenario of a whole domain $\Omega$ split, without loss of generality, into two sub-domains $\Omega_{\text{f}}$, a fluid, and $\Omega_{\text{s}}$, a solid, such that $\Omega=\Omega_{\text{f}}\cup\Omega_{\text{s}}$, and with $\Gamma$ being the solid-fluid interface and $\partial\Omega=\partial\Omega_{\text{f}}\cup\partial\Omega_{\text{s}}$ the external boundaries. See \cref{fig::generalSolidFluidGeometry} for a theoretical sketch of the domain $\Omega$. More in general, in some applications and in our computational framework, we have allowed for an arbitrary number of solid and fluid regions, separated by different interfaces. We consider $N$ ionic species,  with concentrations and valencies $c_i$ and $z_i$ respectively, and $i=1,...,N$. To describe the ionic transport we must define equations for the electric field $\mathbf{E}$, ion concentrations $c_i$ within $\Omega$ and fluid velocity $\mathbf{u}$ within $\Omega_{\text{f}}$.

\subsection{Stokes' flow}\label{sec::StokesFlow}

    Consider $\Omega_{\text{f}}$, with velocity profile $\bu(\bx,t)$ governing the advective dynamics of the ions. Assume a negligible Reynolds numbers $\RE$ defined by the fluid density $\rho_{\text{f}}$, characteristic velocity $U$, characteristic length scale $L$ and dynamic viscosity $\mu$:

        \[
                    \RE := \frac{\rho_{\text{f}}L U}{\mu} \ll 1,
        \]
    such that viscous forces within $\Omega_{\text{f}}$ outweigh the inertial. This common assumption in ionic transport \cite{Mark_Schmu_PNP_Porous,Ionic_Transport_Processes,Constantin2019} leads to linear  Stokes flow.  Furthermore, we assume the fluid to be incompressible, i.e., 
        \begin{align}\label{eqn::StokesFlow}
            \mu\nabla^{2} \bu = \bm{\nabla}{p} - \rho_{\text{el}}\mathbf{E},&\quad \bx\in\Omega_{\text{f}},\\
            \bm{\nabla}\cdot \bu =0,& \quad \bx\in\Omega_{\text{f}},\label{eqn::Continuity}
        \end{align}
    where $p=p(\bx,t)$, $\displaystyle\rho_{\text{el}}=\rho_{\text{el}}(\bx,t)=F\sum_{i=1}^{N}z_{i}c_{i}$ and $\mathbf{E}$ are the static fluid pressure, electric charge density and electric field, respectively, and $F$ is Faraday's constant. Compared to the standard Stokes equation, we have the presence of the body force term $\rho_{\text{el}}\mathbf{E}$, describing the Coulomb forces acted on the fluid by the ions\cite{Ionic_Transport_Processes,Probstein2005}. We neglect any magnetic contribution by assuming our ions move slowly such that $\mathbf{E}$ is irrotational, i.e. $\bm{\nabla}\times\mathbf{E}=\mathbf{0}$. This body force term may be set to zero if the fluid is unaffected by $\mathbf{E}$ or the fluid is approximated as electrically-neutral, $\rho_{\text{el}}=0$. This will be further discussed in \S\ref{sec::dimensionalAnalysis}.
    
    The first coupling term, between the  variables $\bu$, $\mathbf{E}$ and $c_i$ appears here, showing one of the significant difficulties of describing electrokinetic flows. 
    The coupling terms (particularly if non-linear) often add significant numerical difficulties. First of all, they make the velocity field time-dependent. Although by neglecting the time derivative we assume instantaneous relaxation to an equilibrium, and thereby a steady solution, by involving $c_i(\bx,t)$ the relaxation becomes tied to the time scale of Nernst-Planck, which is order magnitudes different. For segregated approaches, the disparity of relaxation time scales between Stokes, Poisson and Nernst-Planck -- leading to mixed parabolic-elliptic systems -- can pose severe instability problems or slow convergence of the coupled system.

\subsection{Poisson's equation}

    To model the electric field $\mathbf{E}$, neglecting magnetic forces, We may then write $\mathbf{E}=-\bm{\nabla}\phi$ and focus on the electric potential $\phi=\phi(\bx,t)$. Assume for each sub-domain their respective electric permittivity $\varepsilon$ is spatially constant. From Maxwell's equations, we obtain Poisson's electrostatic equation  denoting variations in $\phi$ by changes in the charge density $\rho_{\text{el}}$,
        \begin{equation}\label{eqn::Poissons}
            \lapl \phi = -\frac{\rho_{\text{el}}}{\varepsilon} = -\frac{F}{\varepsilon}\sum_{i=1}^N z_ic_i, \qquad \bx\in\Omega_{\{\text{f,s}\}}.
        \end{equation}
    %
    %
where $\varepsilon=\varepsilon_{\text{s}}$ in the solid and $\varepsilon=\varepsilon_{\text{f}}$ in the fluid.    Again we have direct coupling between our variables, here between $c_i$ and $\phi$, although this time the former appears linearly in the source term. Like for Stokes flow, this equation depends on time only through the source/coupling terms, in particular the time-dependent ionic concentration $c_i(\bx,t)$.

For a conductive material (e.g., solid electrode), assuming $\varepsilon\to\infty$ and Ohm's law (linear relation between the potential gradient and the current density), the governing equation for the potential becomes:
    \begin{equation}\label{eqn::Ohm}
        \bm{\nabla}\cdot\sigma\bm{\nabla}\phi=0,
    \end{equation}
    where $\sigma$ is the conductivity of the material.

\subsection{Ionic transport}

    Assume the ionic fluid in $\Omega_{\text{f}}$ is sufficiently dilute to ignore ion-ion interactions and diffusion is isotropic. Under these assumptions we may use the Nernst-Planck flux \cite{Ionic_Transport_Processes,Planck1890,Maex2017,Gagneux2016} as
        \begin{equation}\label{eqn:Nernst_flux}
            \flux_{i}=\flux(c_i,\phi)=
            \begin{cases}
                -D_{i\text{,s}}\left(\bm{\nabla}c_{i} \right) & \qquad \bx\in\Omega_{\text{s}}, \\
    			-D_{i,\text{f}}\left(\bm{\nabla}c_{i} + \frac{Fz_{i}}{RT}c_{i}\bm{\nabla}\phi\right) + \bu c_i & \qquad \bx\in\Omega_{\text{f}},
            \end{cases}
        \end{equation}
    where  denote $R$, $T$, $D_{i,\text{f}}$ and $D_{i,\text{s}}$ are, respectively, the ideal gas constant, absolute temperature and diffusion coefficient of species $i$ in the fluid and in the solid.    
    Taking \cref{eqn:Nernst_flux} in conjunction with the continuity equation for mass conservation we arrive at the set of equations modelling transport of $c_i$,
        \begin{equation}\label{eqn:PNP_mass_contin}
            \pardiff{c_i}{t} + \bm{\nabla}\cdot\flux_i =0, \qquad \bx\in\Omega_{\text{s,f}}.
        \end{equation}
    \\
    As mentioned, we assume a dilute solution to ignore ion-ion interactions. This may not hold true in some cases. Alternatively the Stefan-Maxwell equations \cite{Compare_SM_PNP,Ionic_Transport_Processes,Krishna1997} may be used in lieu of Nernst-Planck. In short, Stefan-Maxwell balances the driving forces exerted on a species with the frictional forces between species. This introduces cross diffusivities  $\mathcal{D}_{ij}$ describing the drag between species $i$ and $j$ due to these frictional forces \cite{Ionic_Transport_Processes}. The added complexity from such explicit coupling between species is difficult to describe, as the exact description of $\mathcal{D}_{ij}$ is hard to determine \cite{Compare_SM_PNP}. Even in situations where the ionic fluid is bordering on dilute, it is common for Nernst-Planck to still be used \cite{Charge_Trans_Modelling_Review_Richard} due to its simplicity.
    
    To close the system of equations listed here, we must include conditions along the boundaries of $\Omega_{\text{f}}$ and $\Omega_{\text{s}}$. We briefly discuss some common choices below and their physical representation. We also require conditions at the interface $\Gamma$ to describe the interaction between sub-domains.  This is discussed further below for $c_i$ and $\phi$ when considering reactions and (non-)conductive interfaces.

\subsection{Boundary conditions}\label{sec::BoundaryConditions}
    
Below we list some relevant conditions for different physical situations often seen in electrokinetic problems:
        \begin{itemize}
            \item \textit{Solid walls:} Either a conductive or insulating wall. The no-penetration ($\bu\cdot\mathbf{n}=0$, and zero normal stresses), where $\mathbf{n}$ is the outward unit normal to the boundary, or no-slip condition ($\bu = \mathbf{0}$) may be used. Ions may not pass through the wall so their fluxes are zero ($\flux_i\cdot\mathbf{n}=0$). If conductive we impose a fixed normal current density $I_\text{w}$ along the wall ($\displaystyle\mathbf{i}\cdot\mathbf{n} = F\sum_i z_i\flux_i\cdot\mathbf{n} = I_{\text{w}}$ for a liquid electrolyte) or fixed potential $\phi_\text{w}$. If insulating (non-conductive) we may simply set $I_\text{w}=0$.
            
            \item \textit{Permeable membrane and reactive boundaries:} For a membrane, we can fix the flux of each species ($\flux_i\cdot\mathbf{n}=j_{\text{m}}$), non-zero for species capable of passing through, zero otherwise. We will discuss the case of reactive boundaries later in \S\ref{sec::Interface_Conditions}.
            
            \item \textit{Inlet/Outlet:} Often used to model an in/out-flow of fluid and/or ions. For a momentum driven inlet, pressure or flow velocity may be fixed ($\bu=\bu_{\text{in}}$, $p=p_{\text{in/out}}$). For ions we may apply fixed fluxes ($\flux_i\cdot\mathbf{n}=j_{\text{in/out}}$), resulting in a fixed electric current ($\displaystyle\mathbf{i}\cdot\mathbf{n} = I_{\text{in/out}}$). More complicated is to determine suitable boundary conditions for the potential. Typically, either a Dirichlet $\phi=\phi_{\text{in/out}}$ or a Neumann condition consistent with the ionic current can be imposed. The latter imposes the total electric current to be either equal to zero or a fixed value.
            
            \item \textit{Periodicity:} when dealing with large (quasi-)periodic (or homogeneous) structures, it is often impossible to solve for the entire domain of interest. In these cases, smaller representative unit cells can be solved with (quasi-)periodic external boundary conditions \cite{boccardo2020computational} are imposed. In these cases, additional driving forces (as a bulk source term or a modification of the periodic BC) need to be added.
        \end{itemize}


\section{Dimensional analysis}\label{sec::dimensionalAnalysis}

\begin{figure}[htbp]
			\centering
			\includegraphics[width=.5\linewidth]{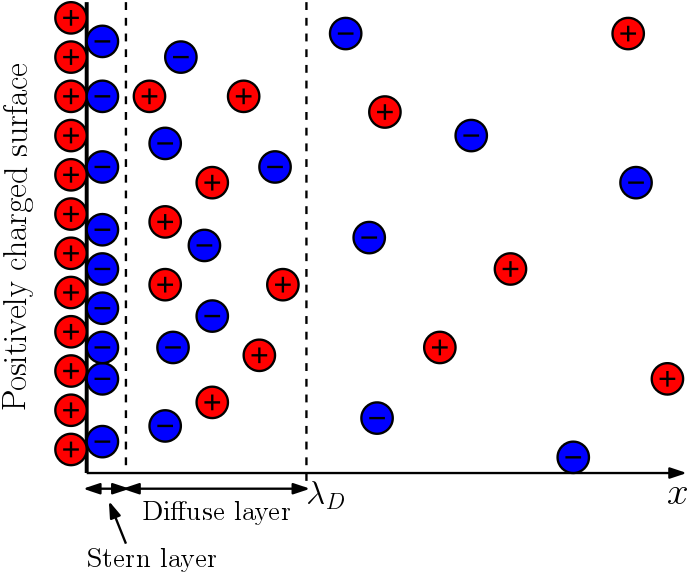}
			\caption{{A graphical representation of the electric double layer along a positively charged surface, comprised of a thin layer of adsorbed negative ions (Stern layer) and loosely connected positive and negative ions (Diffuse layer).}}
			\label{fig::EDL}
		\end{figure}

    When dealing with systems of coupled transport equations it is useful to perform a dimensional analysis to better understand the relationships between the different transport phenomena and identify the limiting regimes and possible approximations. We denote dimensionless variables with a hat symbol, e.g. $\hat{x}$, and reference values with a bar unless stated otherwise. We use the reference values $L$, $U$, $\bar{\phi}$ and $\bar{c}_i$ for the length scale, velocity, electric potential and concentrations respectively. For pressure, we take $\bar{p}=\frac{\mu U}{L}$ as this the appropriate form when under Stokes flow. For time, we take the diffusive timescale $\bar{t}=\frac{L^{2}}{D_{i}}$. With these choices and defining $\bar{C}=\displaystyle\frac{1}{2}\sum_{i=1}^{N}z_{i}^{2}\bar{c_i}$ as the reference total ionic strength, we obtain the dimensionless variables:
		\begin{equation}
			\hat{x} = \frac{x}{L}, \quad \hat{t} = \frac{D_it}{L^2}, \quad \hat{ \bu } = \frac{ \bu }{U}, \quad \hat{\phi} = \frac{\phi}{\bar{\phi}}, \quad \hat{c_{i}} = \frac{c_{i}}{\bar{c_i}}, \quad \hat{p} = \frac{pL}{\mu U}.
		\end{equation}
	\\
	Substituting these dimensionless variables into \cref{eqn::StokesFlow}, \cref{eqn::Poissons} and \cref{eqn:PNP_mass_contin} the dimensionless system of equations are
		\begin{align}
			\frac{D_i\bar{c_i}}{L^2}\pardiff{\hat{c_{i}}}{\hat{t}} &+ \hat{\bm{\nabla}}\cdot\left(-D_{i}\left(\frac{\bar{c_i}}{L^2}\hat{\bm{\nabla}} \hat{c_{i}} + \frac{F\bar{\phi}\bar{c_i}}{L^{2}RT}\hat{c}_{i}z_{i}\hat{\bm{\nabla}}\hat{\phi}\right) + \frac{U\bar{c_i}}{L}\hat{c_{i}}\hat{ \bu }\right) =0, \label{eqn:Non_Dim_Ion_Mass_Cons}\\
			\frac{\bar{\phi}}{L^2}\hat{\nabla}^{2}\hat{\phi} &= -\frac{F}{\varepsilon}\sum_{i=1}^{N}z_i\hat{c_i}\bar{c_i}, \label{eqn:Non_Dim_Poissons} \\
			\frac{\mu U}{L^{2}}\hat{\nabla}^{2}\hat{ \bu } &= \frac{\mu U}{L^{2}}\hat{\bm{\nabla}}\hat{p} + \frac{F\bar{\phi}}{L}\hat{\bm{\nabla}}\hat{\phi}\sum_{i=1}^{N}z_{i}\hat{c}_{i}\bar{c_i}. \label{eqn:Non_Dim_Stokes}
		\end{align}
    The underlying dimensionless numbers may be found by dividing all other terms by the reference values of one term. For  \cref{eqn:Non_Dim_Ion_Mass_Cons,eqn:Non_Dim_Poissons,eqn:Non_Dim_Stokes} we divide by $\frac{D_{i}\bar{c_i}}{L^{2}}$, $\frac{2\bar{C}F}{\varepsilon}$ and $\frac{2\bar{C}\mu U}{2\bar{C}L^{2}}$ respectively, where recall $\bar{C}$ is the reference ionic strength needed to resolve the issue of not being able to factor out the reference concentrations $\bar{c}_i$. As such, the dimensionless equations become:
        \begin{align}
			&\pardiff{\hat{c_i}}{\hat{t}} + \hat{\bm{\nabla}}\cdot\left(-\left(\hat{\bm{\nabla}}\hat{c_{i}} + \N\hat{c_{i}}z_i\hat{\bm{\nabla}}\hat{\phi}\right)+\PE_i \hat{ \bu }\hat{c_{i}}\right) = 0,\label{eqn::NonDimNernstPlanck}\\
			&\L^{2}\hat{\nabla}^{2}\hat{\phi} = -\frac{\sum_{i=1}^{N}\hat{c_i}\bar{c_{i}}z_{i}}{2\bar{C}},\label{eqn::NonDimPoi}\\
			&\hat{\nabla}^{2}\hat{ \bu } = \hat{\bm{\nabla}}\hat{p} + \frac{\P}{2\bar{C}}\sum_{i=1}^{N}\hat{c_i}\bar{c_{i}}z_{i}\hat{\bm{\nabla}}\hat{\phi},\label{eqn::NonDimStokes}
		\end{align}
	\\
    where we have defined four dimensionless numbers $\N$, $\PE$, $\L$ and $\P$. Note that whilst \cref{eqn:Non_Dim_Ion_Mass_Cons} had four terms, we only arrive at two numbers due to our choice of reference time. The same can be said with \cref{eqn::NonDimStokes} and chosen reference pressure. This results in the following four dimensionless numbers for the three equations \cref{eqn::DimensionlessNums} as:

        \begin{equation}\label{eqn::DimensionlessNums}
            \N = \frac{F\bar{\phi}}{RT} = \frac{e\bar{\phi}}{k_{\text{B}} T}, \quad \PE = \frac{UL}{D_i}, \quad \L^2 = \frac{\varepsilon k_{\text{B}} T}{2N_{\text{A}} e^2 L^2 \bar{C}}, \quad \P = \frac{2FL\bar{\phi}\bar{C}}{\mu U},
        \end{equation}
    written in terms of the aforementioned reference values, where we denote $e$, $k_{\text{B}}$ and $N_{\text{A}}$ to be the elementary charge, Boltzmann's constant and Avogadro's number respectively. $\N$ represents the ratio of electrostatic over diffusive forces, with $\N \ll 1$ indicating diffusion is dominant. We also obtain the P\'eclet number $\PE$, i.e., the ratio of advective over diffusive phenomena, and $\L = \frac{\lambda_{\text{D}}}{L}$ the dimensionless Debye length, with the dimensional form $\lambda_{\text{D}}$ defined as
        \begin{equation}
            \lambda_{\text{D}} = \sqrt{\frac{\varepsilon k_{\text{B}}T}{2N_{A}e^2\bar{C}}}\,,
        \end{equation}
     approximating the distance at which a charge's electrostatic effect persists. Typically $\L\ll 1$, stating how the Debye length is much smaller than the reference length $L$. Finally, we have $\P$ denoting the ratio of viscous and electric forces upon our fluid. $\lambda_{\text{D}}$ is also the approximate width of a common electrokinetic phenomenon known as the electric double layer (EDL), see \cref{fig::EDL}, that forms on boundaries. The EDL consists first of ions adsorbed at the boundary, known as the Stern layer, and another of free ions moved by electrical attraction and diffusive motion by the Stern layer, deemed the diffuse layer.
	
    For sufficiently thin EDLs there is a common model reduction known as the electro-neutrality assumption, often employed in the development of macroscopic models \cite{Charge_Trans_Modelling_Review_Richard,Mark_Schmu_PNP_Porous,Newman_Revisit}. This reduction can in fact be determined through dimensional analysis and will be briefly discussed next.
    
\subsection{Asymptotics and electro-neutrality}

    The electro-neutrality approximation states that for a sufficiently dilute electrolyte, all charges of ionic species within the solution roughly cancel each other out, leaving the solution electrically neutral. Electro-neutrality for a solution containing $N$ species is defined as
        \begin{equation}\label{eqn::electroNeutrality}
            \sum_{i=1}^{N}z_i c_i = 0,
        \end{equation}
    and is often used as model reduction when modelling ionic flows. Where \cref{eqn::electroNeutrality} becomes invalid however is in the thin charged double layers, or EDL, mentioned above and often seen in real-world settings.

    To get an understanding of where electro-neutrality comes from and how it relates to EDL formation, we perform an asymptotic analysis. Without any lack of generality, we limit here to a one-dimensional electrolyte and a binary electrolyte. In the numerical framework described above an arbitrary number of ionic species in three dimensions can be considered. Let $x \in \Xi$ where $\Xi = [0,1]$ is a binary electrolyte solution with ionic concentrations $c_1$ and $c_2$ of opposing valencies. For now, we omit any boundary conditions, only requiring these to result in a boundary layer formation near $x=0$. To arrive at \cref{eqn::electroNeutrality} we start by considering the asymptotics of the outer (bulk) layer away from $x=0$. For simplicity, we only consider the leading order terms.
    
    \paragraph{Outer (bulk) layer:}
    To arrive at an asymptotic leading order solution in the bulk layer of $\Xi$ we take the following expansions of $\hat{u}$, $\hat{\phi}$, $\hat{p}$, $\hat{c_1}$ and $\hat{c_2}$. These denote the fluid velocity, electric potential, pressure and ion concentrations respectively, all dimensionless. Expansions are taken in powers of $\epsilon= \mathcal{L}^2$, the squared dimensionless Debye screening length since $\epsilon\ll 1$, or,  $\lambda_{\text{D}} \ll L$:
        \begin{align}
            \hat{\phi} &= \phi^{(0)} + \epsilon\phi^{(1)} + \mathcal{O}(\epsilon^2), \label{eqn::AsymptoticExpansionNonDimPhi}\\
            \hat{c_i} &= c_i^{(0)} + \epsilon c_i^{(1)} + \mathcal{O}(\epsilon^2), \quad i=1,2, \label{eqn::AsymptoticExpansionNonDimCi}\\
            \hat{u} &= u^{(0)} + \epsilon u^{(1)} + \mathcal{O}(\epsilon^2), \label{eqn::AsymptoticExpansionNonDimU} \\
            \hat{p} &= p^{(0)} + \epsilon p^{(1)} + \mathcal{O}(\epsilon^2). \label{eqn::AsymptoticExpansionNonDimP}
        \end{align}
    \\
    Substituting these expansions into the transport \cref{eqn::NonDimNernstPlanck,eqn::NonDimPoi,eqn::NonDimStokes} for one-dimension, alongside the incompressibility condition, we arrive at
        \begin{align}
            \begin{split}
             &\mathcal{O}(\epsilon^2) =\pardiff{}{t}\left(c_i^{(0)} + \epsilon c_i^{(1)}\right) + \pardiff{}{x}\bigg(\PE_i \left[u^{(0)} + \epsilon u^{(1)}\right]\left[c_i^{(0)} + \epsilon c_i^{(1)}\right] - \pardiff{}{x}\left(c_i^{(0)} + \epsilon c_i^{(1)}\right) \\
            &\hspace{4em}- \N z_i\left(c_i^{(0)} +\epsilon c_i^{(1)}\right)\pardiff{}{x}\left(\phi^{(0)} + \epsilon \phi^{(1)}\right)\bigg), \quad i=1,2,
            \end{split}\\
            &\mathcal{O}(\epsilon^2)=\epsilon \secpardiff{}{x}\left(\phi^{(0)} + \epsilon \phi^{(1)}\right) +\frac{1}{2\bar{C}}\sum_{i=1}^2 z_i\bar{c_i}\left(c_i^{(0)} + \epsilon c_i^{(1)}\right), \\
            \begin{split}
            &\mathcal{O}(\epsilon^2)=-\secpardiff{}{x}\left(u^{(0)} + \epsilon u^{(1)}\right) + \pardiff{}{x}\left(p^{(0)} + \epsilon p^{(1)}\right) \\
            &\hspace{4em}+ \frac{\P}{2\bar{C}}\pardiff{}{x}\left(\phi^{(0)} + \epsilon \phi^{(1)}\right)\sum_{i=1}^2 z_i\bar{c_i}\left(c_i^{(0)} + \epsilon c_i^{(1)}\right),
            \end{split}\\
            &\mathcal{O}(\epsilon^2)=\pardiff{}{x}\left(u^{(0)} +\epsilon u^{(1)}\right).
        \end{align}
    \\
    By considering only the leading order terms i.e., terms of $\mathcal{O}(1)$, these equations reduce to,
        \begin{align}
            &\pardiff{}{t}c_i^{(0)} + \pardiff{}{x}\left(\PE_i u^{(0)}c_i^{(0)} - \pardiff{}{x}c_i^{(0)} -\N z_ic_i^{(0)}\pardiff{}{x}\phi^{(0)} \right) = 0, \quad i=1,2, \label{eqn::oneDAsymptoticsNernstPlanck}\\
            & 0= z_1\bar{c_1}c_1^{(0)} + z_2\bar{c_2}c_2^{(0)}, \label{eqn::oneDAsymptoticsElectroNeutralApprox}\\
            &\secpardiff{}{x}u^{(0)} = \pardiff{}{x}p^{(0)}, \\
            &\pardiff{}{x}u^{(0)} =0.
        \end{align}
    Note how \cref{eqn::oneDAsymptoticsElectroNeutralApprox} is the electro-neutrality approximation mentioned before, for a binary solution, and a direct consequence of the asymptotics. This implies said electro-neutrality is only accurate up to leading order. The electric body force of Stokes vanishes as a consequence of \cref{eqn::oneDAsymptoticsElectroNeutralApprox}. To write an equation for the leading order potential $\phi^{(0)}$ we can multiply \cref{eqn::oneDAsymptoticsNernstPlanck} by their respective valencies $z_i$ and reference values $\bar{c_i}$, sum over $i$, and utilise \cref{eqn::oneDAsymptoticsElectroNeutralApprox}, resulting in the equation,
        \begin{equation}\label{eq:pot_en}
            \pardiff{}{x}\left(u^{(0)}\sum_i \PE_i z_i\bar{c_i} c_i^{(0)} - \sum_i z_i\bar{c_i}\pardiff{}{x}c_i^{(0)} -\N\sum_i z_i^2 \bar{c_i}c_i^{(0)} \pardiff{}{x}\phi^{(0)}\right) = 0,
        \end{equation}
which represents a steady equation for $\phi^{(0)}$. This can be interpreted as a modified Ohm's law with new  conductivity accounting for all contributions to the total current: advective, diffusive and electrical, from left to right.
    
    Now that we have constructed the asymptotic solution, up to $\mathcal{O}(1)$, for the outer (bulk) layer of $\Xi$, we move on to determine the inner (boundary) layer asymptotics. As mentioned at the start, we state said inner layer formation is near $x=0$ of $\Xi$. Much like the outer layer, we make no case of boundary conditions, only that a boundary layer near $x=0$ forms as a result of them.
    
    \paragraph{Inner solution:}
    
    To construct a solution for the inner layer near $x=0$ of $\Xi$ we define the variable $y$ to span the inner layer and be a 'fast' variable counterpart to $x$, changing more rapidly,
        \begin{equation}
            y = \frac{x}{\sqrt{\epsilon}}, \quad y \in [0,\infty)\,.
        \end{equation}
    By the definition of $x\in\Xi$, $y$ therefore has domain $[0,\infty)$, as $y=0$ for $x=0$ and $y\rightarrow \infty$ for $x\rightarrow 1$, since $\epsilon \ll 1$. Such change of variable gives derivatives via chain rule as:
        \begin{equation}
            \pardiff{}{x} = \frac{1}{\sqrt{\epsilon}}\pardiff{}{y}, \quad \secpardiff{}{x} = \frac{1}{\epsilon}\secpardiff{}{y}.
        \end{equation}
    Substitution of the above derivatives into the original transport equations \cref{eqn::NonDimNernstPlanck,eqn::NonDimPoi,eqn::NonDimStokes} results in,
        \begin{align}
            &\pardiff{\hat{c_i}}{\hat{t}} + \frac{1}{\sqrt{\epsilon}}\pardiff{}{y}\left(\PE_i \hat{u} \hat{c_i} -\frac{1}{\sqrt{\epsilon}}\pardiff{\hat{c_i}}{y} - \frac{1}{\sqrt{\epsilon}}\N z_i \hat{c_i}\pardiff{\hat{\phi}}{y}\right) = 0, \\
            &\secpardiff{\hat{\phi}}{y} = -\frac{1}{2\bar{C}}\sum_{i=1}^2 z_i\bar{c_i}\hat{c_i}, \\
            &\frac{1}{\epsilon}\secpardiff{\hat{u}}{y} = \frac{1}{\sqrt{\epsilon}}\pardiff{\hat{p}}{y} + \frac{1}{\sqrt{\epsilon}}\frac{\P}{2\bar{C}}\pardiff{\hat{\phi}}{y}\sum_{i=1}^2 z_i\bar{c_i}\hat{c_i}, \\
            &\pardiff{\hat{u}}{y} = 0,
        \end{align}
    where all variables $\hat{u}$, $\hat{p}$, $\hat{\phi}$ and $\hat{c_i}$ are in terms of y e.g., $\hat{c_i} = \hat{c_i}(y,t)$. Just like the outer layer, we substitute the asymptotic expansions \cref{eqn::AsymptoticExpansionNonDimPhi,eqn::AsymptoticExpansionNonDimCi,eqn::AsymptoticExpansionNonDimU,eqn::AsymptoticExpansionNonDimP} into the above equations to obtain,
        \begin{align}
            \begin{split}
                \pardiff{}{t}\left(c_i^{(0)} + \epsilon c_i^{(1)}\right) 
                &+ \frac{1}{\sqrt{\epsilon}}\pardiff{}{y}\bigg(\PE_i\left[u^{(0)} + \epsilon u^{(1)}\right]\left[c_i^{(0)} + \epsilon c_i^{(1)}\right] -\frac{1}{\sqrt{\epsilon}}\pardiff{}{y}\left(c_i^{(0)} + \epsilon c_i^{(1)}\right) \\
                &\quad - \N z_i \frac{1}{\sqrt{\epsilon}}\left(c_i^{(0)} + \epsilon c_i^{(1)}\right)\pardiff{}{y}\left(\phi^{(0)} + \epsilon \phi^{(1)}\right)\bigg) = \mathcal{O}(\epsilon^2), \quad i=1,2,
            \end{split}\\
            &\secpardiff{}{y}\left(\phi^{(0)} + \epsilon \phi^{(1)}\right) = -\frac{1}{2\bar{C}}\sum_{i=1}^2 z_i\bar{c_i}\left(c_i^{(0)} + \epsilon c_i^{(1)}\right) + \mathcal{O}(\epsilon^2), \\
            \begin{split}
                \frac{1}{\epsilon}\secpardiff{}{y}\left(u^{(0)} + \epsilon u^{(1)}\right) 
                &= \frac{1}{\sqrt{\epsilon}}\pardiff{}{y}\left(p^{(0)} + \epsilon p^{(1)}\right) \\
                & \quad + \frac{1}{\sqrt{\epsilon}}\frac{\P}{2\bar{C}}\pardiff{}{y}\left(\phi^{(0)} + \epsilon\phi^{(1)}\right)\sum_{i=1}^2 z_i\bar{c_i}\left(c_i^{(0)} + \epsilon c_i^{(1)}\right) + \mathcal{O}(\epsilon^2),
            \end{split}\label{eqn::oneDInnerAsymptoticsStokesNonDim}\\
            &\pardiff{}{y}\left(u^{(0)} + \epsilon u^{(1)}\right) = \mathcal{O}(\epsilon^2).
        \end{align}
    Considering only the leading order terms of the four equations above we arrive at the following leading order set of equations for the inner layer:
        \begin{align}
            &\pardiff{}{y}\left(-\pardiff{}{y}c_i^{(0)} -\N z_ic_i^{(0)}\pardiff{}{y}\phi^{(0)}\right) =0, \quad i=1,2, \\
            &\secpardiff{}{y}\phi^{(0)} = -\frac{1}{2\bar{C}}\sum_{i=1}^2 z_i\bar{c_i}c_i^{(0)}, \\
            &\secpardiff{}{y}u^{(0)} = 0, \\
            &\pardiff{}{y}u^{(0)} = 0.
        \end{align}
    Like the outer layer equations, to close the system we need another equation, this time for $p^{(0)}$. Considering the $\mathcal{O}(\epsilon^{(1/2)})$ terms of \cref{eqn::oneDInnerAsymptoticsStokesNonDim} we can retrieve such an equation as
        \begin{equation}
            \pardiff{}{y}p^{(0)} + \frac{\P}{2\bar{C}}\pardiff{}{y}\phi^{(0)}\sum_{i=1}^2 z_i\bar{c_i}c_i^{(0)} = 0.
        \end{equation}
    As mentioned earlier, when considering situations involving multiple sub-domains we must also have appropriate interface conditions. In many real applications of electrokinetic flows, there are chemical reactions at interfaces between sub-domains, exchanging mass across the ionic species involved. In these cases, neglecting the EDL might lead to significant errors. In the next section, we will consider such reactions occurring on our interface $\Gamma$, formulating appropriate conditions to capture them whilst retaining mass conservation.

\section{Multi-domain formulation and reactions}\label{sec:multidomain}

    Here we formulate conditions to model heterogeneous reactions which are crucial for many electrochemical and electrokinetic problems. We write a general reaction rate that ensures total mass conservation and apply it to form reactive interface conditions, where we consider the scenarios of species that exist in the whole domain (\textit{unrestricted}) or only in a specific sub-domain (\textit{restricted}). We then discuss conditions on $\phi$ when the interface is conductive or non-conductive.
	
\subsection{Reaction model}\label{sec::reactionModel}

	Consider a general elementary reaction transferring mass between reactants $i=J,...,K$  and products $i=K+1,...,M$ with exchanged molar masses $B_{i}$ and valencies $z_i$. We denote $\nu_i>0$ to be the stoichiometric coefficients determining the number of moles of species $i$ is lost or gained and $n$ the number of released ($n<0$) or absorbed ($n>0$) electrons, alongside the electron mass $e^{-}$. The mass balance of the reaction reads:
		\begin{align}\label{def::general_react_eqn}
			\sum_{i=J}^{K}\nu_{i}B_{i} &+ ne^{-} =\sum_{i=K+1}^{M}\nu_{i}B_{i}\,, \\
			n&=\sum_{i=K+1}^{M}z_{i}\nu_{i} - \sum_{i=J}^{K}\nu_{i}z_{i}. \label{de}
		\end{align}
    To formulate \cref{def::general_react_eqn} as conditions on $\Gamma$ we first determine the rate $r_i$ at which each species in our reaction is exchanging mass. More complex reactions involving several intermediate reactions can be decomposed into elementary steps, i.e. a set of elementary reactions. Typically, the overall reaction rate is then given by the rate of the slowest elementary reaction \cite{RaymondChang2000}. One option to determine the rate is the \textit{law of mass action}, also known as the \textit{rate law} \cite{Generalized_Law_Mass_Action,RaymondChang2000}. Given an ideal solution and reaction involving chemical species $[j]$ with stoichiometric coefficients $\nu_{j}$, at dynamic equilibrium
		\begin{equation}\label{law::Law_of_Mass_Action_Chem_Reaction}
			\sum_{\text{reactants}}\nu_{j}[j] \rightleftharpoons \sum_{\text{products}}{\nu}_{j}[j]\,,
		\end{equation}
	where the net reaction rate (density) $r^\prime$ with units [mol/m$^{2}$\,s], is given by:
		\begin{equation}\label{law::Law_of_Mass_Action_Ideal}
			r^{\prime} = k_{\text{f}}\prod_{j}^{}\bigg(c_{j}^{\nu_{j}}\bigg)_{\text{reac}} - k_{r}\prod_{j}^{}\bigg(c_{j}^{\nu_{j}}\bigg)_{\text{prod}}\,.
		\end{equation}
	Here $k_{\text{f}}$, $k_{\text{r}}$ are the rate constants for the forward and reverse reaction and can be empirically modelled by the Arrhenius equation \cite{Probstein2005,RaymondChang2000}. Note that $r^\prime$ is always positive, so we cannot simply take $r_i=r^\prime$ as we must allow $r_i< 0$ for reactant species. Instead, we define the coefficient $\alpha_i$ as
	    \begin{equation}
            \alpha_i = 
            \begin{cases}
                -1 & \text{reactants,} \\
                +1 & \text{products,}
            \end{cases}
        \end{equation}
	and, multiplying by $r^\prime$ and $\nu_i$, we obtain:
	    \begin{equation}
			r_i = \alpha_i \nu_i r^\prime =
			\alpha_i\nu_{i}\left(k_{\text{f}}\prod_{j=J}^{K}\bigg(c_{j}^{\nu_{j}}\bigg)_{\text{reac}} - k_{r}\prod_{j=K+1}^{M}\bigg(c_{j}^{\nu_{j}}\bigg)_{\text{prod}}\right), \qquad i=J,...,M .
		\end{equation}
	For species not involved in the reaction we set $r_i=0$. Assuming a closed reactive system, i.e no trace ion species, we can represent the mass conservation as the balance across reaction rates $r_i$:
		\begin{equation}\label{eqn::reacMassCons}
			\sum_{i=1}^{N}m_i r_i = 0.
		\end{equation}
    We weigh by the  molar masses $m_i$ to convert from moles (a non-conserved quantity between species) to grams. We use the \textit{rate law} as the example here as it is valid for many reactions \cite{RaymondChang2000}. It is also the more general form for the commonly used Butler-Volmer equation for faradaic reactions, which uses the energy dependence of $k_{\text{f}}$ and $k_{\text{r}}$ to give explicit electric potential dependency. It is important to notice that the units of measure of the reaction rate are [mol/s] for bulk reactions and [m\,mol/s] for surface reactions. Therefore, for the case of linear reactions, the reaction constants $k_i$ can be either [1/s] or [m/s].
    
    
\subsection{Interface conditions}\label{sec::Interface_Conditions}
    
    The reaction models above can be applied in the bulk or on an interface. Here we apply them on $\Gamma$ for $\phi$ and $c_i$ given conductive or non-conductive interface and surface reactions respectively. We employ general reaction rates $r_i$ to balance the ionic fluxes $\flux_i$ through $\Gamma$ with mass exchanged by the reaction. For $\phi$ we use conservation of charge to find conditions on the current passing through $\Gamma$ when acting conductively or non-conductively.
\subsubsection{Interface conditions for ions concentration}\label{sec::InterfaceConditionIons}


    Here we outline reactive conditions along $\Gamma$ to model \cref{def::general_react_eqn}. We use a set of general, possibly non-linear, reaction rates $r_i$ for $i=J,\dots,M$ of ion species involved in \cref{def::general_react_eqn}. We make no assumptions on the form of $r_i$, other than making sure \cref{eqn::reacMassCons} holds true. For species not involved, we simply take flux continuity. The result, jump conditions in ionic flux between $\Omega_{\text{s}}$ and $\Omega_{\text{f}}$, equal to their respective reaction rates $r_i$:
	    \begin{equation}\label{eqn:Surf_Reac_LoMA_Ideal}
	        \evalu{\flux_i \cdot \mathbf{n}}{\Gamma_{\text{s}}}{\Gamma_{\text{f}}} =
	        \begin{cases}
	            r_i	& \bx\in\Gamma \text{ for } i=J,...,M, \\ \\
				0 & \bx\in\Gamma \text{ for } i\neq J,...,M.
	        \end{cases}
	    \end{equation}
	 Here we denote $\big[\flux_i \cdot\mathbf{n}\big]_{\Gamma_{\text{f}}}$ to be the normal flux evaluated at $\Gamma$ from $\Omega_{\text{f}}$'s side. Alongside this we assume continuity of concentrations to provide our second condition:
	    \begin{equation}\label{eqn::unresContinuity}
	        c_{i}\big\rvert_{\Gamma_{\text{s}}} = c_{i}\big\rvert_{\Gamma_{\text{f}}}, \quad \bx\in\Gamma.
	    \end{equation}
	So, for all species involved in the reaction \cref{def::general_react_eqn} we have the difference in ionic flux from $\Omega_{\text{s}}$ and $\Omega_{\text{f}}$ to be the rate $r_i$ of the respective species. 	
	In some scenarios, one or more of the ion species may be restricted to reside in a single sub-domain of $\Omega$. This can however be easily modelled by simply modifying the conditions of those restricted species, setting both $\flux_i$ and $c_i$ to be zero in the inaccessible domains.

\subsubsection{Interface conditions for the potential}

    We consider conditions for the electric potential $\phi$ in different situations. Suppose our interface $\Gamma$ is allowing an electric current to flow through (e.g., an electrochemical reaction on the electrolyte side generating an electrical current in the solid electrode). The conservation of charge can be directly linked to Ohm's law (for a solid electrode), i.e., the current through $\Gamma$ is proportional to the sum of changes in ionic fluxes, to provide an interface condition for the potential on the solid side:
        \begin{equation}
            F\sum_{i=1}^{N}z_{i}\evalu{\flux_i \cdot\mathbf{n}}{\Gamma_{\text{s}}}{\Gamma_{\text{f}}} = \left|{\sigma\bm{\nabla}\phi\cdot\mathbf{n}}\right|_{\Gamma_{\text{s}}}, \quad\bx\in\Gamma.
        \end{equation}
    In the electroneutral (outer expansion) approximation, the same condition can be also applied to the fluid side evaluating \cref{eq:pot_en} at $\Gamma$ and this closes the system.  We see therefore that the conditions on the current are closely linked to the previous reactive conditions. If no reaction occurs at $\Gamma$, we simply have continuity of current due to ionic flux continuity. If a reaction is present, then  \cref{eqn:Surf_Reac_LoMA_Ideal} tells us our jump in current is proportional to the sum of reaction rates $r_i$ alongside their respective charge numbers $z_i$. If $\Gamma$ is non-conductive instead, the current is zero, and we can impose a homogeneous Neumann condition for $\phi$ on conductive sides of $\Gamma$ (e.g., a solid electrode or an electroneutral electrolyte).
    This first interface condition need to be replaced by the  continuity of the electric displacement field $\mathbf{D} = \varepsilon\bm{\nabla}\phi$ in case of non-electroneutral electrolytes (inner expansion) or dielectric media.

    The second interface condition for the potential (not needed in the electroneutral approximation) comes from the definition of the potential jump across the interface. This can be related to the surface charge $\Sigma$ from Gauss' law:    
        \begin{equation}
            \evalu{\phi}{\Gamma_{\text{s}}} {\Gamma_{\text{f}}} = \frac{\Sigma}{\varepsilon_\Gamma}, \quad\bx\in\Gamma.
        \end{equation}
    with $\varphi_\Gamma$ the permittivity of the interface. In this case the surface charge needs to be computed from the microscale local interface properties or assumed to be negligible, therefore resulting in the continuity of the potential across the interface.

\section{Numerical implementation}\label{sec:numericalImplementation}

    As we have mentioned our goal is to construct simple, effective numerical solvers to the Stoke-Poisson-Nernst-Planck (SPNP) model at the pore scale. To do so, we employ the computational fluid dynamics (CFD) package \of due to its open-source nature, large active community and robust handling of complex geometries.
    \of is a finite volume library for general unstructured mesh. Coupled equations are solved iteratively in a segregated/splitted approach, solving sequentially each discretized equation. This removes the explicit coupling between equations, as well as the need of linearizing multi-linear terms (terms linear in each variable but where multiple variables appear) at the expense of having internal iterations. These would be unavoidable also if we adopted a monolithic approach due to the non-linearities of the model. 
    
    Two separate solvers have been implemented. The first, \textit{pnpFoam}, models electrokinetic flow of a single ionic fluid, modelled by SPNP. The second, \textit{pnpMultiFoam}, is more generalized, modelling a general set of ionic fluids and solids following SPNP and diffusion respectively. Furthermore, we developed the numerical counterpart, named \textit{mappedChemicalKinetics}, to the reactive conditions in \S\ref{sec::InterfaceConditionIons} for the case of a binary reaction. In this section, we will describe in detail the structure of the solvers and the corresponding boundary conditions.
    
\subsection{Single- and multi-domain solvers}\label{sec::NumericImplementationSolvers}

\of discretizes each equation using a sparse matrix $\mathbf{M}$ with entries relating to the cell centres of the mesh. For example, for the fluid momentum, the discrete form reads: 
        \begin{equation}
            \mathbf{M}[\mathbf{u}] = -\bm{\nabla}p -\rho_{\text{el}}\bm{\nabla}\phi,
        \end{equation}
    where $[\mathbf{u}]$ is the vector of unknown velocities at cell centres and $\mathbf{M}$ is the coefficient matrix scaling the effect of neighbouring cell velocities. Note the right-hand terms are left as source terms and are computed explicitly and the gradient operators can be discretized with different schemes (details about the schemes will be presented in the results section).    
    To solve for the pressure-velocity coupling we employ the PIMPLE algorithm with an extra term due to the electric body force. This works by decomposing $\mathbf{M}$ into its diagonal, $A$, and off-diagonal, $\mathbf{H}$, parts $\mathbf{M}[\mathbf{u}] = A\mathbf{u} - \mathbf{H}$. This leads to the velocity correction equation:
        \begin{equation}
            \mathbf{u} = \frac{\mathbf{H}}{A} - \frac{1}{A}\bm{\nabla}p - \frac{1}{A}\rho_{\text{el}}\bm{\nabla}\phi.
        \end{equation}
    Interpolating $\mathbf{u}$ to the cell faces and taking the dot product with cell face area vectors $\mathbf{S}_{\text{f}}$ leads to the flux $U$ correction equation,
        \begin{equation}
            U = \mathbf{u}_{\text{f}}\cdot\mathbf{S}_{\text{f}} = \left(\frac{\bm{H}}{A}\right)_{\text{f}}\cdot\bm{S}_{\text{f}} - \left(\frac{1}{A}\right)_{\text{f}}\bm{S}_{\text{f}}\cdot\bm{\nabla}^{\perp}_{\text{f}}p^{n+1} - \left(\frac{\rho_{\text{el}}\bm{\nabla}\phi}{A}\right)_{\text{f}}\cdot\bm{S}_{\text{f}},
        \end{equation}
    where subscript $\text{f}$ denotes values at cell faces. Discretising the incompressibility condition $\bm{\nabla}\cdot\mathbf{u}=0$ gives $\bm{\nabla}\cdot U =0$ which when applied to the flux correction equation forms the pressure correction equation:
        \begin{equation}
            \bm{\nabla}\cdot\left(\left[\frac{\mathbf{H}}{A}\right]_{\text{f}}\bm{\nabla}p \right) = \bm{\nabla}\cdot\left(\left[\frac{\mathbf{H}}{A}\right]_{\text{f}} - \left[\frac{\rho_{\text{el}}\bm{\nabla}\phi}{A}\right]_{\text{f}}\right).
        \end{equation}
\\
This is solved iteratively until convergence. In an external loop, momentum and pressure equations are coupled with the concentration and potential equations. When the mesh is highly skewed, or for complex discretisation schemes with implicit-explicit terms, additional iterations can be added for each single equation.
The pseudo-code algorithm for our single ionic fluid solver \textit{pnpFoam} is presented in Algorithm 1.
    	\begin{algorithm}[htbp]\label{alg:pnp}
    	    \SetAlgoLined
    	    Construct fields: $\mathbf{u}$, $p$, $\phi$, $c_{i}$\;
    	    \While{$t<$ end time}{
        	    \While{!converged}{
        	        Initialize $\mathbf{u}^{n+1}$, $p^{n+1}$ to previous values: $\mathbf{u}^{n+1} = \mathbf{u}^{n}$, $p^{n+1} = p^{n}$\;
        	        Solve momentum equation for $\mathbf{u}^{n+1}$: $-\mu\nabla^{2}\mathbf{u}^{n+1} +\rho_{\text{el}}\bm{\nabla}\phi= -\bm{\nabla}p^{n+1}$\;
        	        \While{Pressure correction}{
        	            Solve pressure correction equation: $\bm{\nabla}\cdot\left[\left(\frac{1}{A}\right)_{\text{f}}\bm{\nabla}p^{n+1}\right] = \bm{\nabla}\cdot\left(\left(\frac{\bm{H}}{A}\right)_{\text{f}} - \left(\frac{\rho_{\text{el}}\bm{\nabla}\phi}{A}\right)_{\text{f}}\right)$\;
        	            Correct flux: $U^{n+1} = \mathbf{u}_{\text{f}}\cdot\bm{S}_{\text{f}} = \left(\frac{\bm{H}}{A}\right)_{\text{f}}\cdot\bm{S}_{\text{f}} - \left(\frac{1}{A}\right)_{\text{f}}\bm{S}_{\text{f}}\cdot\bm{\nabla}^{\perp}_{\text{f}}p^{n+1} - \left(\frac{\rho_{\text{el}}\bm{\nabla}\phi}{A}\right)_{\text{f}}\cdot\bm{S}_{\text{f}}$\;
        	            Correct velocity: $\mathbf{u}^{n+1} = \frac{\bm{H}}{A} - \frac{1}{A}\bm{\nabla}p^{n+1} - \frac{\rho_{\text{el}}\bm{\nabla}\phi}{A}$ 
        	        }
        	        \While{Non-orthogonal correction}{
        	            Solve Poisson's: $-\varepsilon\nabla^{2}\phi = \rho_{\text{el}}$\;
        	        }
        	        \For{All $c_{i}$'s}{
            	        \While{Non-orthogonal correction}{
            	            Solve Nernst-Planck: $\pardiff{c_{i}}{t} + \bm{\nabla}\cdot\left(\mathbf{u}c_{i} -D_{i}\left(\bm{\nabla} c_i + \frac{F}{RT}z_ic_i\bm{\nabla}\phi\right)\right) =0$\;
            	        }
        	        }
        	    }
        	    Increment time: $t\rightarrow t + \Delta t $
    	    }
    	    \caption{\textit{pnpFoam} algorithm}
    	\end{algorithm}
	
	The same method is employed within the algorithm of \textit{pnpMultiFoam}, our solver for ionic transport over a general set of ionic fluids and solids, where each region is solved separately, and an additional outer loop is added to ensure the coupling between regions. The algorithm is presented in Algorithm 2.
	Both solvers here are presented for the case of a non-electroneutral solution. In case electroneutrality is assumed, the algorithm is slightly modified to solve a modified
	potential equation \cref{eq:pot_en} (with ionic conductivity instead of permittivity) and with the last species calculated to ensure electro-neutrality.
	
        \begin{algorithm}[htbp]\label{alg:pnpMulti}
            \SetAlgoLined
            Construct fields: $\mathbf{u}$, $p$, $\phi_{\text{f}}$, $\phi_{\text{s}}$, $c_{i,\text{f}}$, $c_{i,\text{s}}$\;
            \While{$t<$ end time}{
                \For{All fluid regions}{
                    Initialize $\mathbf{u}^{n+1}$, $p^{n+1}$ to previous values: $\mathbf{u}^{n+1} = \mathbf{u}^{n}$, $p^{n+1} = p^{n}$\;
        	        Solve momentum equation for $\mathbf{u}^{n+1}$ prediction : $-\mu\nabla^{2}\mathbf{u}^{n+1} +\rho_{\text{el}}\bm{\nabla}\phi= -\bm{\nabla}p^{n+1}$\;
        	        \While{Pressure correction}{
        	            Solve pressure correction equation: $\bm{\nabla}\cdot\left[\left(\frac{1}{A}\right)_{\text{f}}\bm{\nabla}p^{n+1}\right] = \bm{\nabla}\cdot\left(\frac{\bm{H}}{A}\right)_{\text{f}}$\;
        	            Correct flux: $\mathbf{U}^{n+1} = \mathbf{u}_{\text{f}}\cdot\bm{S}_{\text{f}} = \left(\frac{\bm{H}}{A}\right)_{\text{f}}\cdot\bm{S}_{\text{f}} - \left(\frac{1}{A}\right)_{\text{f}}\bm{S}_{\text{f}}\cdot\bm{\nabla}^{\perp}_{\text{f}}p^{n+1}$\;
        	            Correct velocity: $\mathbf{u}^{n+1} = \frac{\bm{H}}{A} - \frac{1}{A}\bm{\nabla}p^{n+1}$ 
        	        }
        	        \While{Non-orthogonal correction}{
        	            Solve Poisson's: $-\varepsilon\nabla^{2}\phi_{\text{f}} = \rho_{\text{el}}$\;
        	        }
        	        \For{All $c_{i,\textnormal{f}}$'s}{
            	        \While{Non-orthogonal correction}{
            	            Solve Nernst-Planck: $\pardiff{c_{i,\textnormal{f}}}{t} + \bm{\nabla}\cdot\left(\mathbf{u}c_{i,\textnormal{f}} -D_{i}\left(\bm{\nabla} c_{i,\textnormal{f}} + \frac{F}{RT}z_ic_{i,\textnormal{f}}\bm{\nabla}\phi_{\text{f}}\right)\right) =0$\;
            	        }
        	        }
                }
                \For{All solid regions}{
                    \While{Non-orthogonal correction}{
        	            Solve Poisson's: $-\varepsilon\nabla^{2}\phi_{\text{s}} = \rho_{\textnormal{el}}$\;
        	        }
        	        \For{All $c_{i,s}$'s}{
            	        \While{Non-orthogonal correction}{
            	            Solve Diffusion equation: $\pardiff{c_{i,\textnormal{s}}}{t} - D_{i}\nabla^{2}c_{i,\textnormal{s}} =0$\;
            	        }
        	        }
                }
                Increment time: $t\rightarrow t + \Delta t $\;
            }
            \caption{\textit{pnpMultiFoam} algorithm}
        \end{algorithm}
        
\subsection{Boundary and interface conditions}\label{sec::NumericalImplementationBoundary}

    To make use of the object orientation of \of, all conditions are first reformulated as effective Robin conditions. We consider here, as an example, the inhomogeneous Robin BC for a variable $c$ with coefficients $D^\star$, $K^\star$ and $F^\star$ along a boundary $\Gamma$ with normal $\mathbf{n}$:
        \begin{equation}\label{eqn::GeneralRobinEffectiveCondition}
            D^\star\bm{\nabla}_{\mathbf{n}}c\rvert_{\Gamma} = -K^\star c\rvert_{\Gamma} - F^\star\,.
        \end{equation}
    \\
    In \of\, the boundary values $c\rvert_{\Gamma}$ and $\bm{\nabla}_{\mathbf{n}}c\rvert_{\Gamma}$ are approximated using the value of the cell centres with faces along $\Gamma$:
        \begin{align}\label{eqn::OFInterpolationField}
            c\rvert_{\Gamma} &\approx \alpha_{1} c_{c} + \alpha_2, \\
            \bm{\nabla}_{\mathbf{n}}c\rvert_{\Gamma} &\approx \alpha_3 c_{c} + \alpha_4 \approx B(c_{\text{f}} - c_{c})\label{eqn::OFInterpolationGradient}\,.
        \end{align}
    Here $c_{\text{c}}$ and $c_{\text{f}}$ are the values of $c$ at the cell centre and face respectively. The $\alpha$'s are the interpolation weights, and $B$ is the inverse distance between the cell centre and the boundary. We can then rearrange \cref{eqn::GeneralRobinEffectiveCondition} using \cref{eqn::OFInterpolationField,eqn::OFInterpolationGradient} to find the $\alpha$ values that allow us to approximate \cref{eqn::GeneralRobinEffectiveCondition} using the cell centres $c_{\text{c}}$:
        \begin{equation}\label{eqn::GeneralRobinConditionAlphas}
            \alpha_1 = \frac{D^\star B}{D^\star B - K^\star}\:, \alpha_2 = \frac{F^\star}{D^\star B - K^\star}\:, \alpha_3 = \frac{K^\star B}{D^\star B - K^\star}\:, \alpha_4 = \frac{F^\star B}{D^\star B - K^\star}\,.
        \end{equation}
    This forms the basis of the \textit{Robin} BC implemented. Reformulating all other conditions into a form like \cref{eqn::GeneralRobinEffectiveCondition} lets us reuse the same equations in \cref{eqn::GeneralRobinConditionAlphas} to approximate all other conditions.
    
    The interface conditions are implemented as a derived class named \textit{mappedChemicalKinetics}. If we consider here the limiting case of two reacting species $c_{\text{s}}$ and $c_{\text{f}}$, each restricted to their respective sub-domains $\Omega_{\text{s}}$ and $\Omega_{\text{f}}$ respectively. Along the interface $\Gamma$ connecting the two sub-domains is the reaction:
        \begin{equation}
            \ce{\beta_{\text{s}}X_{\text{s}} <-> \beta_{\text{f}}X_{\text{f}}}\,.
        \end{equation}
    The general reactive conditions \cref{eqn:Surf_Reac_LoMA_Ideal} in this case reduces to the following condition along $\Gamma$:
        \begin{equation}\label{eqn::ReactionConditionMathematicalForm}
            \flux_{\text{f}}\rvert_{\Gamma_{\text{f}}}\cdot\mathbf{n} = \flux_{\text{s}}\rvert_{\Gamma_{\text{s}}}\cdot\mathbf{n} = r(c_{\text{f}}\rvert_{\Gamma_{\text{f}}},c_{\text{s}}\rvert_{\Gamma_{\text{s}}}),
        \end{equation}
    where we denote here $\Gamma_{\text{s}}$ to be evaluation at the $\Omega_{\text{s}}$ side of $\Gamma$ and $\mathbf{n}$ the unit normal of $\Gamma$ facing into $\Omega_{\text{f}}$. Note the reaction rate $r$ is kept general and not necessarily linear.
    \Cref{eqn::ReactionConditionMathematicalForm} is solved iteratively using the Newton-Raphson method, linearised about the previous solution $\mathbf{c}^{N} = (c_{\text{s}}^N \hspace{.5em}  c_{\text{f}}^N)^\top$. Here the previous solution could be the solution at the previous time step (if an explicit time stepping is chosen) or at the previous internal iteration (for fully implicit time stepping). This allows us to rewrite \cref{eqn::ReactionConditionMathematicalForm} into two decoupled effective Robin conditions for $c_{\text{s}}^{N+1}$ and $c_{\text{f}}^{N+1}$ that we solve separately. The coupling between the two sides of the interface conditions (and therefore the two domains) is achieved through the internal iterations but, for stiff reactions, additional sub-looping to update both boundary values whenever one of the two domains is solved for.    In algorithm 3 we detail the \textit{mappedChemicalKinetics} pseudo-code to solve the effective Robin conditions. As both resulting effective Robin conditions are solved in the same manner we only write here solving of the condition for $c_{\text{f}}$.

        \begin{algorithm}[htbp]
            \SetAlgoLined
            Set tolerance: $\text{tol}$\;
            Set max Newton iterations: $N_{\text{max}}$\;
            Reset iterator: $N=0$\;
            \While{$||c_{\textnormal{f}}^{N} - c_{\textnormal{f}}^{N+1}||$ > $\textnormal{tol}$ $\cap$ $N\leq N_{\textnormal{max}}$}{
                Store previous iteration: $c_{\textnormal{f}}^{N} = c_{\textnormal{f}}^{N+1}$\;
                Extract the face and internal values of the neighbouring field (solid side)\;
                Update $D^{*}_{\textnormal{f}}$, $K^*_{\textnormal{f}}$ and $F^*_{\textnormal{f}}$: see \cref{bcon:D_f_Eff_GenReacBinary,bcon:K_f_Eff_GenReacBinary,bcon:F_f_Eff_GenReacBinary} \;
                Evaluate iterative effective Robin BC: $-D^*_{\textnormal{f}} \bm{\nabla}_{\mathbf{n}}c_{\textnormal{f}}^{N+1} = -K^*_{\textnormal{f}} c_{\textnormal{f}}^{N+1} - F^*_{\textnormal{f}}$\;
                Apply the boundary condition to obtain the new value: $c_{\textnormal{f}}^{N+1}$\;
                \If{Stiff} {Call \textit{mappedChemicalKinetics} BC on the solid side\;}
                Increment iterator: $N \rightarrow N+1$
            }
            \caption{mappedChemicalKinetics}
        \end{algorithm}

More details about the linearization can be found in \cref{app:int}. Alongside the reactive conditions, we have also implemented a number of simpler conditions, such as the continuity of total fluxes, continuity of value or continuity of derivatives, often used within applications. Just as with the non-linear reactive condition, we rewrite all conditions into effective Robin conditions.

\section{Numerical examples}\label{sec:results}

    Here we present four numerical examples of electrokinetic flows. To verify the accuracy of results we compare with the spectral  \textsc{Matlab}\textsuperscript\textregistered \hspace{.4em}toolbox \textit{Chebfun} \cite{Driscoll2014} with machine-precision accuracy. The first case verifies the accuracy of the flow description given a single ion species using a pressure-driven infinite ion channel similar to \cite{Berg2011}. Next, we verify accuracy when considering multiple ion species. Afterwards, we verify the implemented reactive interface conditions counterpart to \S\ref{sec::InterfaceConditionIons}. The final case displays the capabilities of our solver(s), simulating ionic transport within a randomized solid-fluid porous medium.
    
    To show spatial convergence of \of results we use the following normalized  $L^2$ error point norm. We use here the dummy variables $v$ and $v_{\text{cref}}$ to denote \of and \textit{Chebfun} results respectively as an example:
    
        \begin{equation}
            L^2 \text{-norm error }=\frac{\lVert v - v_{\text{cref}}\rVert_2}{\lVert v_{\text{cref}}\rVert_2}.
        \end{equation}

\subsection{SPNP in an infinite channel}
    
    To verify \textit{pnpFoam} and \textit{pnpMultiFoam} for a single ion species take $\Omega$, of length $L$, with two boundaries $\Gamma_{\text{t}}$ and $\Gamma_{\text{b}}$ denoting the top and bottom channel walls respectively. Take boundaries $\Gamma_{\text{in}}$ and $\Gamma_{\text{out}}$ denoting the inlet and outlet of fluid in $\Omega$. To mathematically describe this channel as infinite we take periodic conditions on $\Gamma_{\text{in}}$ and $\Gamma_{\text{out}}$.
	
		\begin{figure}[htbp]
			\centering
			\includegraphics[width=.55\linewidth]{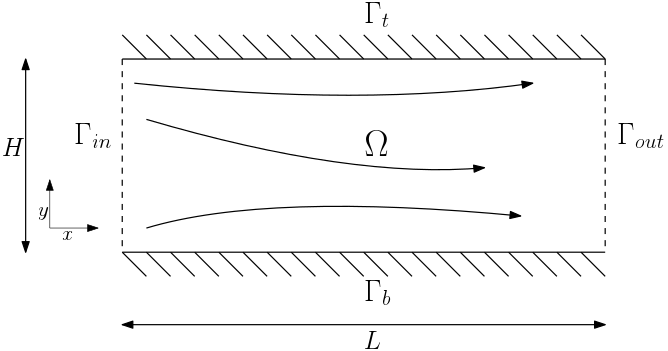}
			\caption{{Graphical representation of the infinite ion channel with channel width $H$.}}
		\end{figure}
    Take a monovalent ionic fluid, i.e a single ion species $c$ with valency $z=1$.  Enforce a fixed external force $\mathbf{J} =\bigg(\begin{smallmatrix}
		-\mu/H^2 \\
		0
	\end{smallmatrix}\bigg)$ to induce transport across the channel length. To allow electro-migration of $c$ between the two channel walls, fix $\phi=0.1$ along $\Gamma_{\text{t}}$ and $\phi=0$ on $\Gamma_{\text{b}}$. For $c$ apply no-flux along $\Gamma_{\text{t}}$ and $\Gamma_{\text{b}}$, denoting $\mathbf{n}$ to be the outward unit normal of either $\Gamma_{\text{t}}$ or $\Gamma_{\text{b}}$. We consider the channel at steady state and, by assuming $L$ is sufficiently large, take all derivatives in $x$ to be zero, i.e. $\pardiff{a}{x}=0$ for $a\in\{ \bu , \phi, c, p\}$. Just like Poiseuille flow we assume a uni-directional velocity, such that $\bu = \left(\begin{smallmatrix} u_1 \\ 0 \end{smallmatrix}\right)$. These assumptions in turn produce a one-dimensional reduce system of equations in $y$.
	
	
	The SPNP system with a single ionic species in an infinite channel is considered with an added pressure-driven driving force $\mathbf{J}$:
        \begin{align}
			\mu\nabla^{2} \bu  - \mathbf{J} - \bm{\nabla}p - Fc\bm{\nabla}\phi = 0\,, 
			 \qquad &\bx\in\Omega ,
			\label{eqn::Stoke}\\
			\bm{\nabla}\cdot \bu  =0\,, 
			 \qquad &\bx\in\Omega ,\label{eqn::InfChannelIncompressible}\\
			\epsilon\nabla^{2}\phi + Fc = 0 \,, 
			 \qquad &\bx\in\Omega ,
			\label{eqn::Poisson}\\
			\bm{\nabla}\cdot\left[-D\left(\bm{\nabla}c  +  \frac{F}{RT}c\bm{\nabla}\phi\right) +  \bu c\right] = 0\,, 
			 \qquad &\bx\in\Omega ,
			\label{eqn::NernstPlanck}\\
			 \bu =\mathbf{0}\,, \quad \phi=0.05\,, 
			  \qquad &\bx\in\Gamma_{\text{t}} \label{bc::Inf_Channel_Top}\\
			 \bu =\mathbf{0}\,,  \quad \phi=0\,, 
			  \qquad &\bx\in\Gamma_{\text{b}} ,\label{bc::Inf_Channel_Bot}\\
			-D\bigg(\bm{\nabla}_{\mathbf{n}}c + \frac{F}{RT}c\bm{\nabla}_{\mathbf{n}}\phi\bigg) + c \bu \cdot\mathbf{n}=0\,,
			 \qquad &\bx\in\Gamma_{\text{b}}\cup\Gamma_{\text{t}} ,\label{bc::InfChannelNoFlux}\\
			{a}\bigg\rvert_{\Gamma_{\text{in}}} = {a}\bigg\rvert_{\Gamma_{\text{out}}}\,,
			\pardiff{a}{x}\bigg\rvert_{\Gamma_{\text{in}}} = \pardiff{a}{x}\bigg\rvert_{\Gamma_{\text{out}}}\,,
			 \qquad & a\in\{c, \bu ,\phi,p\} .\label{bc::Periodic}
		\end{align}
\\
The complete system is solved in \of in a periodic channel, and compared with the equivalent one-dimensional model. By applying the assumptions of uni-directional flow and zero derivatives in $x$ we obtain a one-dimensional reduced set of equations in $y$:
        \begin{align}
			\mu\secdiff{u_1}{y} = -\frac{\mu}{H^2},&\qquad \bx\in\Omega ,\label{eqn::Stokes_Reduced}\\
			0 = \diff{p}{y} + Fc\diff{\phi}{y},& \qquad \bx\in\Omega, \label{eqn::Pressure_Explicit}\\
			\varepsilon\secdiff{\phi}{y} = -Fc,& \qquad \bx\in\Omega, \label{eqn::Poisson_Reduced}\\
			0 = \secdiff{c}{y} + \frac{F}{RT}\left(\diff{c}{y}\diff{\phi}{y} + c\secdiff{\phi}{y}\right),& \qquad \bx\in\Omega, \label{eqn::NernstPlanck_Reduced}\\
			u_1 =0,\, \phi=0.05,& \qquad \bx\in\Gamma_{\text{t}},\\
			u_1 =0,\, \phi=0,& \qquad \bx\in\Gamma_{\text{b}},\\
			\diff{c}{y} + \frac{F}{RT}c\diff{\phi}{y} = 0,& \qquad \bx\in\Gamma_{\text{b}}\cup\Gamma_{\text{t}},
		\end{align}
where we excluded the incompressibility condition as this is trivially satisfied.
	

	To compare results between \textit{Chebfun} and our single ionic fluid solver \textit{pnpFoam} we the following dimensionless numbers, geometrical parameters and transport properties, according to \Cref{tab::InfiniteChannelProperties}.
        
\begin{figure}
            \centering
            \begin{subfigure}{.65\textwidth}
                \includegraphics[width=1.1\linewidth]{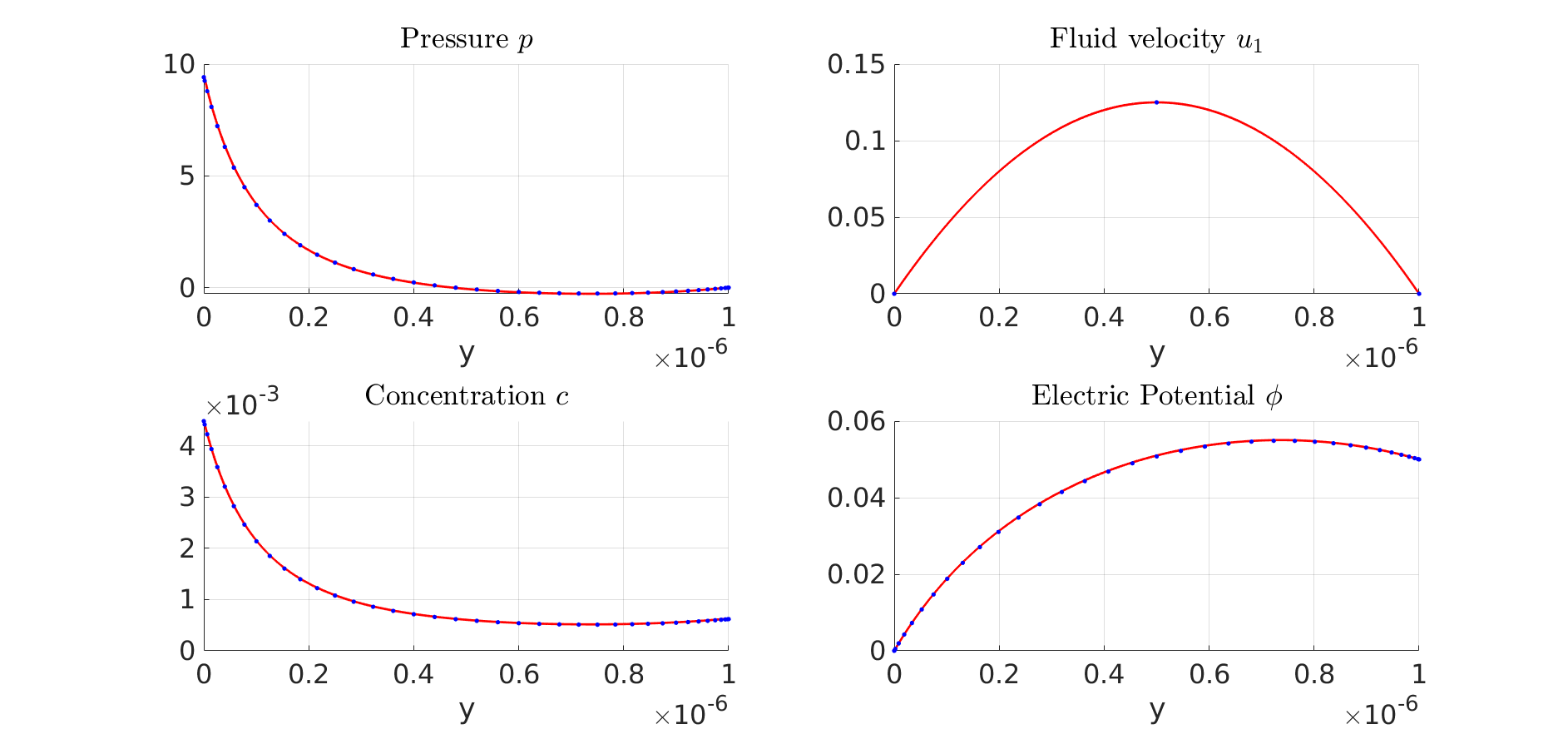}
                \caption{Comparison of \textit{Chebfun} (blue dots) and \textit{pnpFoam} (red line) of the pressure, fluid velocity, ion concentration, and electric potential respectively, along the channel width i.e., in $y$.}
            \end{subfigure}
            \begin{subfigure}{.34\textwidth}
                \includegraphics[width=1.1\linewidth]{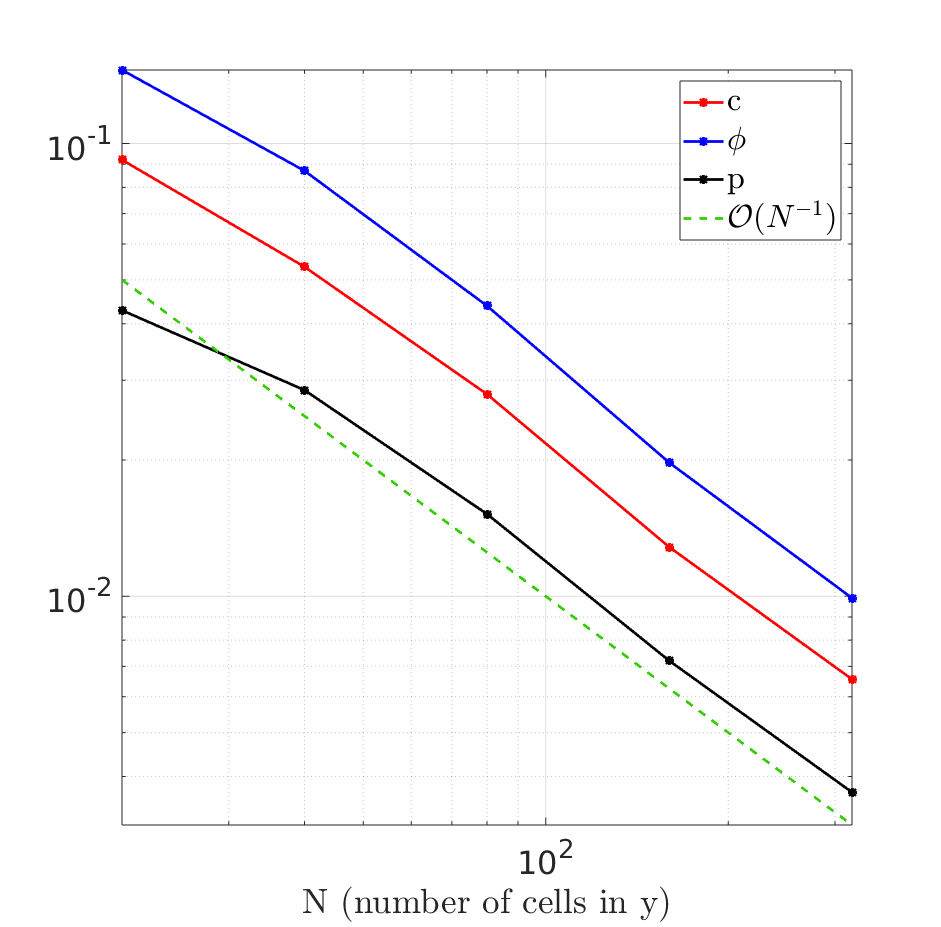}
                \caption{$L^2$ error norm convergence plot of \textit{pnpFoam} results, showing spatial convergence of $\mathcal{O}(N^{-1})$}
                \label{fig::CompareL2ErrorInfChannel}
            \end{subfigure}
        \end{figure}

        \begin{table}[htbp]
    	    \centering
    	    \begin{tabular}{cccccccc}\hline
    	       Symbol & $\mu$ & $H$ & $\varepsilon$ & $\N$ & $\PE$ & $\L^{2}$ & $\P$ \\ \hline
    	       Value & $10$ & $1\times10^{-6}$ & $40\varepsilon_0$ & $1.95$ & $1.25\times10^5$ & $0.09$ & $3.85\times10^{-6}$ \\\hline
    	    \end{tabular}
    	    \caption{Dimensionless numbers, domain and transport properties for the SPNP equation in the infinite channel. We denote $\varepsilon_0 = 8.85\times10^{-12}$ to be the permittivity of vacuum.}
    	    \label{tab::InfiniteChannelProperties}
    	\end{table}
    When running \textit{pnpFoam} we solve it in a steady state and take a large number of PIMPLE corrector iterations to ensure convergence of the coupled system. The mesh used in \of is up to $N=420$ cells ($100$ in $x$, $320$ in $y$). Comparison of results and the $L^{2}$ error norm convergence are shown in \cref{fig::CompareL2ErrorInfChannel}. We observe an accumulation of $c$ along $\Gamma_{\text{b}}$, causing a large pressure gradient to form. The velocity profile follows a parabolic arc, as expected, whilst we see a non-linear profile for $\phi$. Overall we find good agreement, with linear order spatial convergence $\mathcal{O}(N^{-1})$.
    
\subsection{Multi-component ionic fluid}

    To verify \textit{pnpFoam} and \textit{pnpMultiFoam} for a multi-component fluid we take two species $c_1$ and $c_2$ with opposite valencies $z_1=-z_2=1$ over the domain $\Omega\in[0,L]$ where $L=1\times10^{-6}$. We set a zero-flux BC, i.e., $j(c_{i})=0$ at $x=0,L$ for both species and fix $\phi=0.1$ at $x=L$ and $\phi=0$ at $x=0$. We fix here a constant velocity $ \bu =0.001$ for simplicity, and we bypass the solution of the Stokes system (\textit{frozenFlow} flag). Concentrations and potential are initially set to $c_i=1\times 10^{-3}$ and $\phi=0.05(1-\cos{\frac{\pi x}{L}})$.


    We consider here the  dimensionless number values in \Cref{tab::DimNumsCase2Properties}, indicating electrostatic forces dominate. From $\L$ we find the Debye length as approximately $\lambda_{\text{D}} =\SI{308}{\nano\meter}$. Since we are fixing $\bu$ and neglecting Stokes we have $\P=0$.
    
        \begin{table}[htbp]
    	    \centering
    	    \begin{tabular}{ccccc}\hline
    	       Symbol & $\N$ & $\PE$ & $\L$ & $\P$ \\ \hline
    	       Value & $3.868$ & $1\times10^{-3}$ & $9\times10^{-3}$ & $0$ \\\hline
    	    \end{tabular}
    	    \caption{Dimensionless numbers and parameters for the multi-component ionic fluid testcase.}
    	    \label{tab::DimNumsCase2Properties}
    	\end{table}

\begin{figure}[htbp]
            \centering
            \begin{subfigure}{0.34\linewidth}
                \centering
                \includegraphics[width=1.1\linewidth]{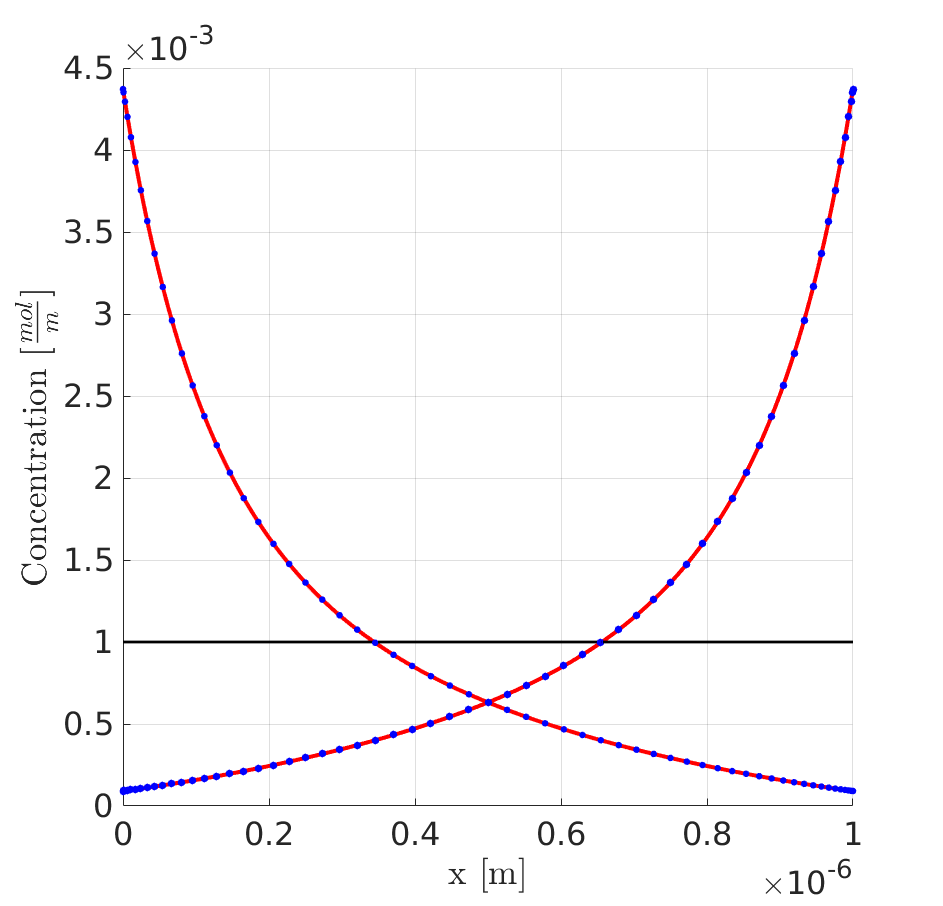}
                \caption{Results of $c_1$ (left curve) and $c_2$ (right) at $t=1\times10^{-5}$s. Blue and red denote \textit{Chebfun} and \textit{pnpFoam} results respectively. Black is the initial condition.}
            \end{subfigure}
            \begin{subfigure}{0.31\linewidth}
                \centering
                \includegraphics[width=1.1\linewidth]{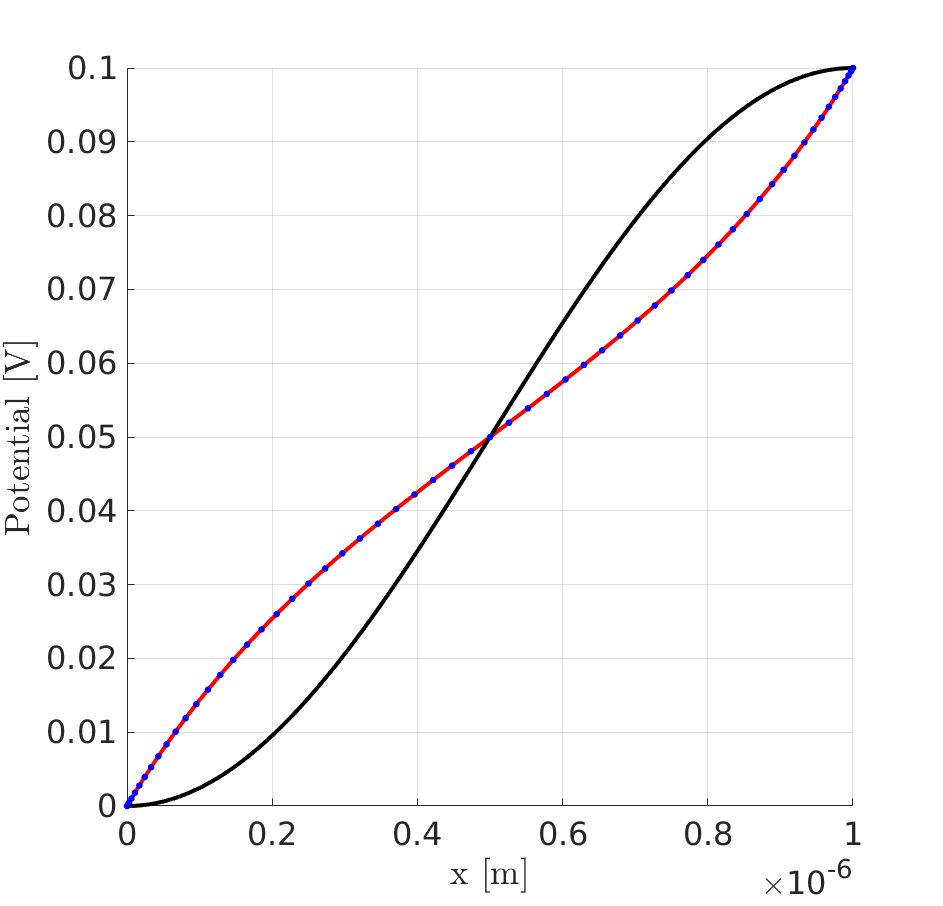}
                \caption{Results of $\phi$ for \textit{Chebfun} (blue) and \textit{pnpFoam} (red) at $t=1\times10^{-5}$s, with initial condition in black.}
            \end{subfigure}
            \begin{subfigure}{.31\textwidth}
                \centering
                \includegraphics[width=1.1\linewidth]{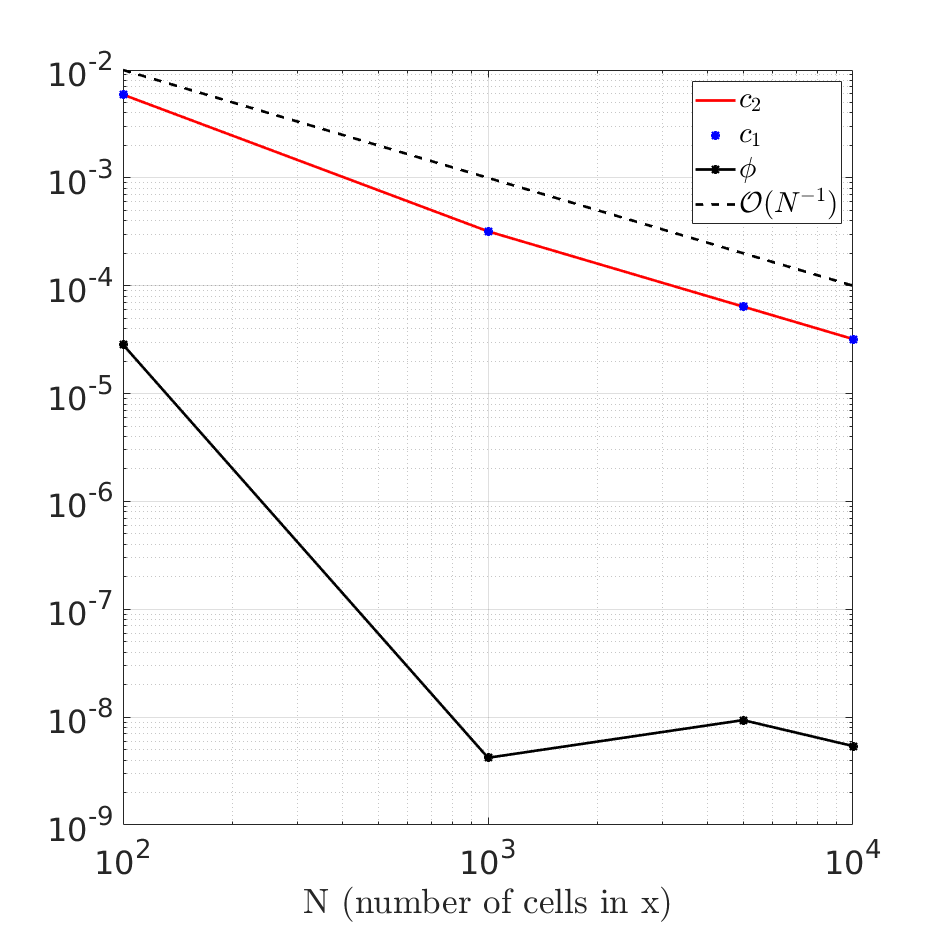}
                \caption{$L^2$ point error norm convergence plot, showing $\mathcal{O}(N^{-1})$ spatial convergence of \textit{pnpFoam}.}
            \end{subfigure}
            \label{fig::CompareL2Normcase2}
        \end{figure}

    When running \textit{pnpFoam} we use a time step $\Delta t=10^{-7}$s, end time $t=1\times10^{-5}$s, third-order implicit time scheme (\textit{backward} keyword) and a mesh of up to $N=1000$ cells. Results and $L^{2}$ norm convergence plots are depicted in \cref{fig::CompareL2Normcase2} where we find good agreement of results and linear spatial convergence of all fields.
	
    We see the ions are transported to the outer walls due to the high electric potential gradient. Most of the ions then accumulate within $\lambda_{\text{D}}$ from the walls forming two overlapping EDLs. Non-linear behaviour between $c_i$ and $\phi$ is observed through the slight shift in $\phi$'s profile due to the clustering of ions at the walls.

\subsection{ Reactive interface}
    Here we verify the accuracy of \textit{mappedChemicalKinetics}, the numerical counterpart to the conditions found in \S\ref{sec::InterfaceConditionIons}. Consider the domain $\Omega=\Omega_{\text{s}}\cup\Omega_{\text{f}}=[-1,1]$ split into fluid $\Omega_{\text{f}}\in[-1,0]$ and solid $\Omega_{\text{s}}\in[0,1]$, each containing the species $c_{\text{f}}$ and $c_{\text{s}}$ respectively. Both species are restricted to their respective sub-domain, i.e. $c_{\text{s}}=0$ in $\Omega_{\text{f}}$ and vice versa. Consider the following elementary reaction at the $x=0$ interface:
    
        \begin{equation}
            [\text{s}] \leftrightharpoons [\text{f}].
        \end{equation}
    As $c_{\text{s}}$ and $c_{\text{f}}$ are restricted to $\Omega_{\text{s}}$ and $\Omega_{\text{f}}$ respectively the reactive conditions of \cref{eqn:Surf_Reac_LoMA_Ideal,eqn::unresContinuity} become
    
        \begin{align}
            j(c_{\text{f}}) &= j(c_{\text{s}}), \quad x=0, \label{eqn::case3fluxcontinuity}\\
            j(c_{\text{f}}) &= k_{\text{f}} c_{\text{s}} - k_{\text{r}} c_{\text{f}}, \quad x=0, \label{eqn::case3reactiveflux}
        \end{align}
    where we have used the linear \textit{rate law} to model the reaction rate. Note $c_{\text{f}}$ and $c_{\text{s}}$ in \cref{eqn::case3fluxcontinuity,eqn::case3reactiveflux} are evaluated on $\Omega_{\text{f}}$ and $\Omega_{\text{s}}$'s side of $x=0$ respectively. We set both species as uncharged, i.e. $z_i=0$, to remove electrostatic effects and $u=1$ in $\Omega_{\text{f}}$ to ignore Stokes. We initially set $c_{\text{s}}=0$ in $\Omega_{\text{s}}$ and $c_{\text{f}} = e^{-200(x+0.5)^{2}}$ in $\Omega_{\text{f}}$ so that a clear increase in $c_{\text{s}}$ can be seen.
    
 		\begin{figure}[htbp]
 			\centering
 			\includegraphics[width=.5\linewidth]{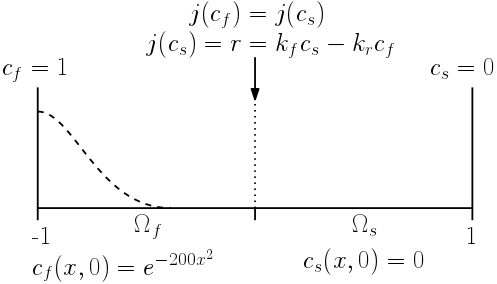}
 			\caption{{Graphical representation of domain $\Omega$ containing a solid and fluid with zero flux outer boundaries and reactive flux conditions at the shared interface.}}
 		\end{figure}
	
	Looking at the dimensionless numbers for the problem in \Cref{tab::DimNumsCase3Properties} we can notice here that advection dominates. We also include two more numbers from the reactive interface conditions at $x=0$. These are $\DA_{\text{II,f}}$ and $\DA_{\text{II,r}}$, the Damk\"ohler numbers given as the ratio between the diffusion rate and reaction rates, forward and reverse, respectively.
		
\begin{figure}[htbp]
            \centering
            \begin{subfigure}{.45\textwidth}
                \includegraphics[width=\linewidth]{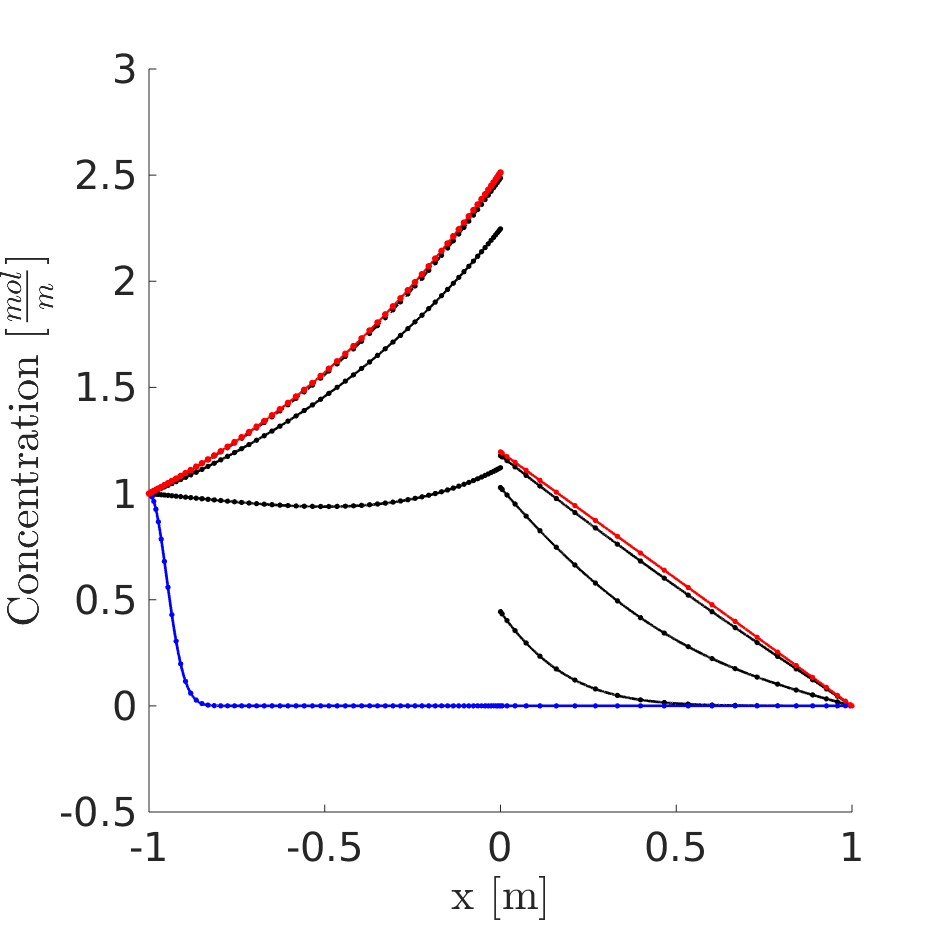}
            \end{subfigure}
            \begin{subfigure}{.45\textwidth}
                \includegraphics[width=\linewidth]{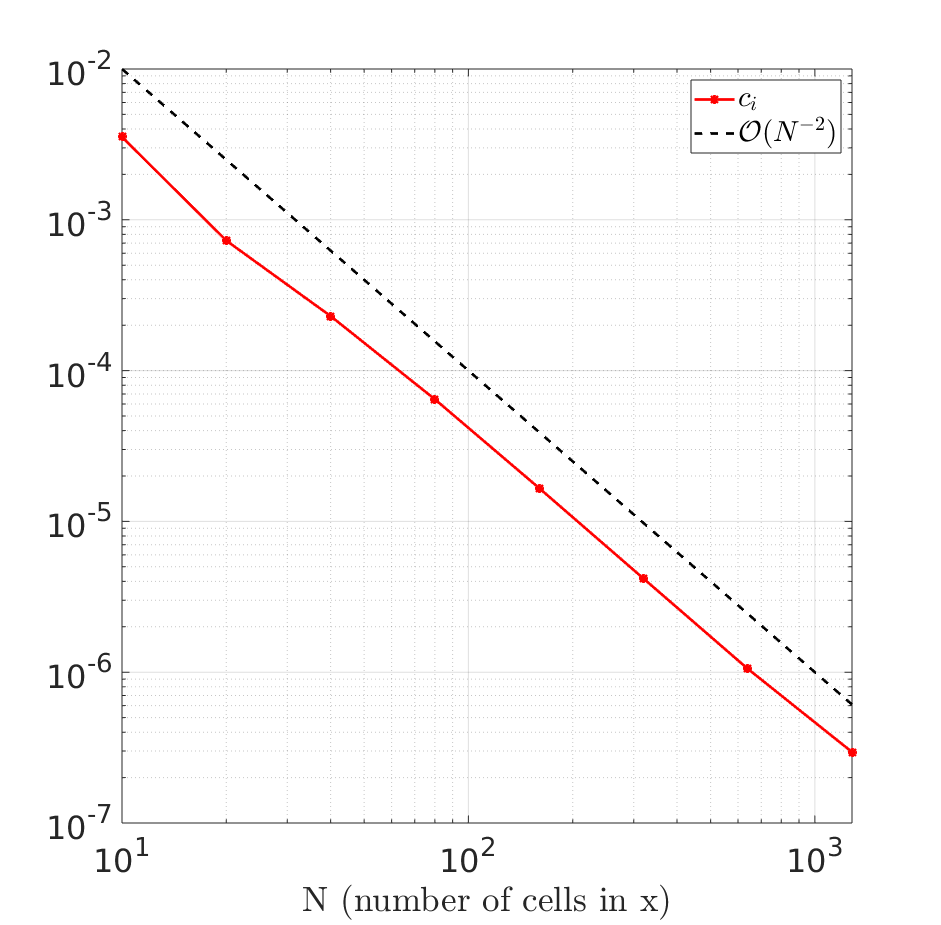}
            \end{subfigure}
            \caption{Left: {Molar concentration results of $c_{\text{f}}$ in fluid (dashed for \textit{pnpMultiFoam}, markers for \textit{Chebfun}) and $c_{\text{s}}$ in solid (solid for \textit{pnpMultiFoam}, markers for \textit{Chebfun}). Blue for initial concentrations, black at intermediate times, red for final time (i.e steady state)}, Right: {$L^2$ norm convergence plot between \textit{Chebfun} and mappedChemicalKinetics/pnpMultiFoam}}
            \label{fig::mappedChemL2Norm&Plot}
        \end{figure}

		\begin{table}[htbp]
    	    \centering
    	    \begin{tabular}{ccccccc}\hline
    	       Symbol & $\N$ & $\PE$ & $\L$ & $\P$ & $\DA_{\text{II,f}}$ & $\DA_{\text{II,r}}$ \\ \hline
    	       Value & $0$ & $1$ & $0$ & $0$ & $10$ & $100$ \\\hline
    	    \end{tabular}
    	    \caption{Table of dimensionless values}
    	    \label{tab::DimNumsCase3Properties}
    	\end{table}
    When running \textit{pnpMultiFoam} to show qualitative comparison we use a time step $\Delta t = 10^{-3}$s, second-order implicit time scheme \textit{backwards} and mesh discretisation of $N=1000$ cells (500 in $\Omega_{\text{f}}$, 500 in $\Omega_{\text{s}}$). To show convergence of results and compute the $L^2$ error point norm we move to the \textit{steadyState} time scheme to compare steady states. Results of the transient case are shown on the left in \cref{fig::mappedChemL2Norm&Plot} where we find good agreement. First an initial diffusion and advection of $c_{\text{f}}$ is seen towards the interface, once $c_{\text{f}}$ reaches $x=0$ it reacts to form $c_{\text{s}}$ in $\Omega_{\text{s}}$, where afterwards $c_{\text{s}}$ diffuses through $\Omega_{\text{s}}$. After $t=15$s a steady state is reached implying chemical equilibrium of the reaction.
 As for the $L^{2}$ error norm of the steady case, seen in \cref{fig::mappedChemL2Norm&Plot}, we find a second-order convergence of \textit{mappedChemicalKinetics} to \textit{Chebfun}. This is because the interface conditions are linear and therefore exactly approximated by the linear approximation.

\subsection{Random porous REV}
    
\begin{figure}[htbp]
            \centering
            \includegraphics[width=.45\linewidth]{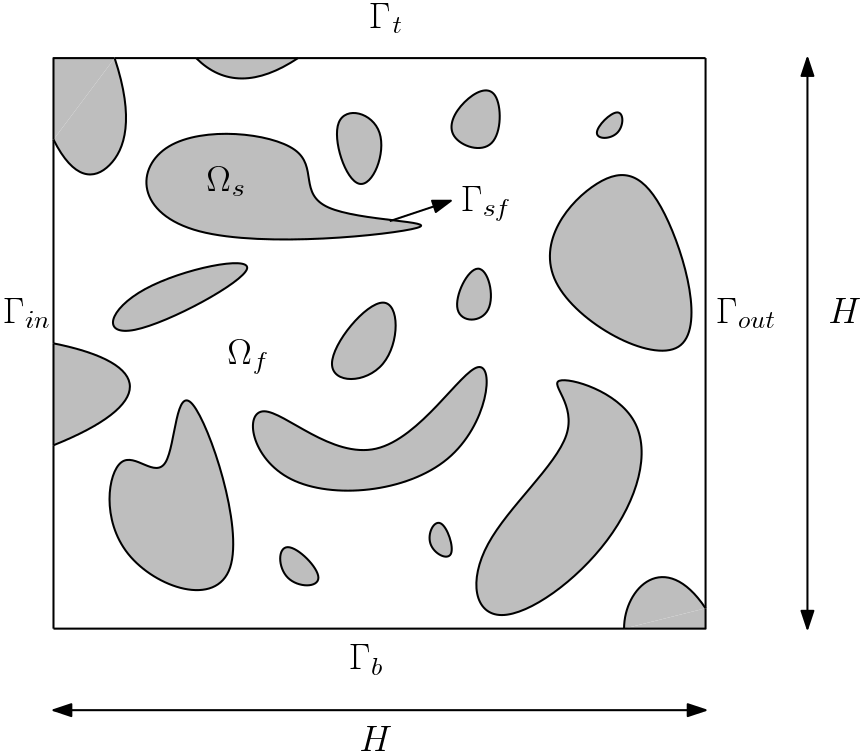}
            \caption{{Visual representation of the random porous domain containing solid $\Omega_{\text{s}}$ (grey) and fluid $\Omega_{\text{f}}$ (white) with shared interface $\Gamma_{\text{sf}}$.}}
            \label{fig::RandomPorousMedium}
        \end{figure}

    To demonstrate the capabilities of our solvers, we consider a randomly generated porous solid-fluid domain $\Omega=\Omega_{\text{s}}\cup\Omega_{\text{f}}$, see \cref{fig::RandomPorousMedium}. The random generation is done by a truncated Gaussian random field \cite{Icardi:2022aa,Municchi2021}, sampled at the mesh points and categorized into two bins, denoting solid and fluid cells, by a threshold. Raising or lowering this threshold alters the porosity of the domain.    
    
    To have $\Omega$ be a microscopic representative elementary volume (REV) of a much larger macroscopic porous medium, we apply periodic conditions along the outer boundaries $\Gamma_{\text{ext}}= \Gamma_{\text{in}}\cup\Gamma_{\text{out}}\cup\Gamma_{\text{t}}\cup\Gamma_{\text{b}}$. To generate movement, we fix a jump in $\phi$ between $\Gamma_{\text{in,s}}$ and $\Gamma_{\text{out,s}}$, where the subscript $s$ denotes the section of $\Gamma_{\text{in}}$ neighbouring $\Omega_{\text{s}}$. This is to mimic a fixed applied potential difference across the macroscopic medium, such as an applied voltage across a battery cell. We consider two ion species $c_1$ and $c_2$ with opposite valencies $z_1=-z_2=-1$. The height and length of the region is set as $H=1\times 10^{-4}$m.

    In $\Omega_{\text{f}}$ we start with a uniform concentration of $c_2$, the same is done for $\Omega_{\text{s}}$ with $c_1$. Fluid is initially taken at rest. Along $\Gamma_{\text{sf}}$ we set no flux for $c_1$ and flux continuity for $c_2$. For $\phi$ we assume continuity of the electric displacement $\mathbf{D}=\varepsilon\mathbf{E}$ along $\Gamma_{\text{sf}}$. In \Cref{tab::DimNumsCase4Properties} we list the dimensionless number values of the case.
    
            \begin{table}[htbp]
    	    \centering
    	    \begin{tabular}{ccccc}\hline
    	       Symbol & $\N$ & $\PE$ & $\L$ & $\P$ \\ \hline
    	       Value & $0.4$ & $9.65\times10^{-4}$ & $0.95$ & $1$\\\hline
    	    \end{tabular}
    	    \caption{Dimensionless numbers for the random porous media testcase.}
    	    \label{tab::DimNumsCase4Properties}
    	\end{table}
	

    The system is initialized with constant $c_1$ and $c_2$ in $\Omega_{\text{f}}$ and $\Omega_{\text{s}}$ respectively. All fields are periodic on the outer boundaries $\Gamma_{\text{ext}}=\Gamma_{\text{in}}\cup\Gamma_{\text{out}}\cup\Gamma_{\text{t}}\cup\Gamma_{\text{b}}$, apart from $\phi$ where we introduce a quasi-periodic condition with a fixed jump between $\Gamma_{\text{in,s}}$ and $\Gamma_{\text{out,s}}$, denoting the parts of $\Gamma_{\text{in}}$ and $\Gamma_{\text{out}}$ neighbouring $\Omega_{\text{s}}$. We apply no flux and flux continuity for $c_1$ and $c_{\text{s}}$ along $\Gamma_{\text{sf}}$ respectively. For $\phi$ we assume continuity of the electric displacement $\mathbf{D}=\varepsilon\mathbf{E}$ along $\Gamma_{\text{sf}}$:
    
        \begin{align}
            c_{1}(\bx,0) = 0\,,\quad c_{2}(\bx,0) = 1\times 10^{-8}\,,\quad \phi(\bx,0)=0, \quad &\bx\in\Omega_{\text{s}}; \\
            c_{1}(\bx,0) = 1\times 10^{-8}\,,\quad c_{2}(\bx,0) = 0\,,\quad
              \phi(\bx,0)=0\,,\quad  \bu (\bx,0)=\mathbf{0}\,,\quad  p(\bx,0)=0, \quad &\bx\in\Omega_{\text{f}}; \\
            \evalu{\bm{\nabla}_{\mathbf{n}}a}{\Gamma_{\text{in,s}}}{\Gamma_{\text{out,s}}} =0
            \,,\quad  \evalu{a}{\Gamma_{\text{in,s}}}{\Gamma_{\text{out,s}}} =0\,,\quad \evalu{\bm{\nabla}_{\mathbf{n}}\phi}{\Gamma_{\text{in,s}}}{\Gamma_{\text{out,s}}} = 0\,,\quad
            \evalu{\phi}{\Gamma_{\text{in,s}}}{\Gamma_{\text{out,s}}} = 1\times 10^{-2}, \quad & a\in\{c_{1},c_{2}\}; \\
            \evalu{\bm{\nabla}_{\mathbf{n}}a}{\Gamma_{\text{in,f}}}{\Gamma_{out,f}} =0
            \,,\quad  \evalu{a}{\Gamma_{\text{in,f}}}{\Gamma_{out,f}} =0, \quad & a\in\{c_{1},c_{2},\phi,p,\bu\}; \\
            \evalu{\bm{\nabla}_{\mathbf{n}}a}{\Gamma_{\text{t}}}{\Gamma_{\text{b}}} =0
            \,,\quad  \evalu{a}{\Gamma_{\text{t}}}{\Gamma_{\text{b}}} =0, \quad & a\in\{c_{1},c_{2},\phi,p,\bu\} ;\\
            \flux_{1}\cdot\mathbf{n} = 0\,,\quad \evalu{\flux_{2}\cdot\mathbf{n}}{\Gamma_{\text{sf,f}}}{\Gamma_{\text{sf,s}}}=0\,,\quad
            \bm{\nabla}p\cdot\mathbf{n}=0\,,\quad \mathbf{u}\cdot\mathbf{n}=0\,,\quad
            \evalu{\varepsilon\bm{\nabla}_{\mathbf{n}}\phi}{\Gamma_{\text{sf,f}}}{\Gamma_{\text{sf,s}}}=0, \quad &\bx\in\Gamma_{\text{sf}}.
        \end{align}
    	
\begin{figure}[htbp]
        \centering
        \begin{subfigure}{.32\textwidth}
            \includegraphics[width=\linewidth]{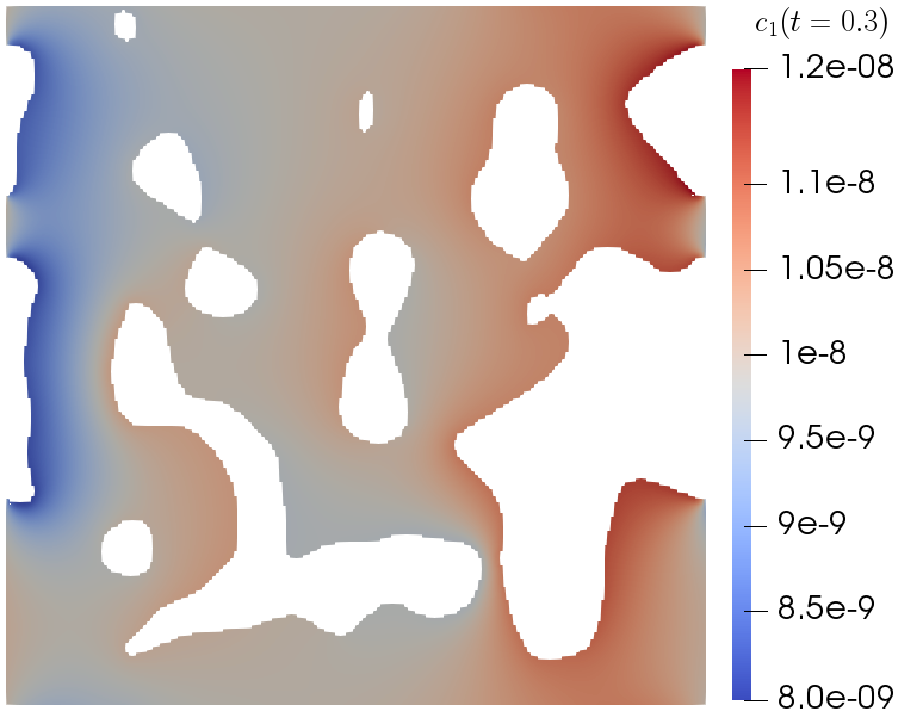}
            \caption{}
        \end{subfigure}
        \begin{subfigure}{.32\textwidth}
            \includegraphics[width=\linewidth]{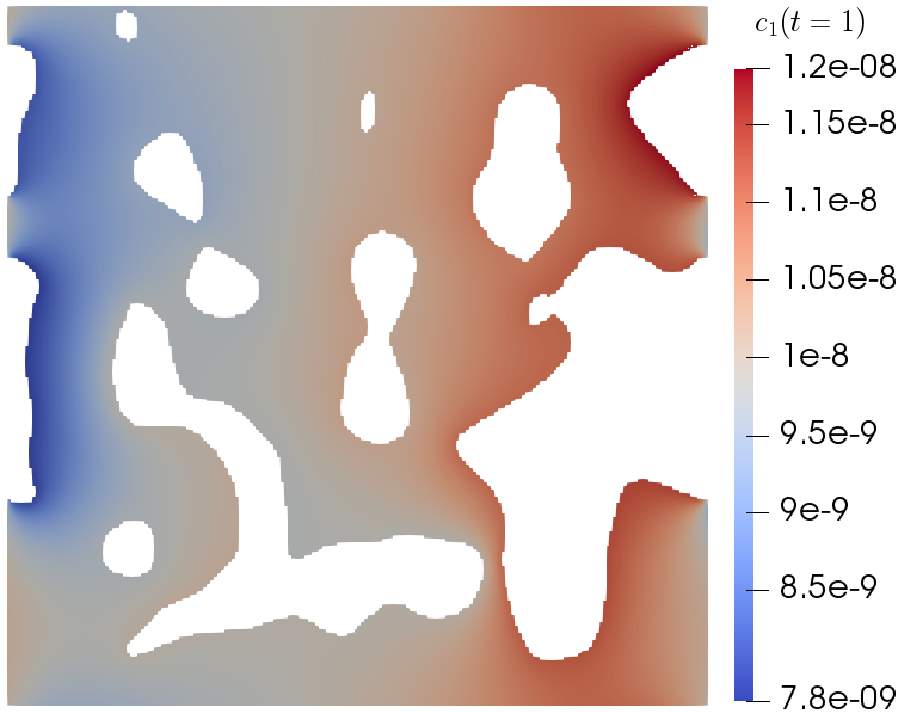}
            \caption{}
        \end{subfigure}
        \begin{subfigure}{.32\textwidth}
            \includegraphics[width=\linewidth]{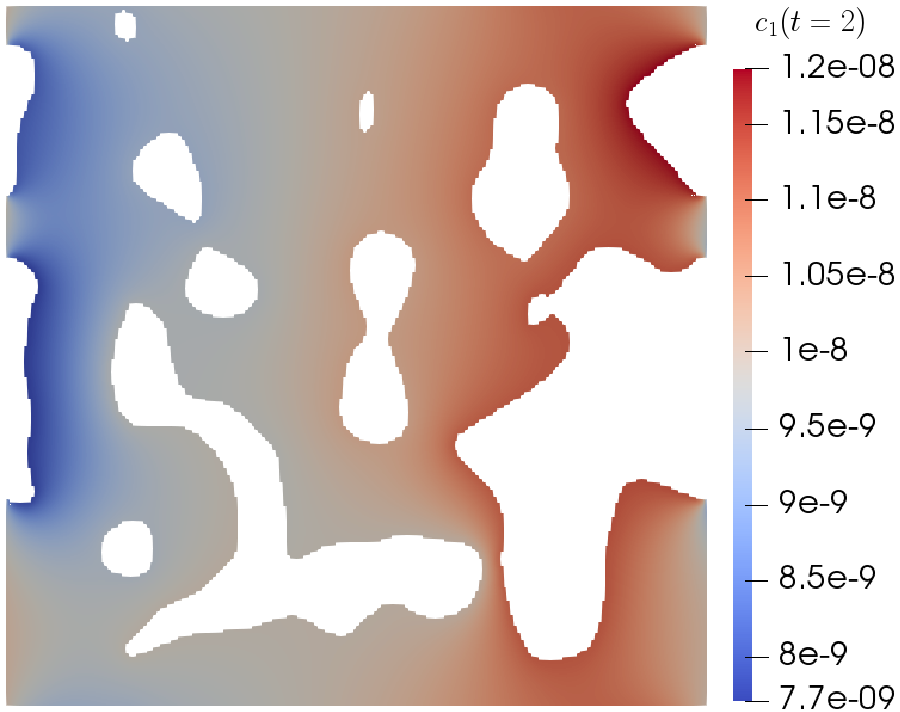}
            \caption{}
        \end{subfigure}
        \begin{subfigure}{.32\textwidth}
            \includegraphics[width=\linewidth]{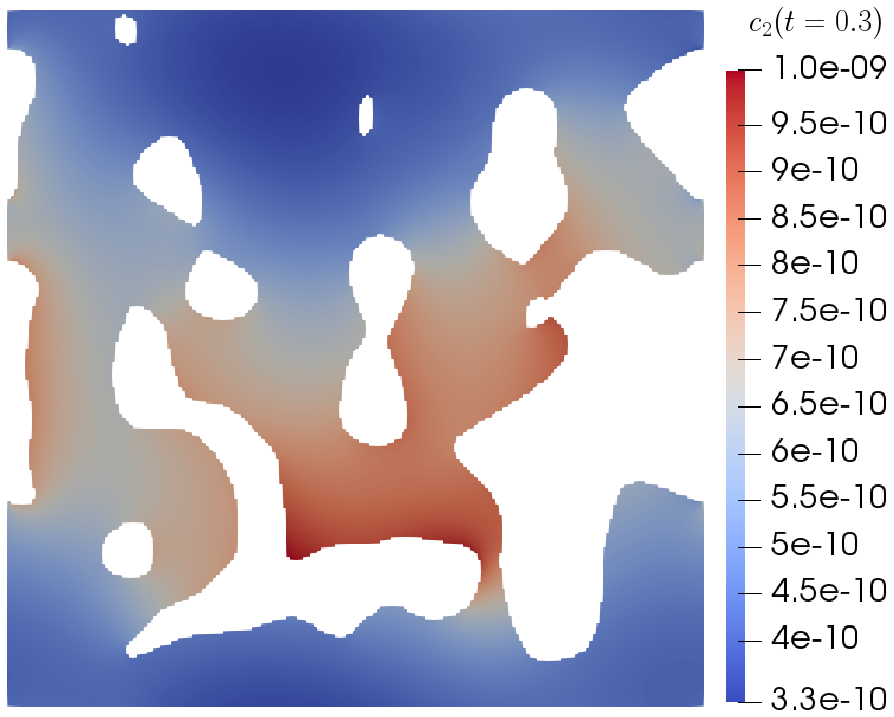}
            \caption{}
        \end{subfigure}
        \begin{subfigure}{.32\textwidth}
            \includegraphics[width=\linewidth]{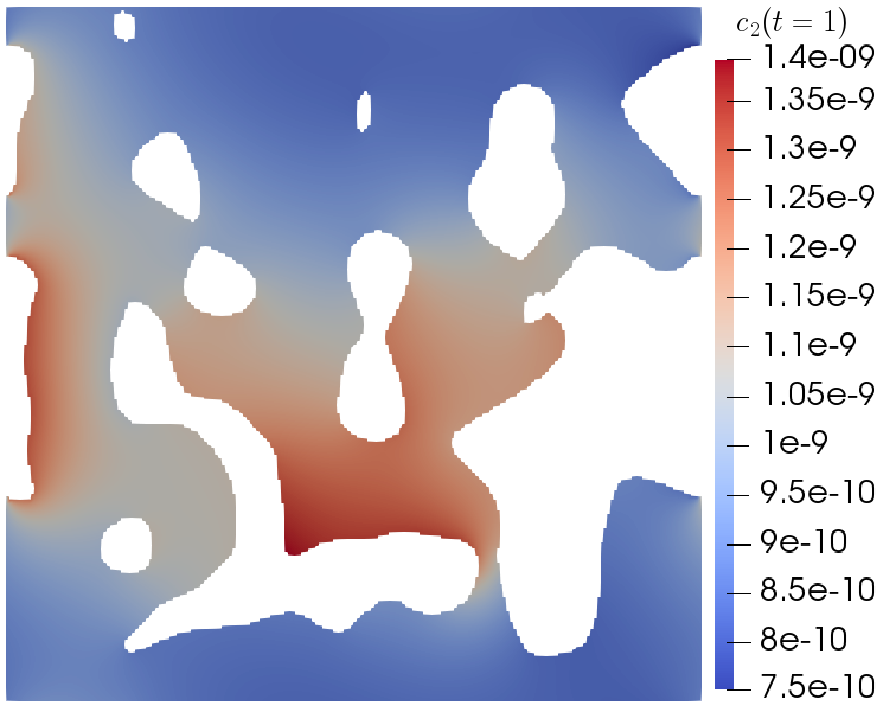}
            \caption{}
        \end{subfigure}
        \begin{subfigure}{.32\textwidth}
            \includegraphics[width=\linewidth]{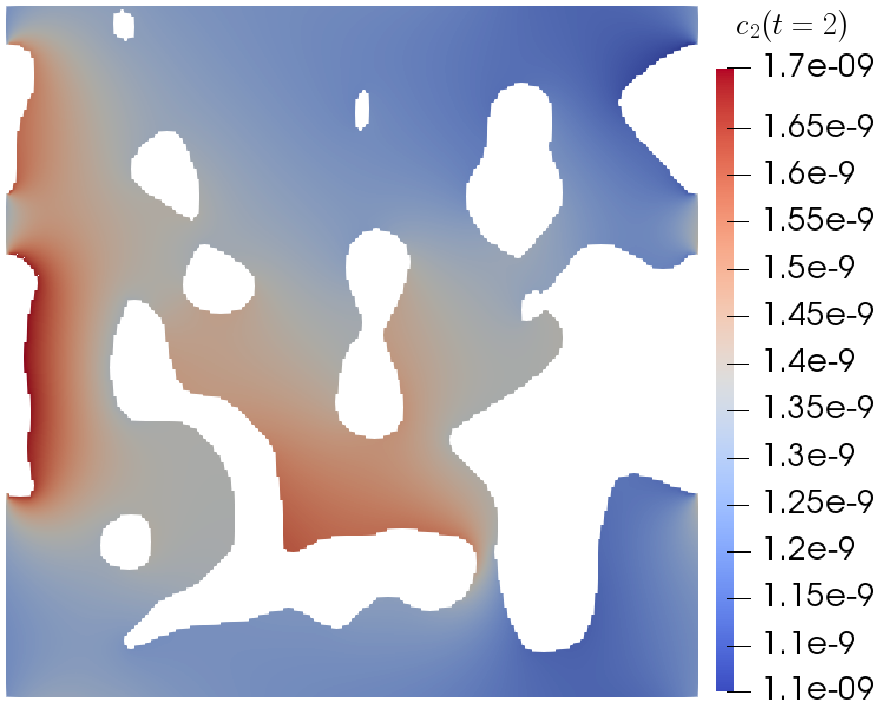}
            \caption{}
        \end{subfigure}
        \begin{subfigure}{.32\textwidth}
            \includegraphics[width=\linewidth]{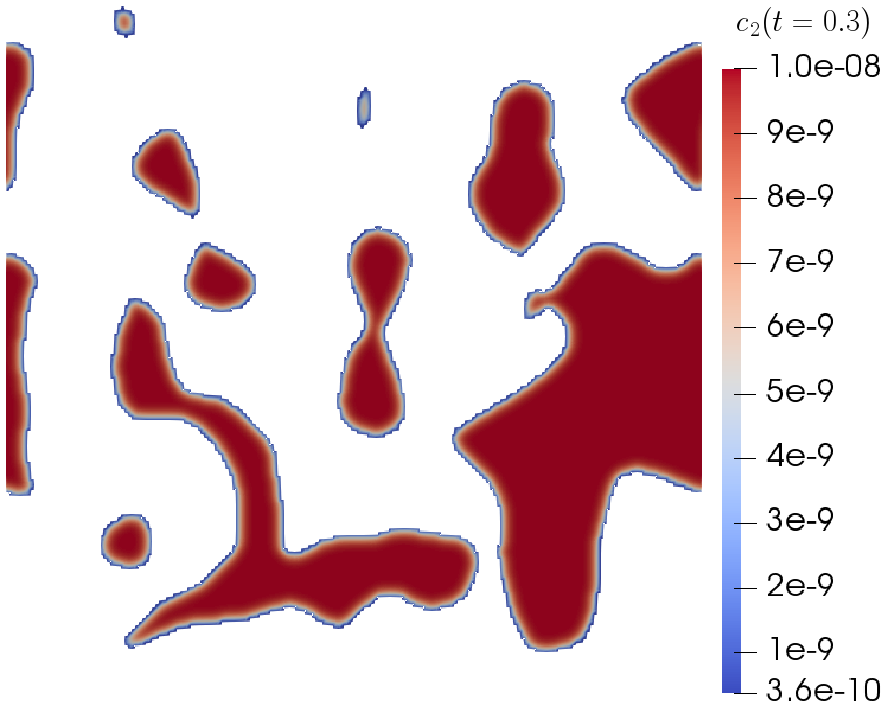}
            \caption{}
        \end{subfigure}
        \begin{subfigure}{.32\textwidth}
            \includegraphics[width=\linewidth]{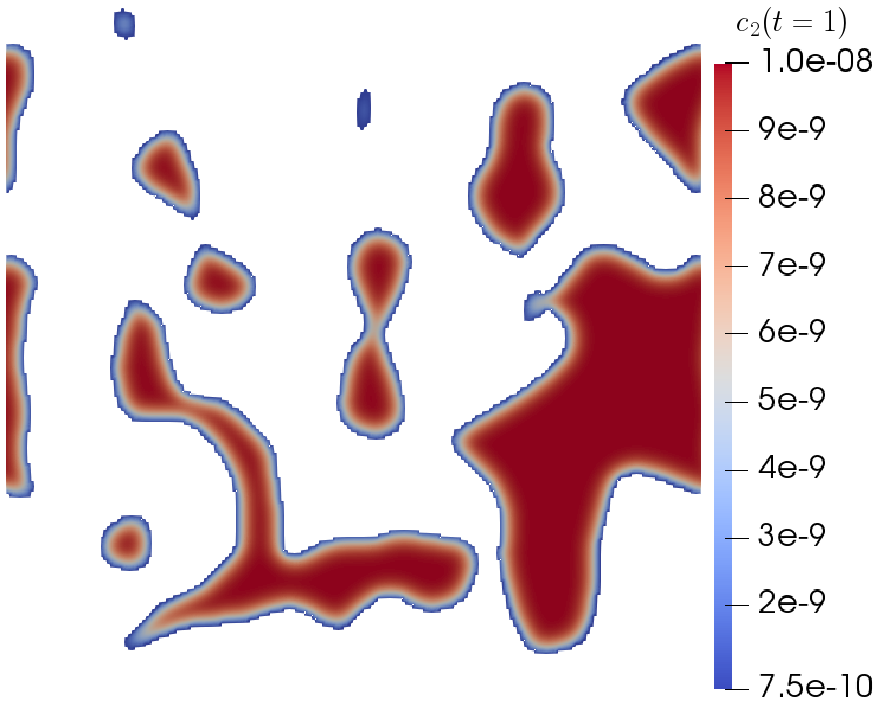}
            \caption{}
        \end{subfigure}
        \begin{subfigure}{.32\textwidth}
            \includegraphics[width=\linewidth]{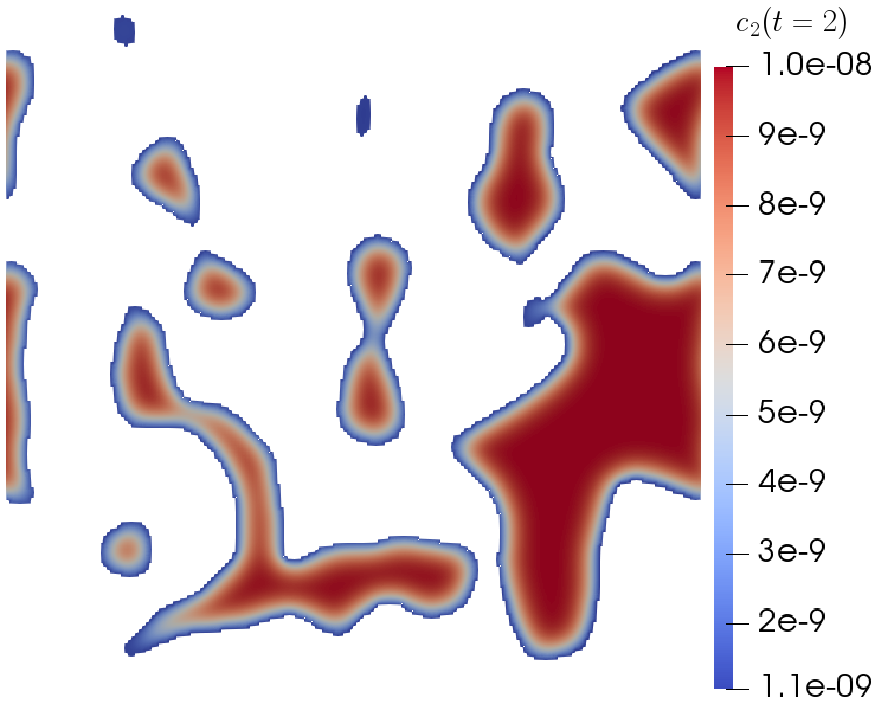}
            \caption{}
        \end{subfigure}
        \begin{subfigure}{.32\textwidth}
            \includegraphics[width=\linewidth]{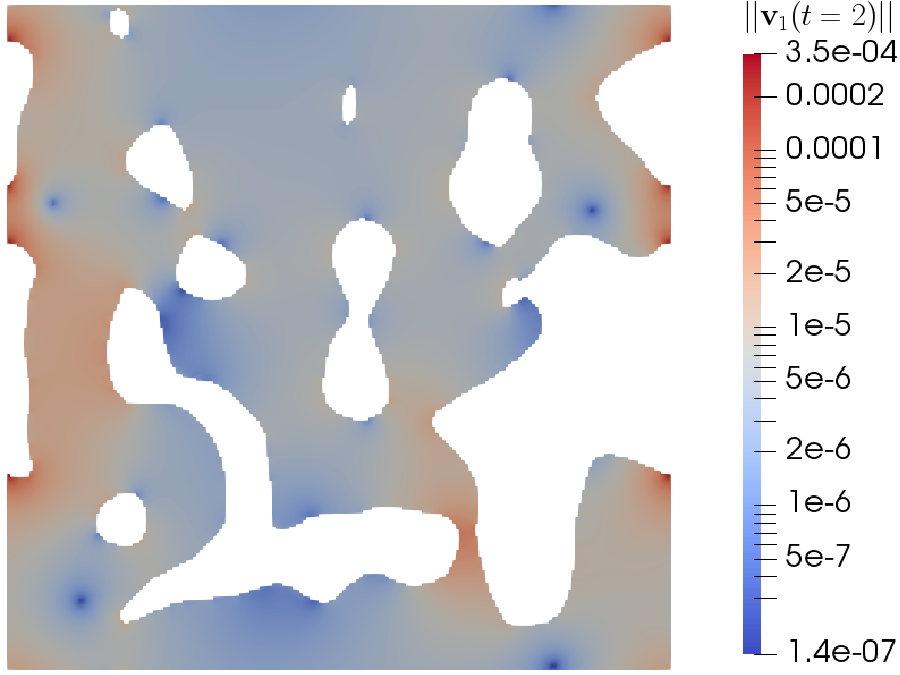}
            \caption{}
        \end{subfigure}
        \begin{subfigure}{.32\textwidth}
            \includegraphics[width=\linewidth]{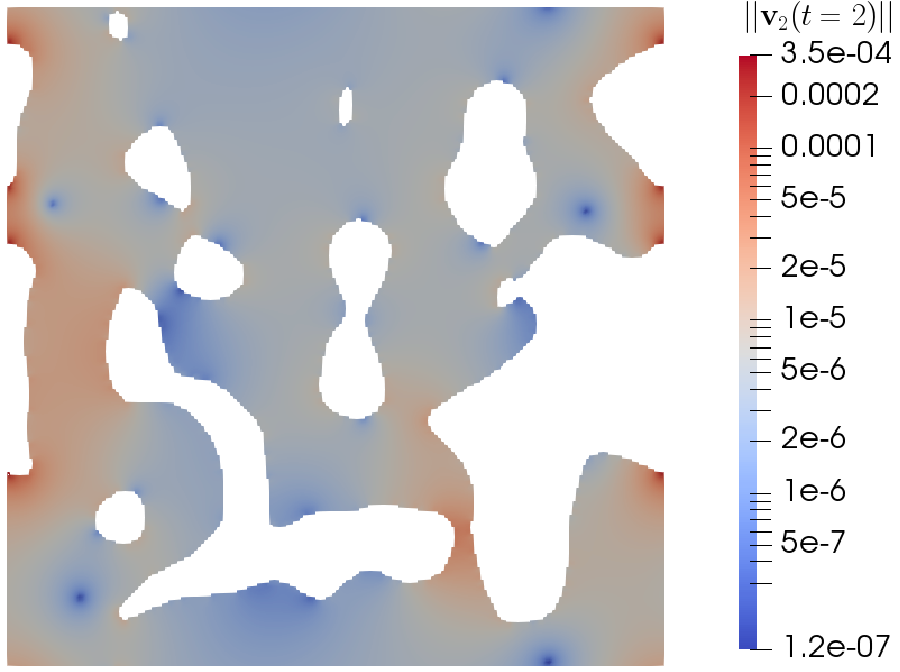}
            \caption{}
        \end{subfigure}
        \begin{subfigure}{.32\textwidth}
            \includegraphics[width=\linewidth]{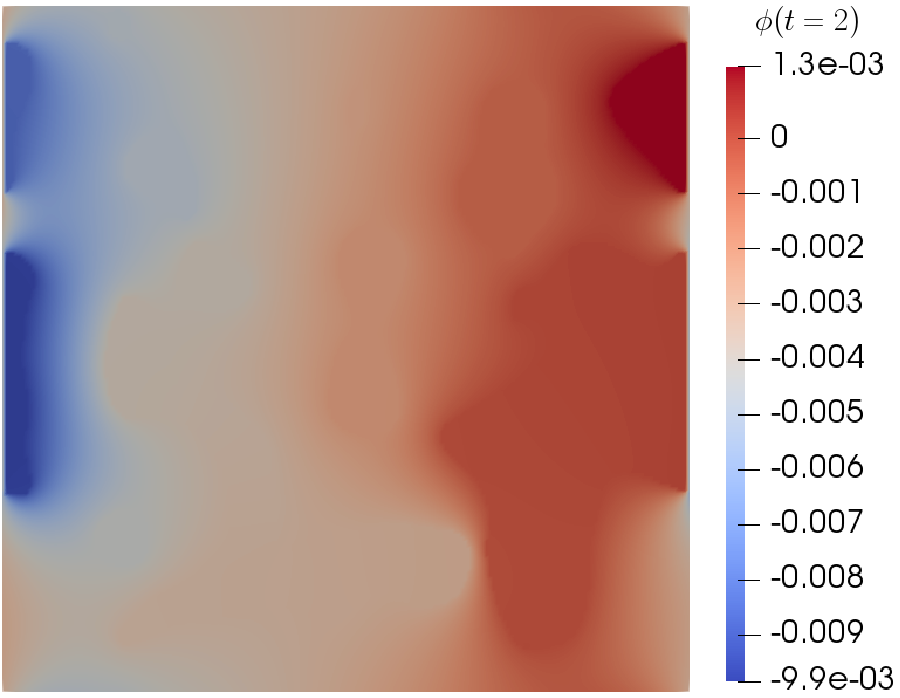}
            \caption{}
        \end{subfigure}
        \caption{{Results of the first random solid-fluid geometry, discretized using $4\times10^{4}$ cells.} (a)-(f): {Profiles of $c_1$ and $c_2$ at $t=0.3$, $1$, $2$s within $\Omega_{\text{f}}$.} (g)-(i): {Profiles of $c_2$ at $t=0.3$, $1$, $2$s within $\Omega_{\text{s}}$.} (j)-(l): {Profiles of $||\mathbf{v}_1||$, $||\mathbf{v}_2||$ and $\phi$ at $t=2$s.}}
        \label{fig::seed2_200x200Cells}
    \end{figure}

%
    To show consistency in the results we run four simulations in total. This is done for two randomly generated meshes, each discretized over two levels of refinement for $N=1\times10^{4}$ and $N=4\times10^{4}$ cells. We will use letters $A$ and $B$ to differentiate between the random generations and subscripts $1$ and $2$ for the levels of refinement. E.g., $A_{2}$ denotes the first random mesh, refined using $4\times10^{4}$ cells. All runs use a time step $\Delta t = 5\times10^{-4}$ and second-order implicit time scheme \textit{backwards}. To observe the advective and electrostatic effects on our ions we define the following velocities in $\Omega_{\text{f}}$:
    
        \begin{equation}
            \mathbf{v}_{i} = \bu - \frac{D_i z_i F}{RT}\bm{\nabla}\phi,
        \end{equation}
    such that the continuity \cref{eqn:PNP_mass_contin} in $\Omega_{\text{f}}$ may be written as
    
        \begin{equation}
            \pardiff{c_i}{t} + \dive{\mathbf{v}_i c_i - D_{i}\bm{\nabla}c_i} = 0.
        \end{equation}
    \Cref{fig::seed2_200x200Cells} shows the results of the first random realization with a finer mesh ($N=4\times 10 ^{4}$), where we observe a gradient in $\phi$ formed between regions of $\Omega_{\text{s}}$ connected at $\Gamma_{\text{in}}$ and $\Gamma_{\text{out}}$. This is primarily due to our jump condition for $\phi$ applied along $\Gamma_{\text{in,s}}$ and $\Gamma_{\text{out,s}}$. Whilst the movement of our ion species has an effect on this gradient, for the most part, it remains steady.
    The applied electrical driving force causes the uniform concentration $c_1$ in $\Omega_{\text{f}}$ to accumulate around the regions of $\Omega_{\text{s}}$ connected at $\Gamma_{\text{out}}$, where a higher electric potential is present. Conversely, around $\Omega_{\text{s}}$ connected at $\Gamma_{\text{in}}$ we find a much lower, but not zero, concentration of $c_1$. A similar observation is found with $c_2$, once diffused out of $\Omega_{\text{s}}$, accumulating where there is a lower potential. As with $c_1$, at the region of $\Omega_{\text{s}}$ connected at $\Gamma_{\text{out}}$ there is a much lower, but again not zero, concentration of $c_2$. We observe therefore the formation of EDLs, see \S\ref{sec::dimensionalAnalysis}, around these connected solids. Given more time to evolve these EDLs would become further apparent, with more of $c_2$ diffusing out of $\Omega_{\text{s}}$. As for $\mathbf{v}_1$ and $\mathbf{v}_2$ we find their magnitudes $||\mathbf{v}_1||$ and $||\mathbf{v}_2||$ are highest within these EDLs located at the connected solids due to higher gradients in $\phi$.

\begin{figure}[htbp]
        \centering
        \begin{subfigure}{.32\textwidth}
            \includegraphics[width=\linewidth]{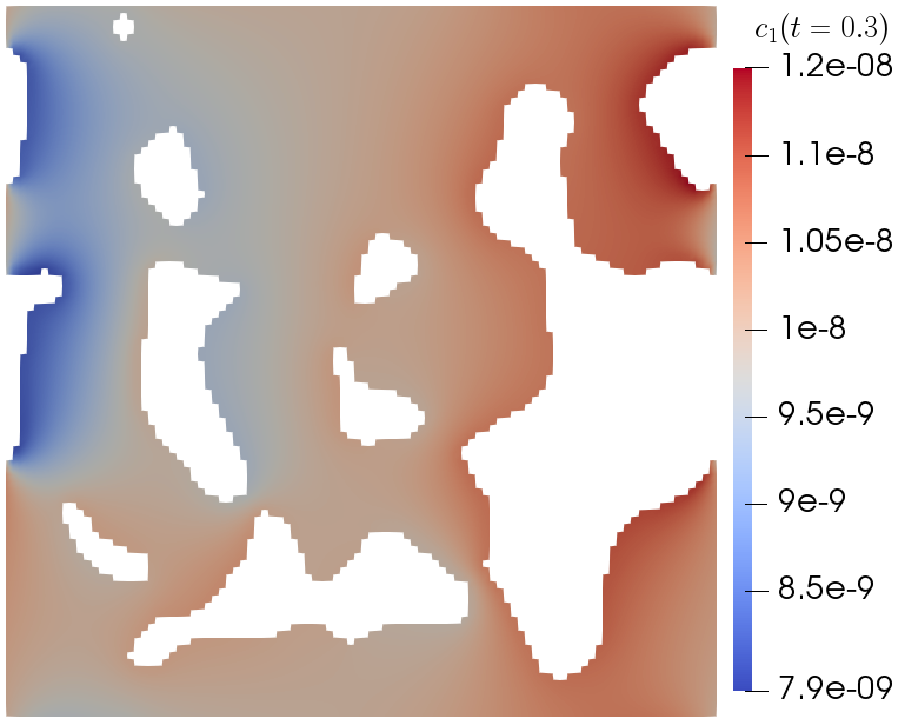}
        \end{subfigure}
        \begin{subfigure}{.32\textwidth}
            \includegraphics[width=\linewidth]{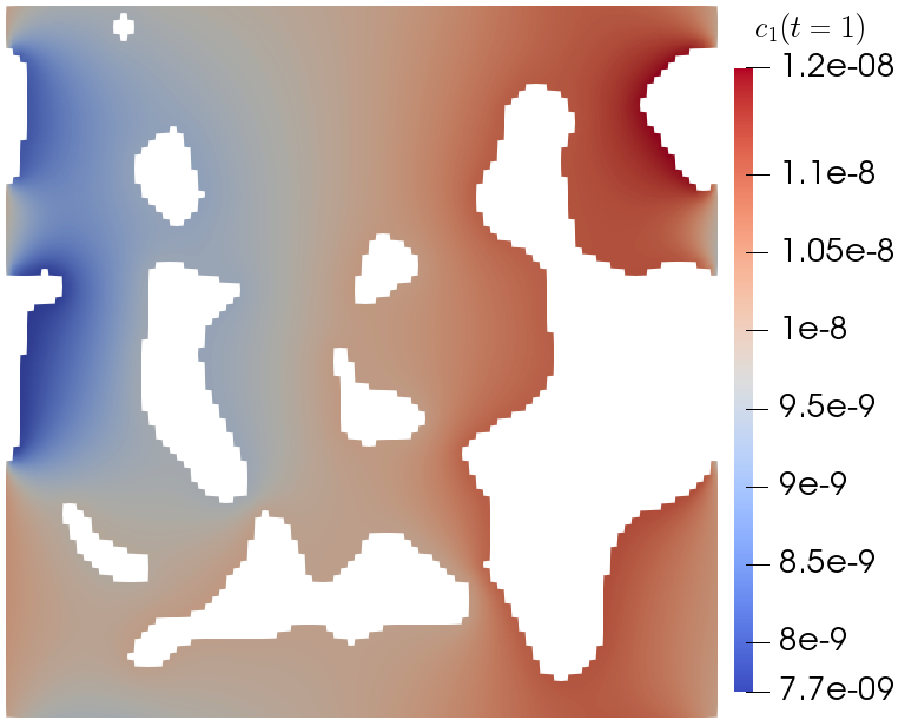}
        \end{subfigure}
        \begin{subfigure}{.32\textwidth}
            \includegraphics[width=\linewidth]{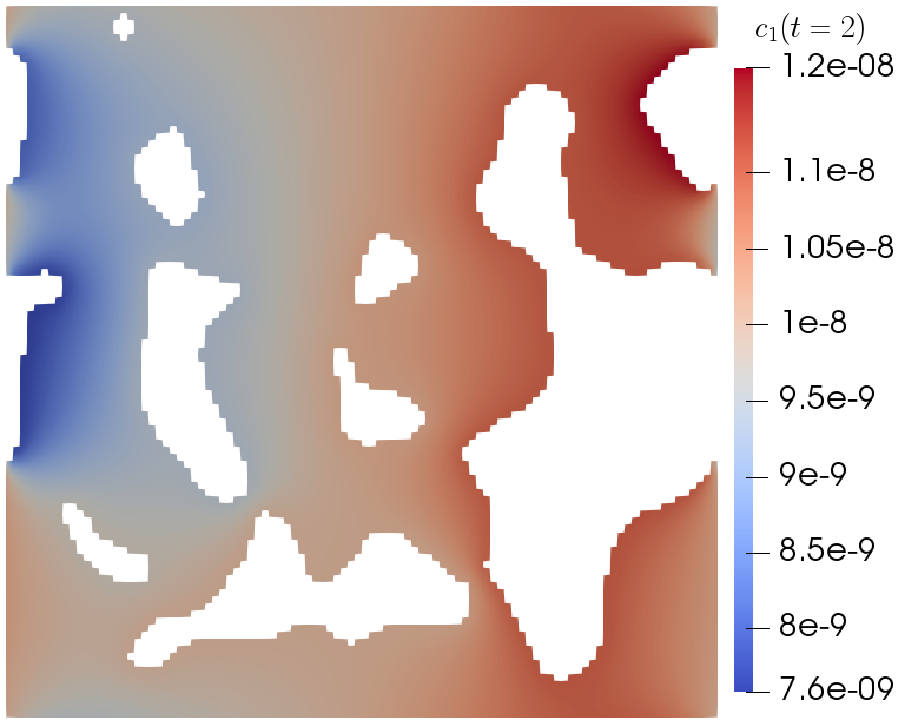}
        \end{subfigure}
        \begin{subfigure}{.32\textwidth}
            \includegraphics[width=\linewidth]{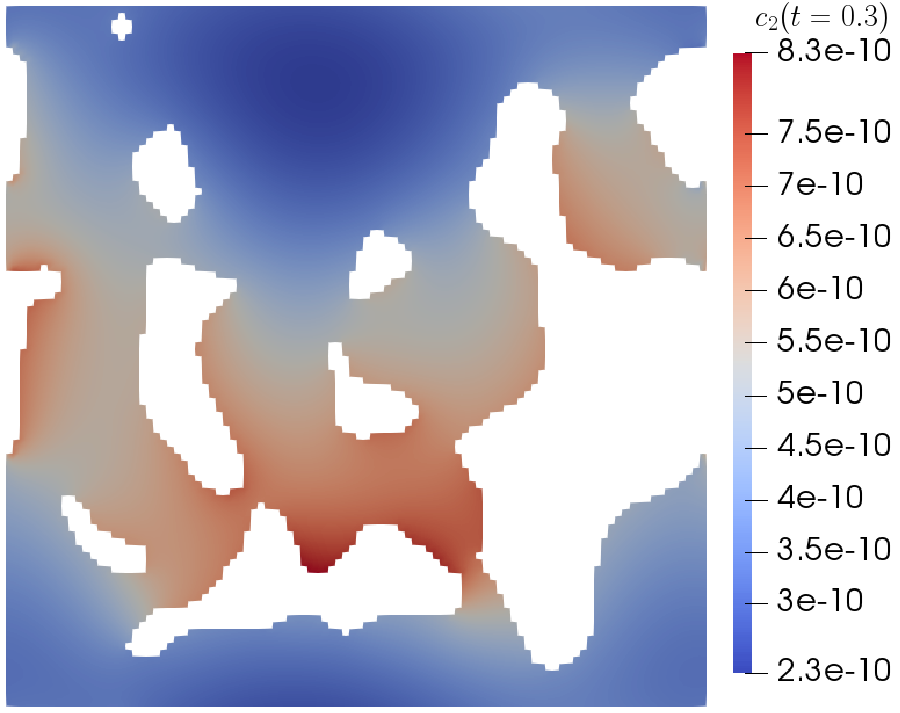}
        \end{subfigure}
        \begin{subfigure}{.32\textwidth}
            \includegraphics[width=\linewidth]{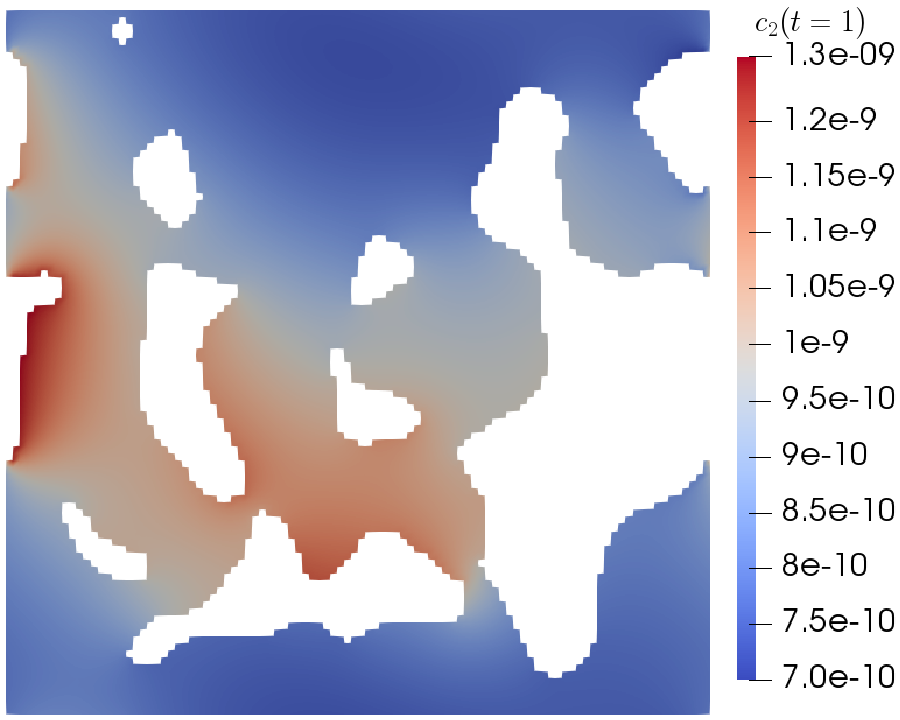}
        \end{subfigure}
        \begin{subfigure}{.32\textwidth}
            \includegraphics[width=\linewidth]{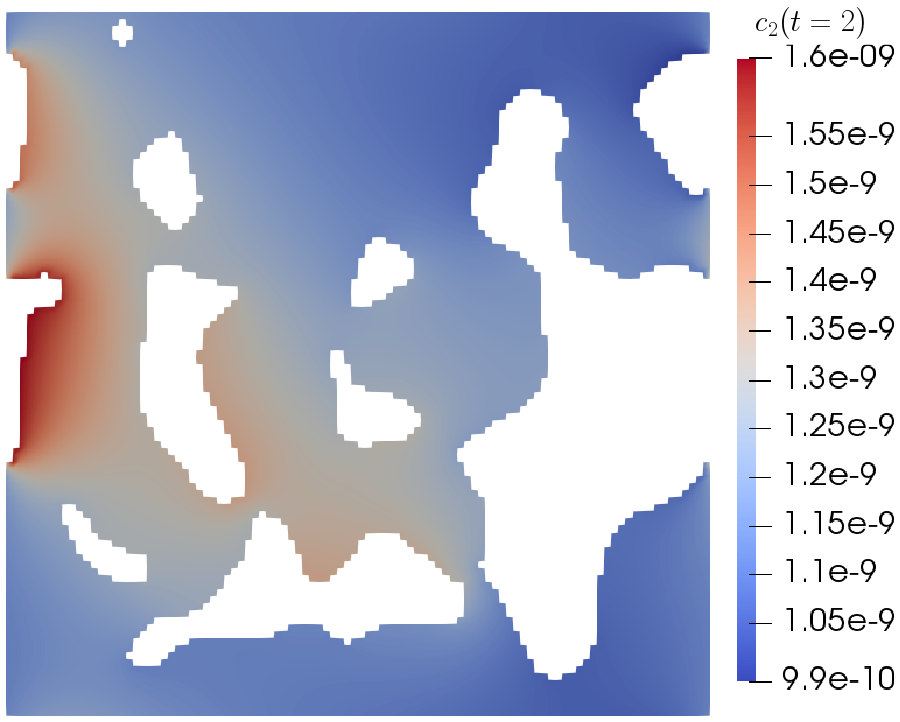}
        \end{subfigure}
        \begin{subfigure}{.32\textwidth}
            \includegraphics[width=\linewidth]{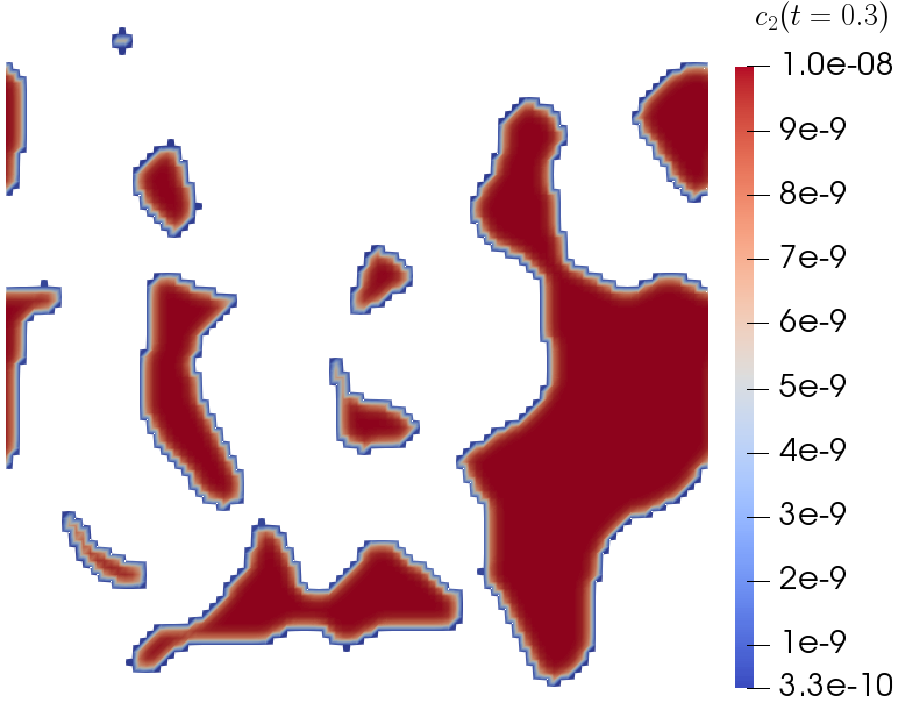}
        \end{subfigure}
        \begin{subfigure}{.32\textwidth}
            \includegraphics[width=\linewidth]{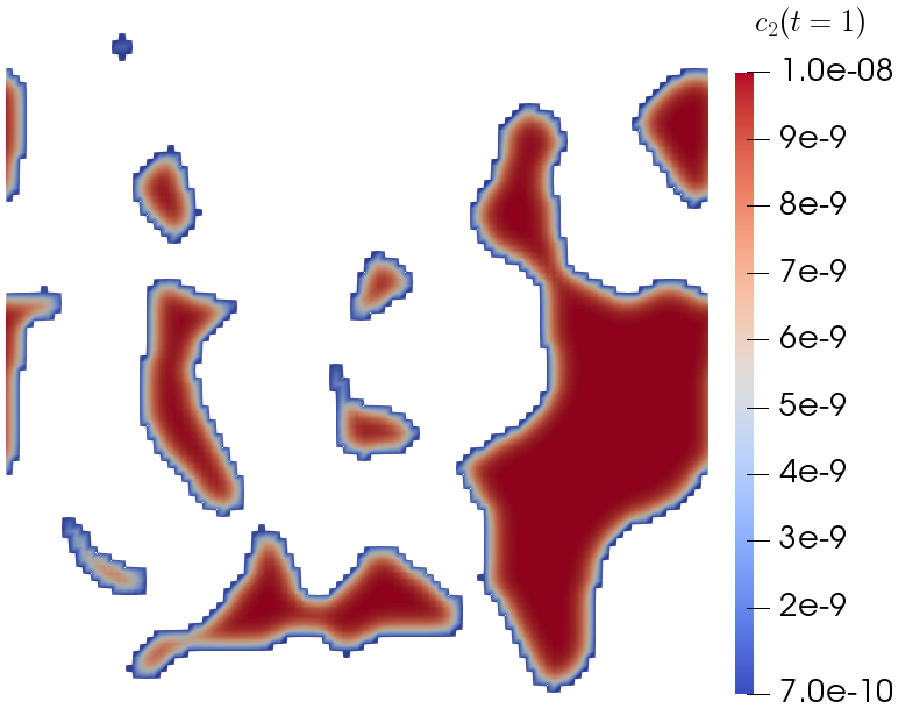}
        \end{subfigure}
        \begin{subfigure}{.32\textwidth}
            \includegraphics[width=\linewidth]{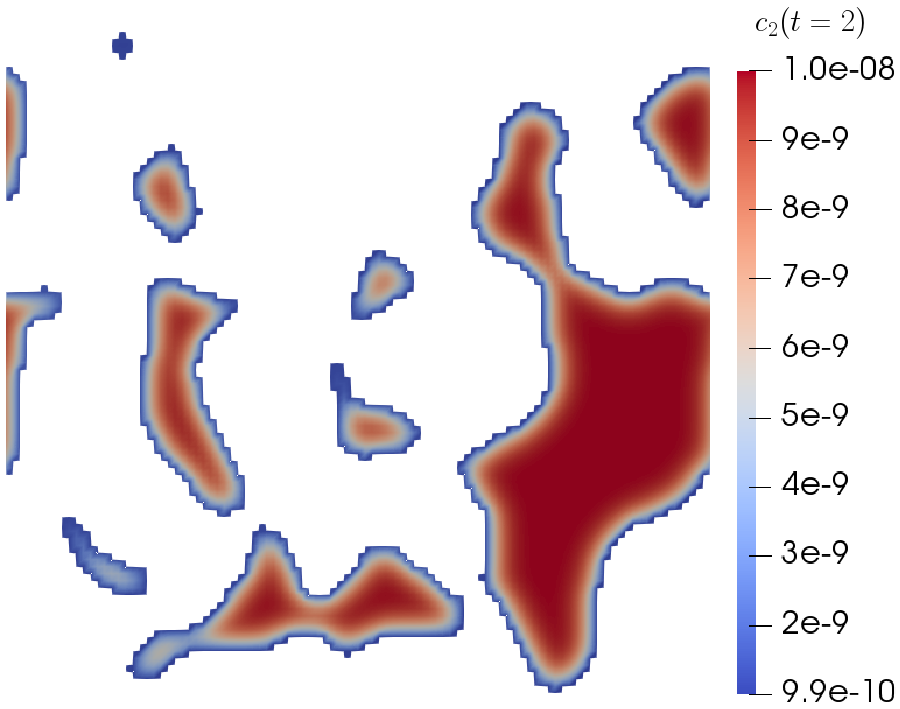}
        \end{subfigure}
        \begin{subfigure}{.32\textwidth}
            \includegraphics[width=\linewidth]{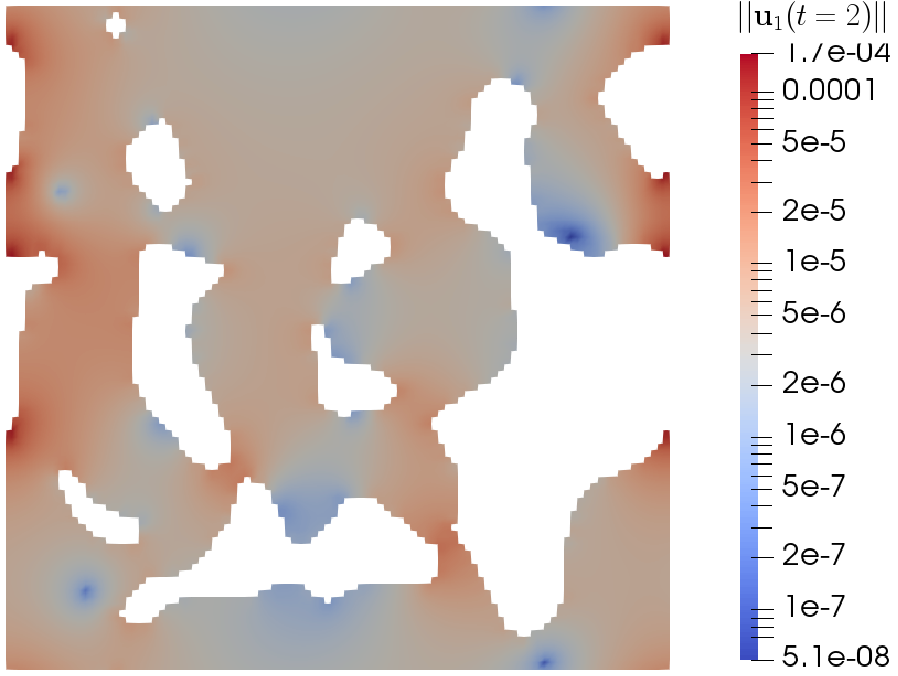}
        \end{subfigure}
        \begin{subfigure}{.32\textwidth}
            \includegraphics[width=\linewidth]{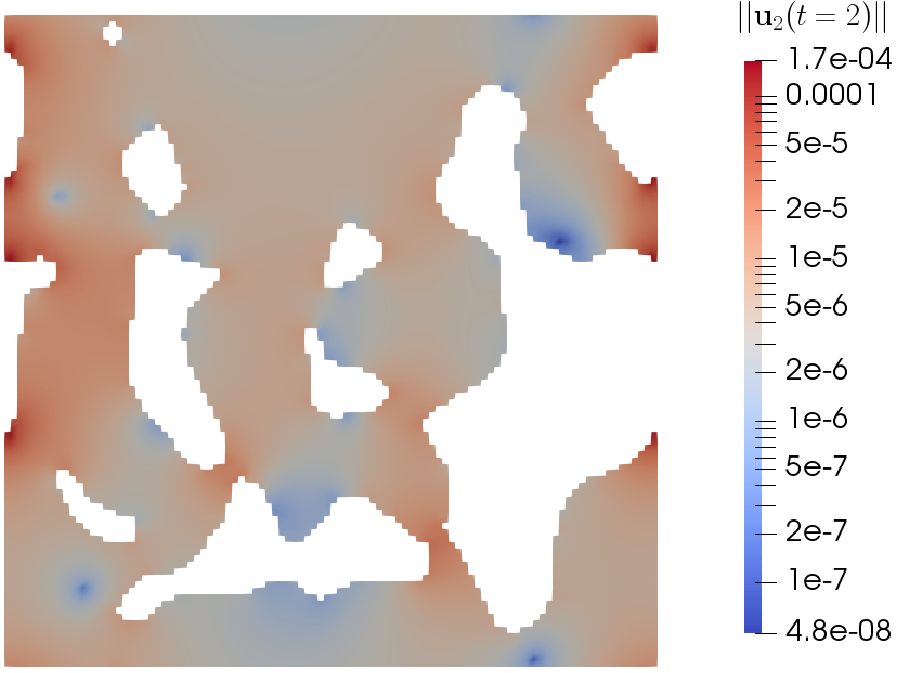}
        \end{subfigure}
        \begin{subfigure}{.32\textwidth}
            \includegraphics[width=\linewidth]{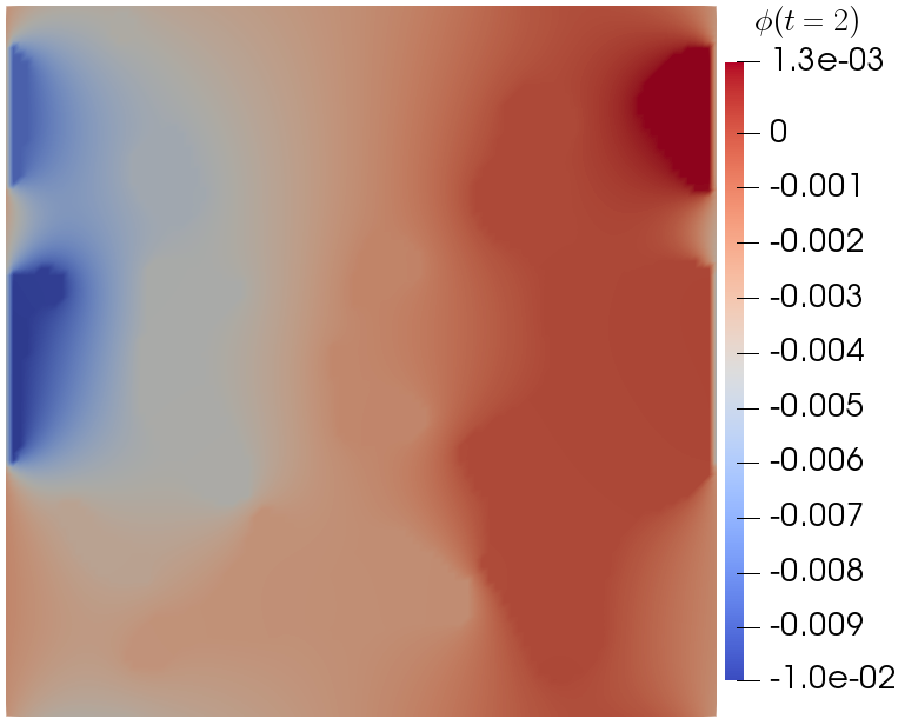}
        \end{subfigure}
        \caption{Results of the first random solid-fluid geometry, discretized using $N=1\times10^{4}$ cells. Top two rows: {Profiles of $c_1$ and $c_2$ at $t=0.3$, $1$, $2$s within $\Omega_{\text{f}}$.} Third row: {Profiles of $c_2$ at $t=0.3$, $1$, $2$s within $\Omega_{\text{s}}$.} Bottom row: {Profiles of $||\mathbf{u}_1||$, $||\mathbf{u}_2||$ and $\phi$ at $t=2$s.}}
        \label{fig::seed2_100x100Cells}
    \end{figure}

\begin{figure}[htbp]
        \centering
        \begin{subfigure}{.32\textwidth}
            \includegraphics[width=\linewidth]{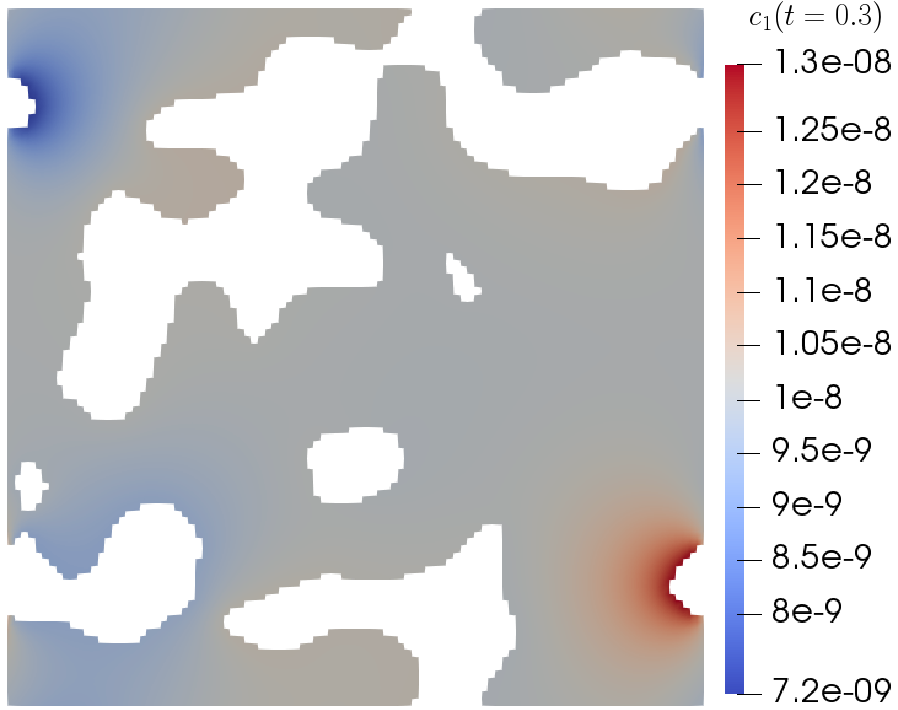}
        \end{subfigure}
        \begin{subfigure}{.32\textwidth}
            \includegraphics[width=\linewidth]{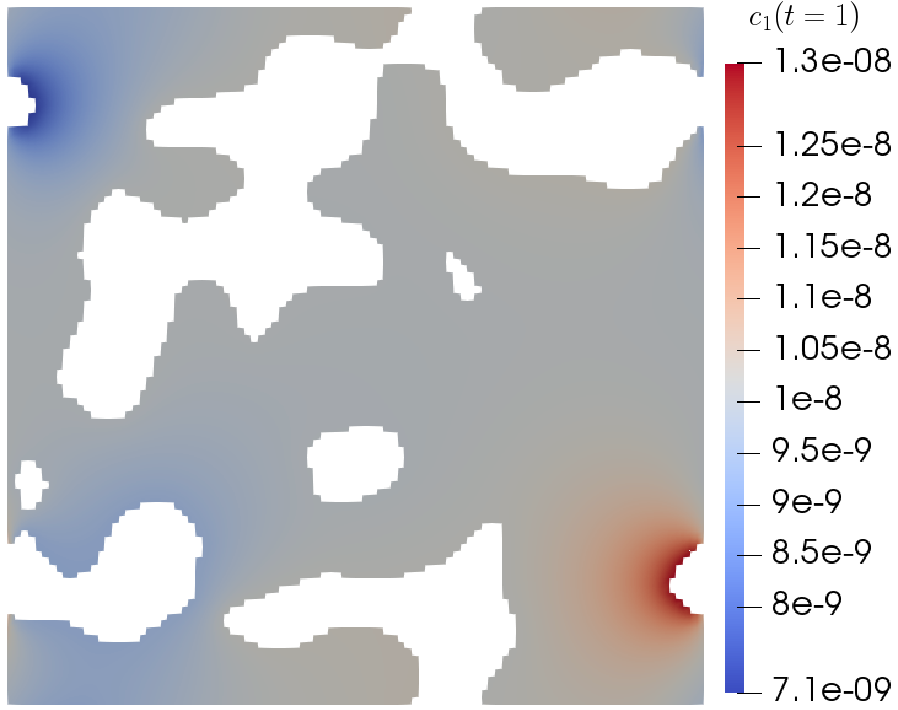}
        \end{subfigure}
        \begin{subfigure}{.32\textwidth}
            \includegraphics[width=\linewidth]{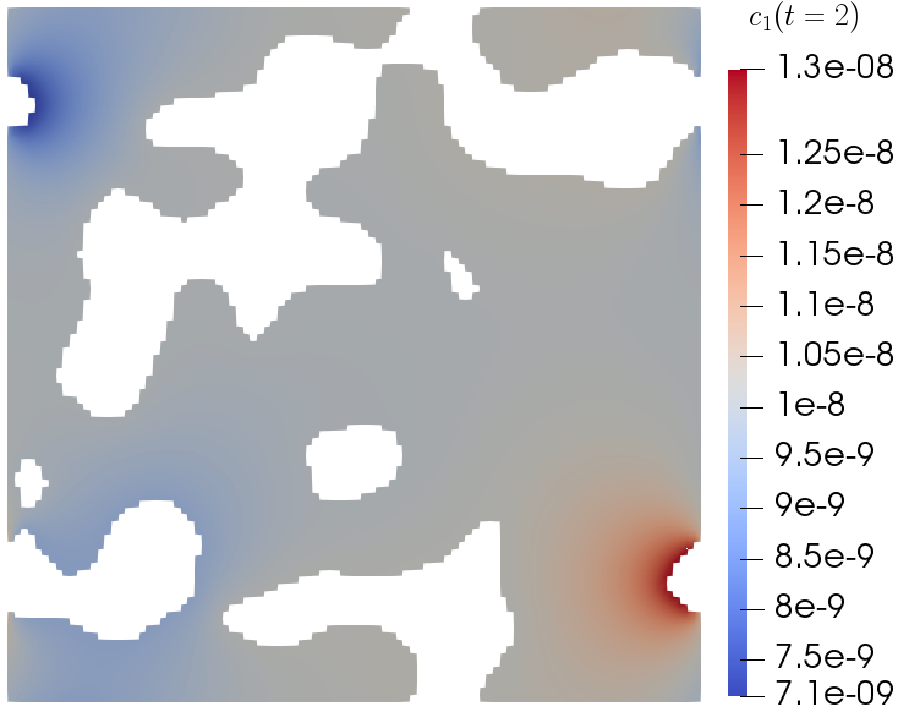}
        \end{subfigure}
        \begin{subfigure}{.32\textwidth}
            \includegraphics[width=\linewidth]{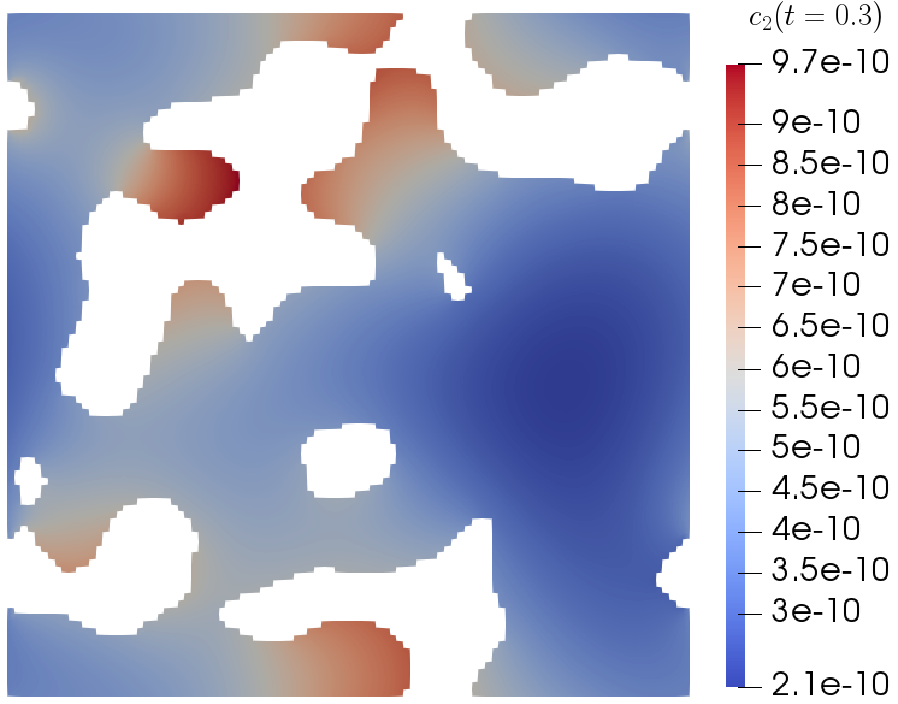}
        \end{subfigure}
        \begin{subfigure}{.32\textwidth}
            \includegraphics[width=\linewidth]{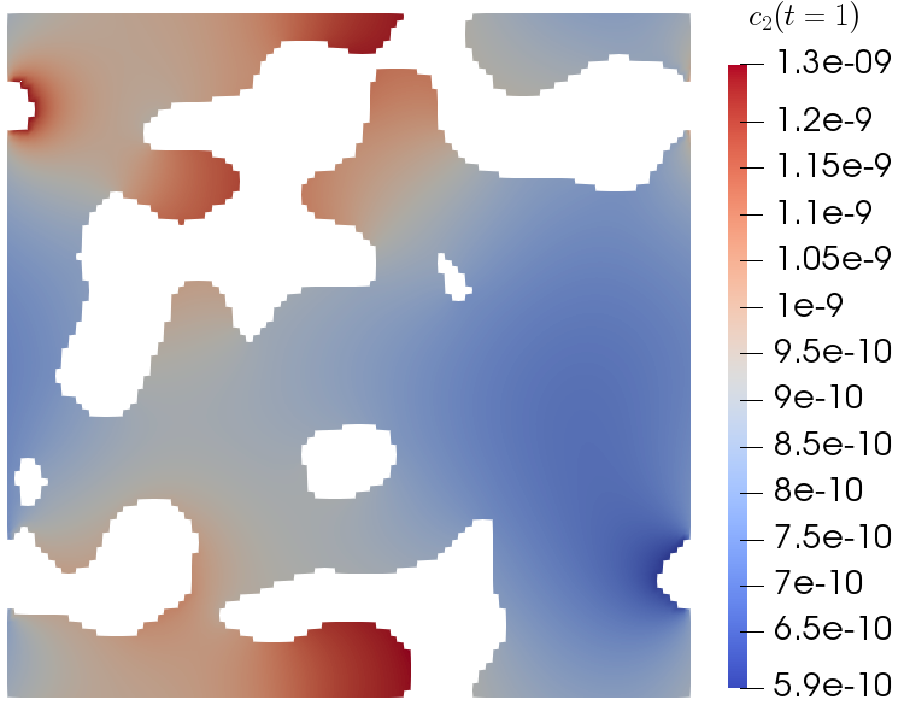}
        \end{subfigure}
        \begin{subfigure}{.32\textwidth}
            \includegraphics[width=\linewidth]{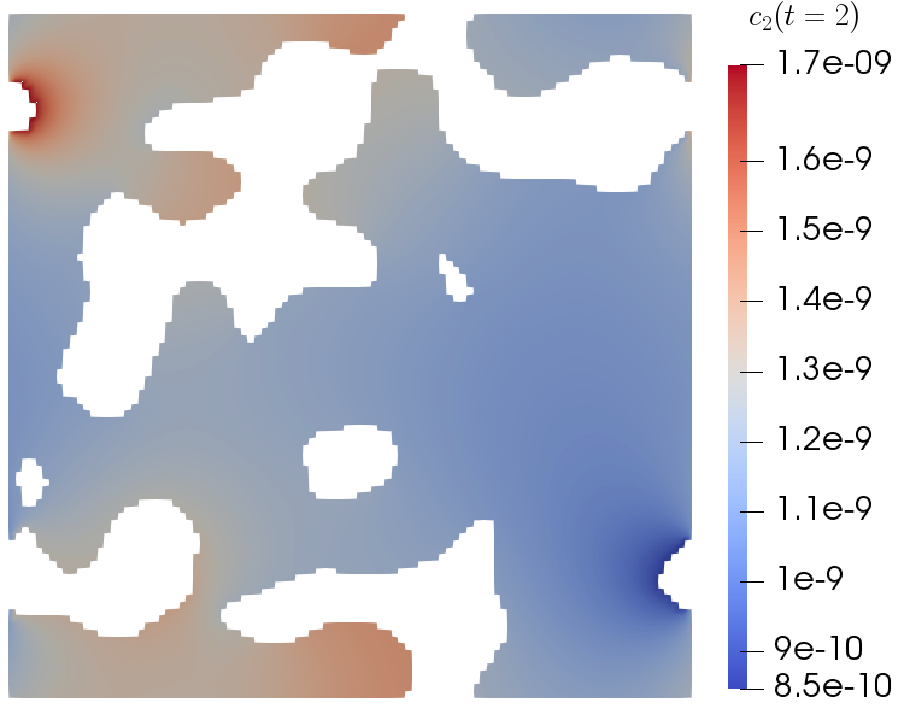}
        \end{subfigure}
        \begin{subfigure}{.32\textwidth}
            \includegraphics[width=\linewidth]{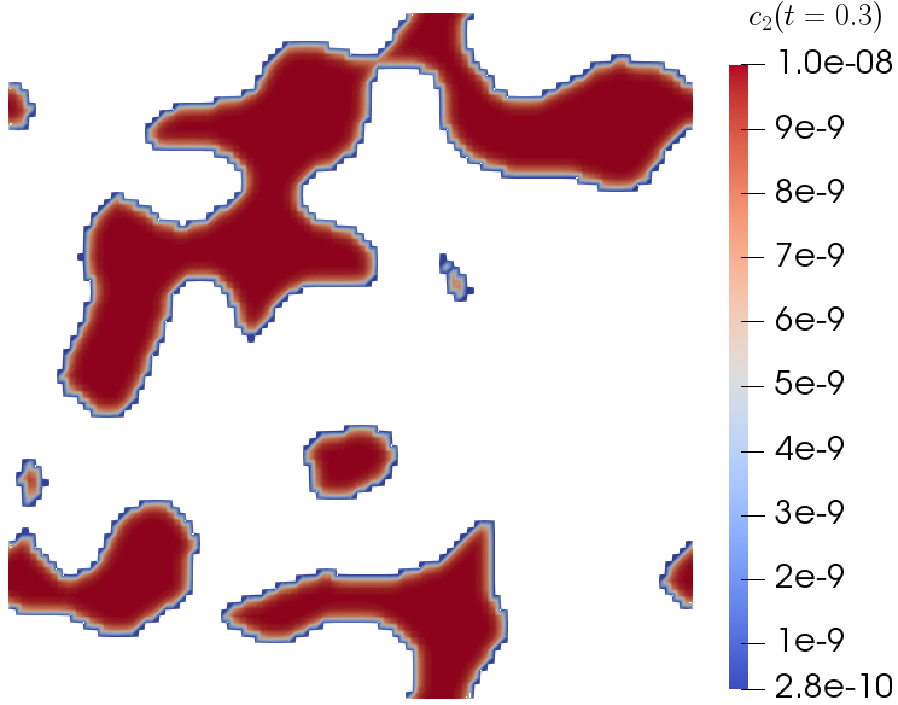}
        \end{subfigure}
        \begin{subfigure}{.32\textwidth}
            \includegraphics[width=\linewidth]{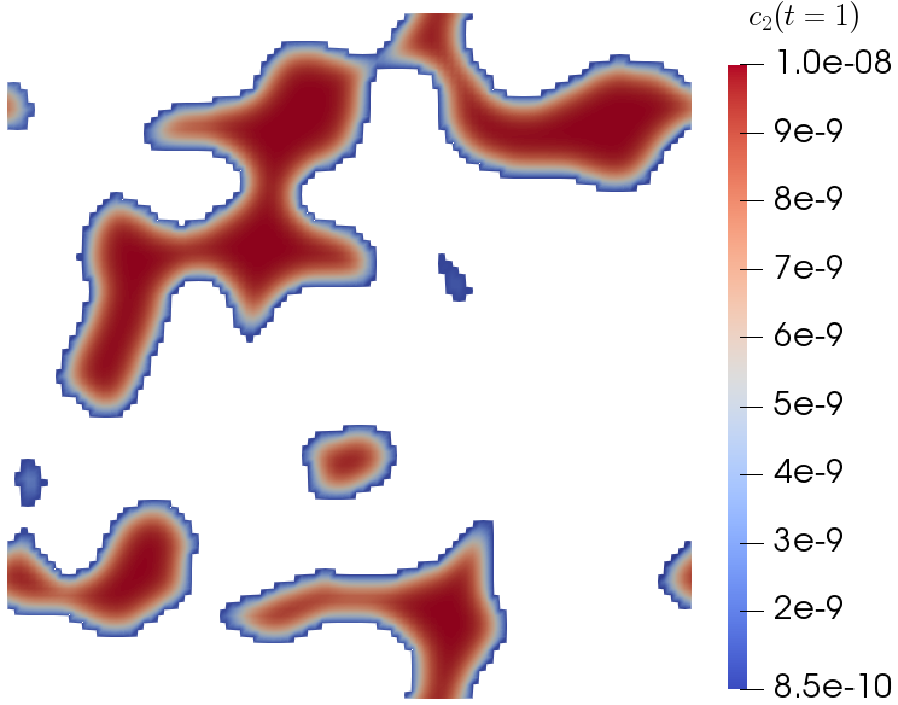}
        \end{subfigure}
        \begin{subfigure}{.32\textwidth}
            \includegraphics[width=\linewidth]{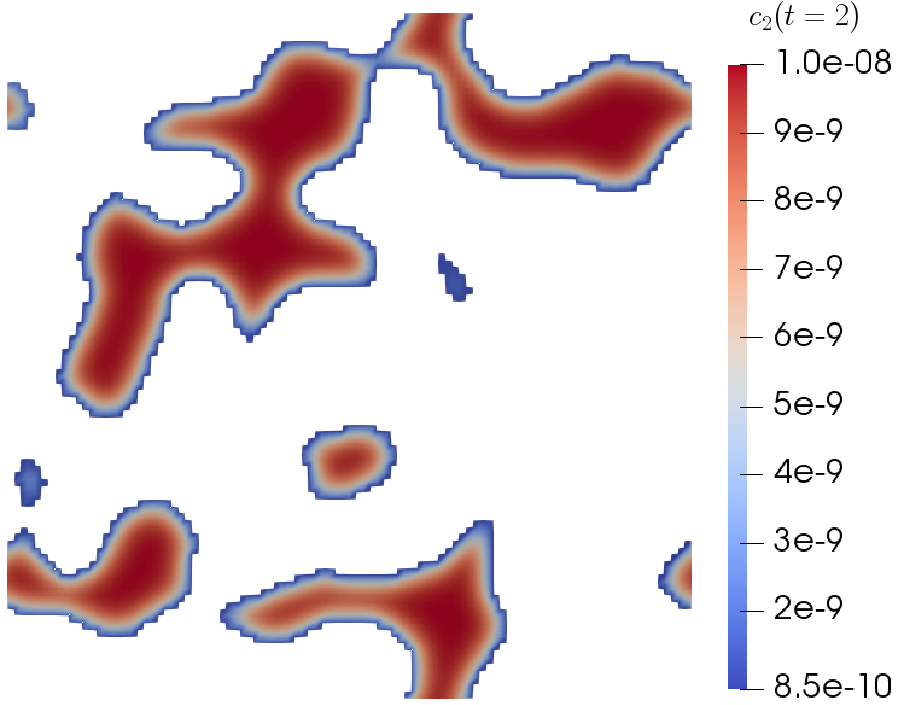}
        \end{subfigure}
        \begin{subfigure}{.32\textwidth}
            \includegraphics[width=\linewidth]{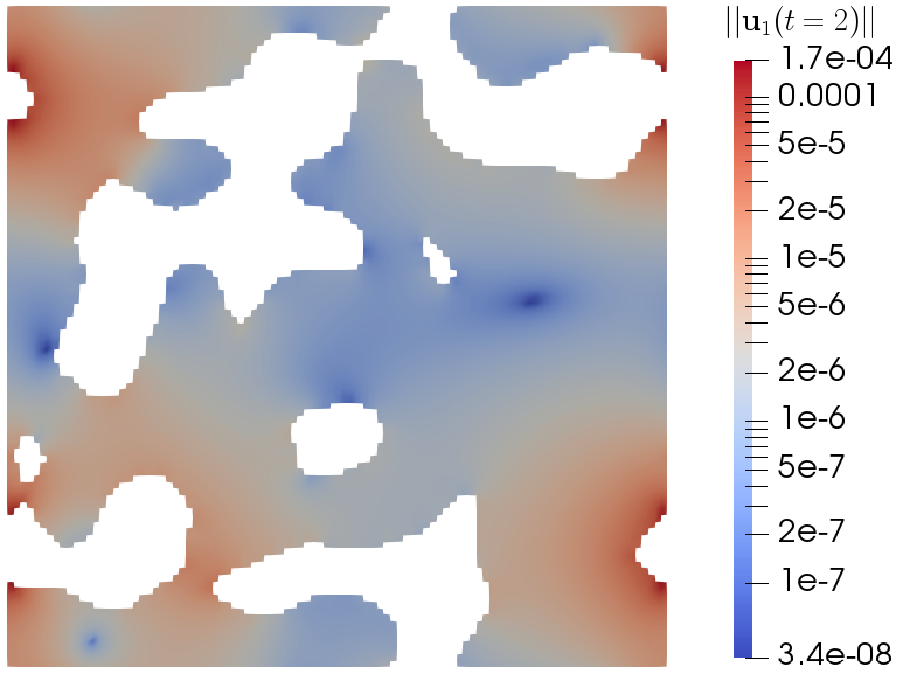}
        \end{subfigure}
        \begin{subfigure}{.32\textwidth}
            \includegraphics[width=\linewidth]{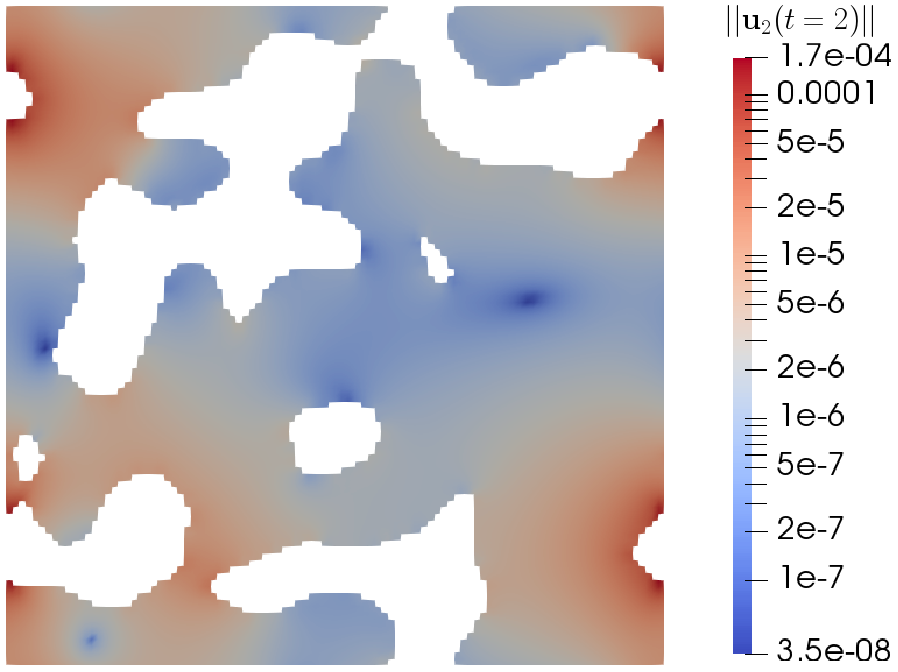}
        \end{subfigure}
        \begin{subfigure}{.32\textwidth}
            \includegraphics[width=\linewidth]{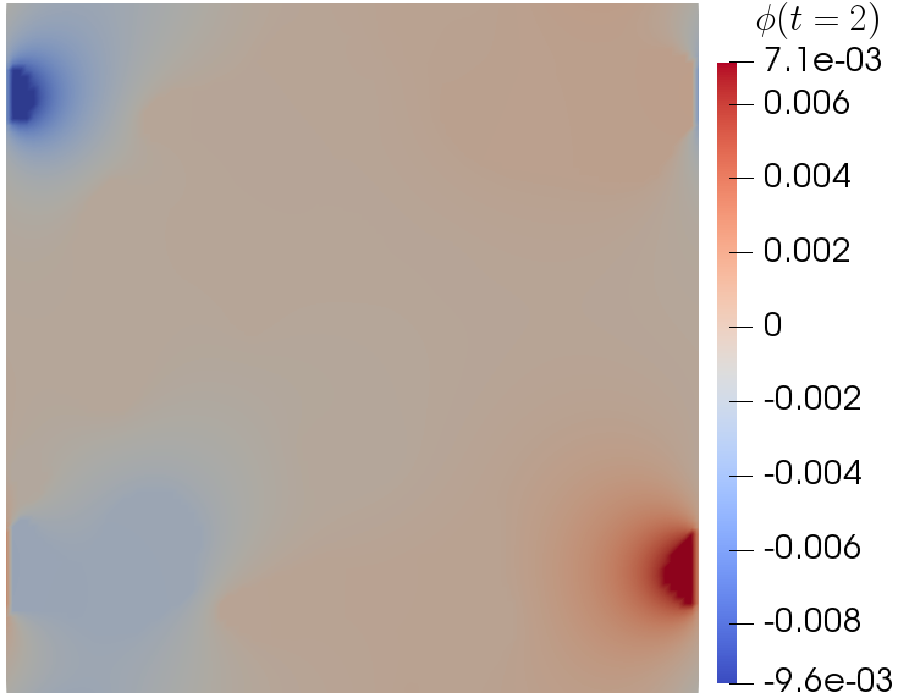}
        \end{subfigure}
        \caption{Results of the second random solid-fluid geometry, discretized using $N=1\times10^{4}$ cells. Top two rows: {Profiles of $c_1$ and $c_2$ at $t=0.3$, $1$, $2$s within $\Omega_{\text{f}}$.} Third row: {Profiles of $c_2$ at $t=0.3$, $1$, $2$s within $\Omega_{\text{s}}$.} Bottom row: {Profiles of $||\mathbf{u}_1||$, $||\mathbf{u}_2||$ and $\phi$ at $t=2$s.}}
        \label{fig::seed1_100x100Cells}
    \end{figure}

\begin{figure}[htbp]
        \centering
        \begin{subfigure}{.32\textwidth}
            \includegraphics[width=\linewidth]{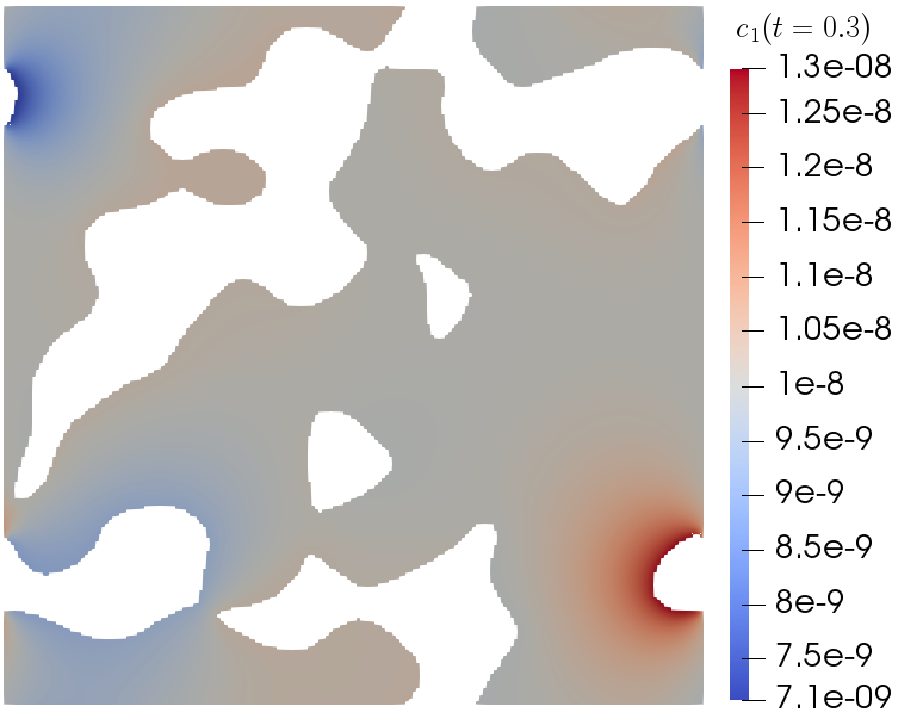}
        \end{subfigure}
        \begin{subfigure}{.32\textwidth}
            \includegraphics[width=\linewidth]{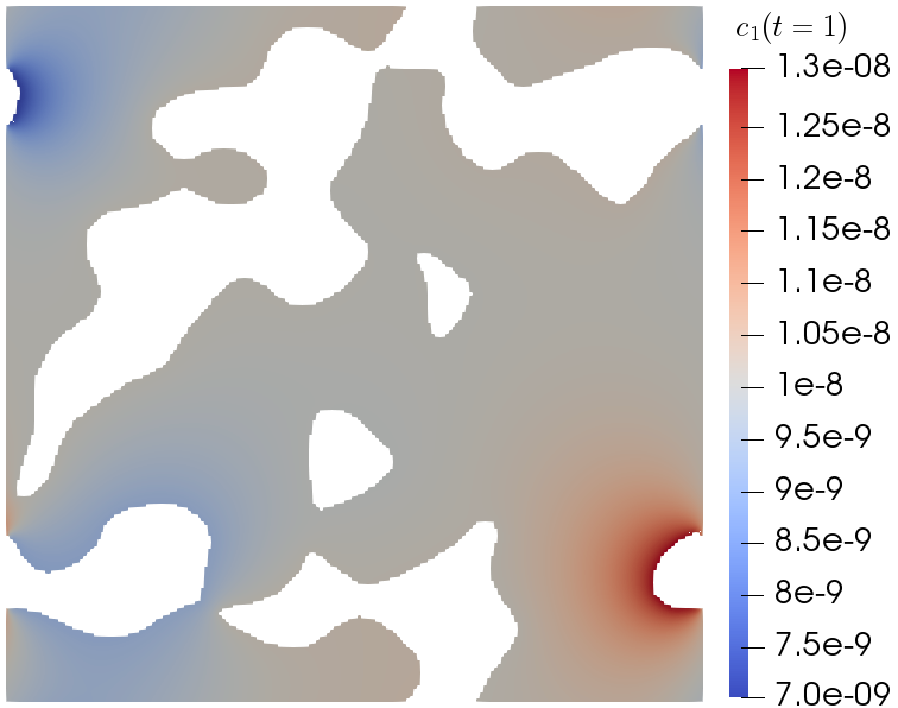}
        \end{subfigure}
        \begin{subfigure}{.32\textwidth}
            \includegraphics[width=\linewidth]{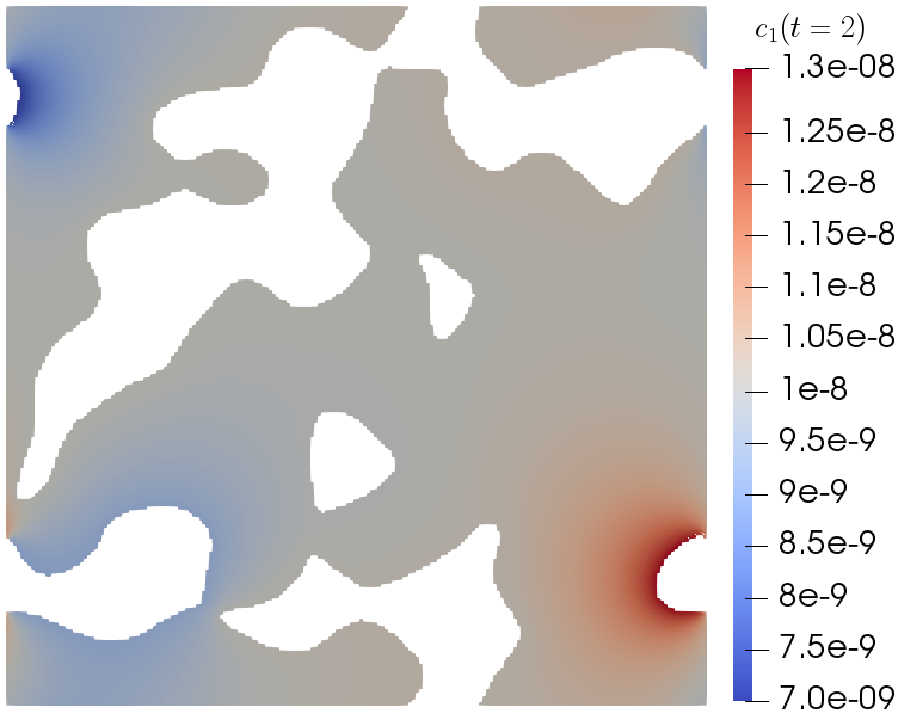}
        \end{subfigure}
        \begin{subfigure}{.32\textwidth}
            \includegraphics[width=\linewidth]{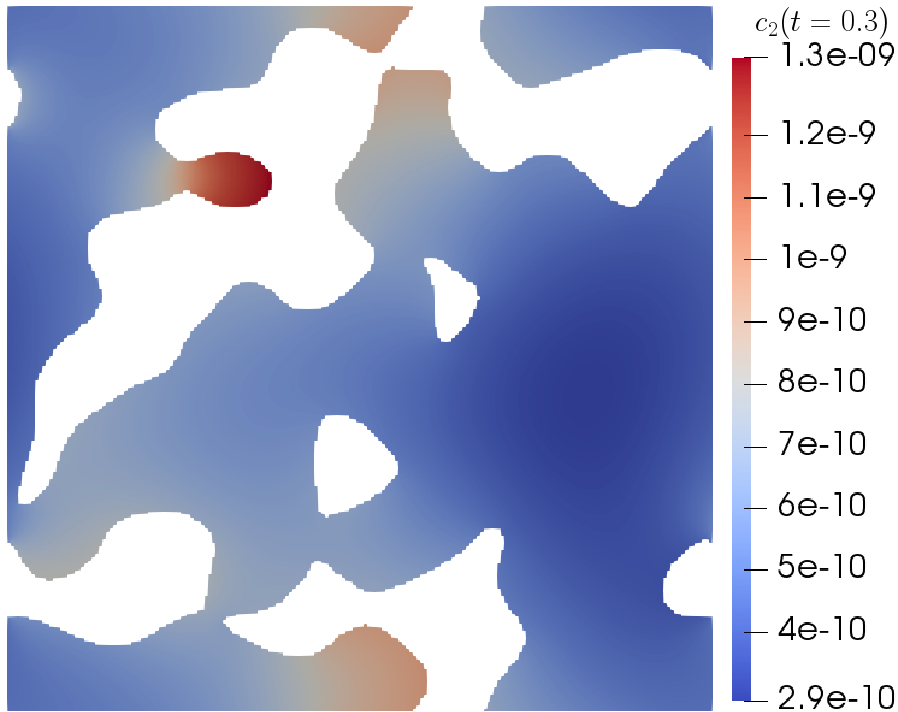}
        \end{subfigure}
        \begin{subfigure}{.32\textwidth}
            \includegraphics[width=\linewidth]{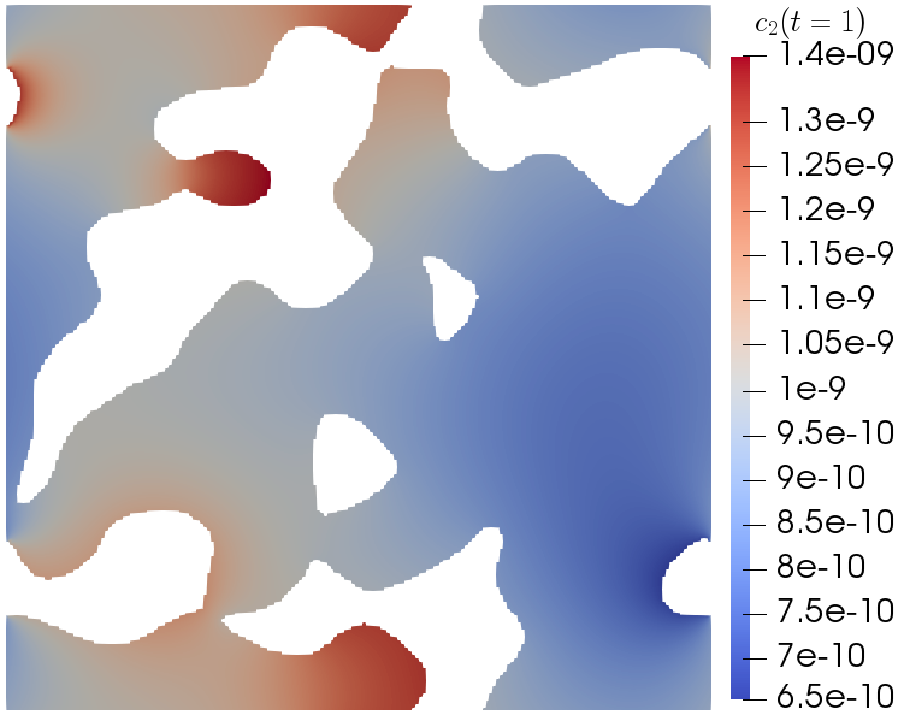}
        \end{subfigure}
        \begin{subfigure}{.32\textwidth}
            \includegraphics[width=\linewidth]{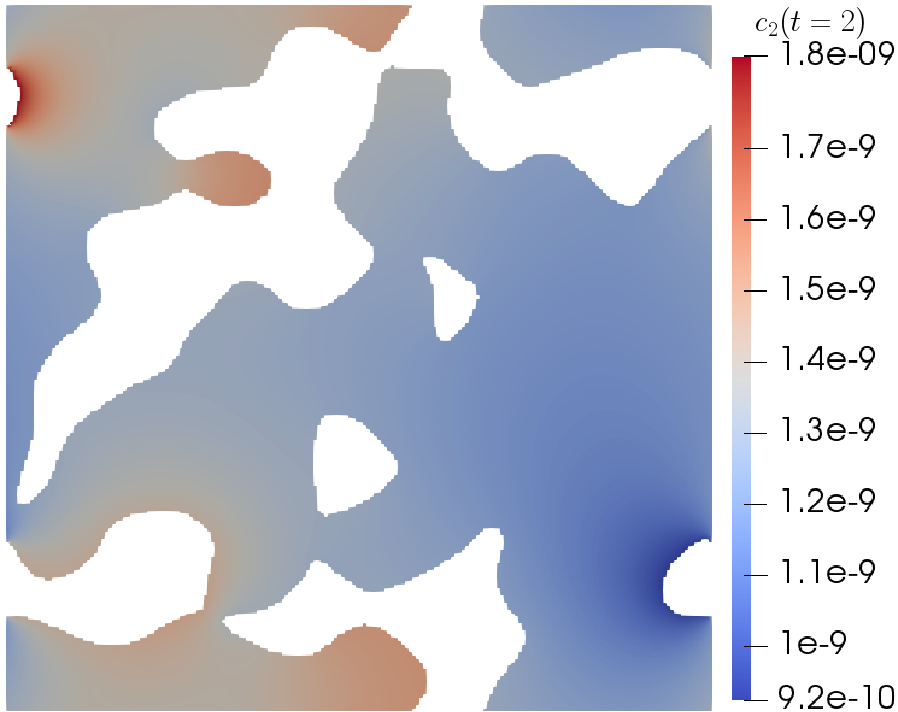}
        \end{subfigure}
        \begin{subfigure}{.32\textwidth}
            \includegraphics[width=\linewidth]{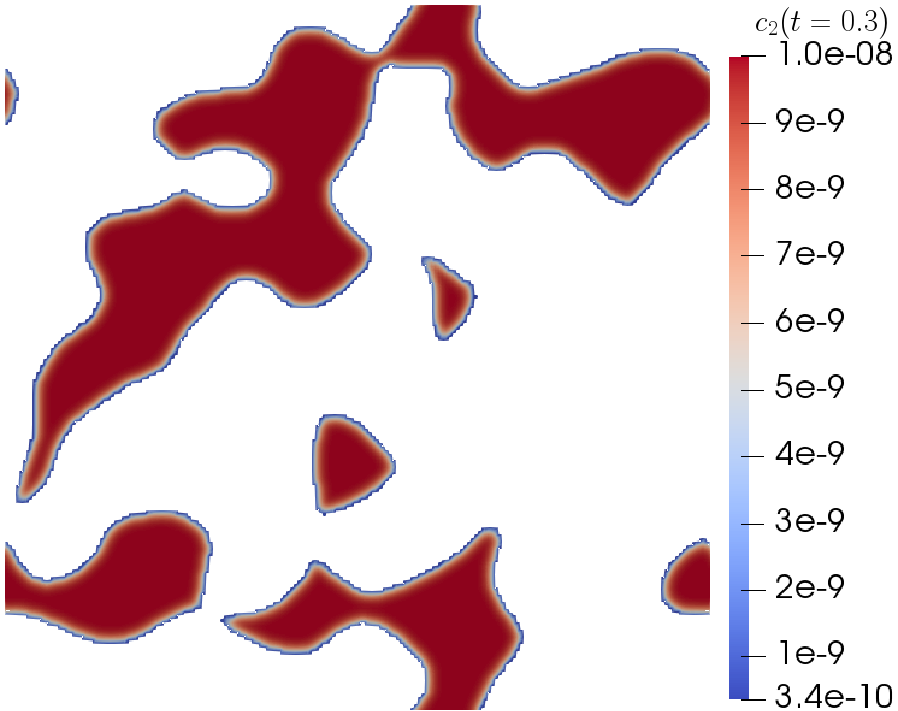}
        \end{subfigure}
        \begin{subfigure}{.32\textwidth}
            \includegraphics[width=\linewidth]{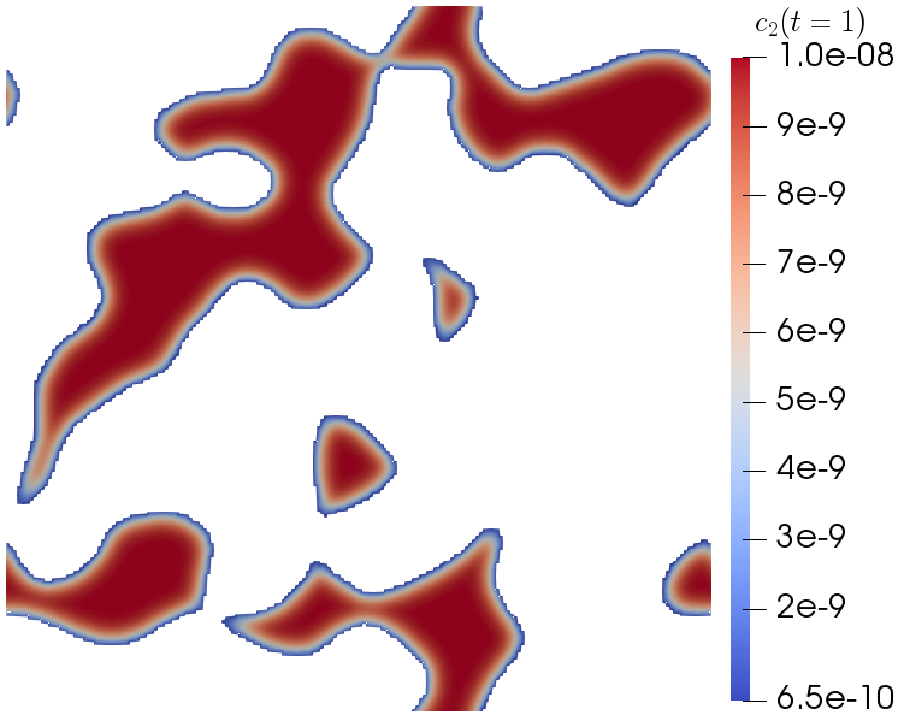}
        \end{subfigure}
        \begin{subfigure}{.32\textwidth}
            \includegraphics[width=\linewidth]{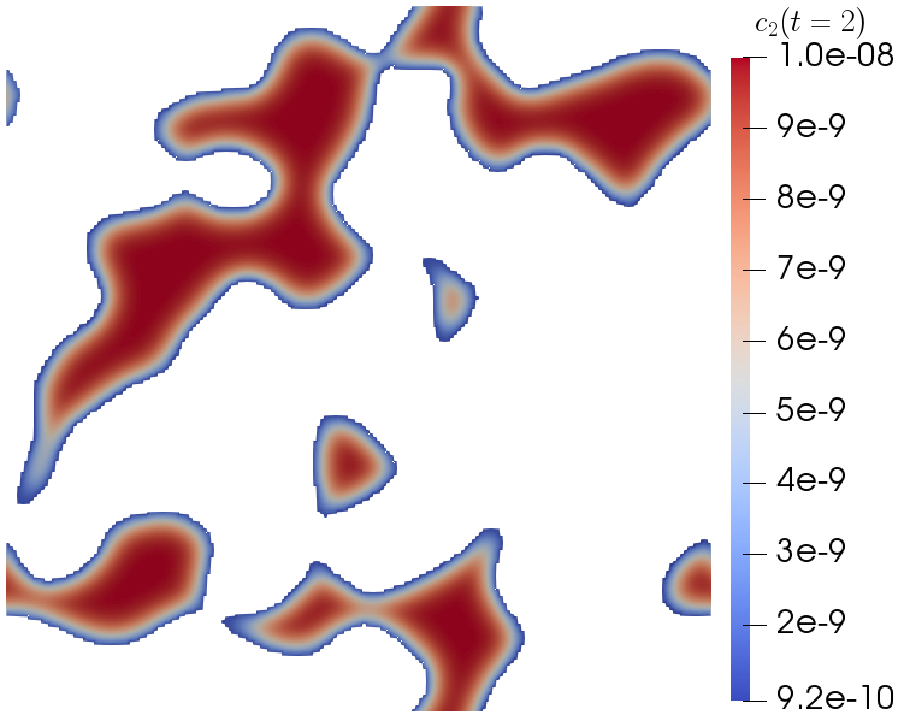}
        \end{subfigure}
        \begin{subfigure}{.32\textwidth}
            \includegraphics[width=\linewidth]{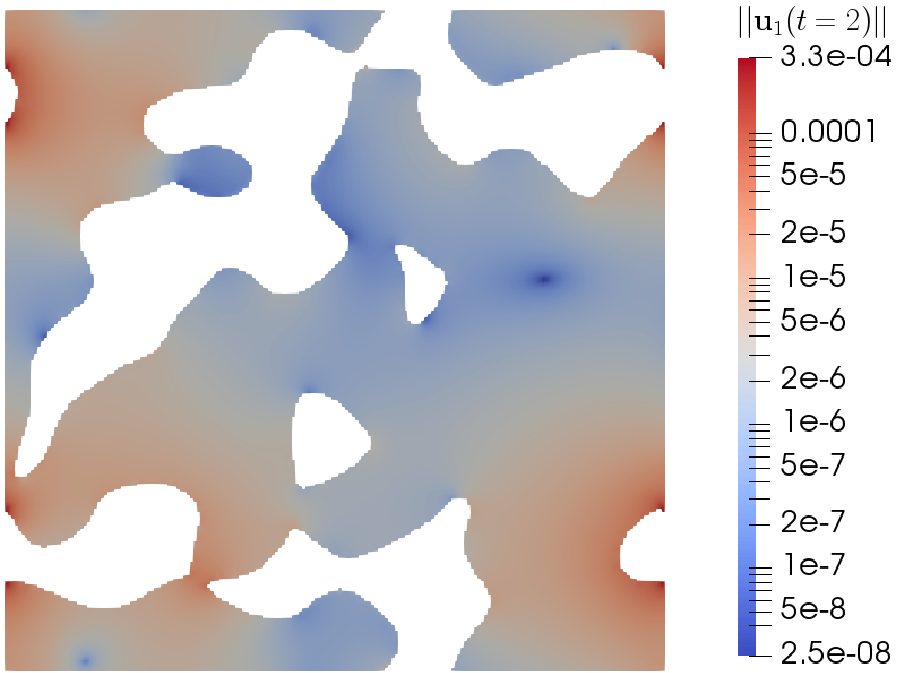}
        \end{subfigure}
        \begin{subfigure}{.32\textwidth}
            \includegraphics[width=\linewidth]{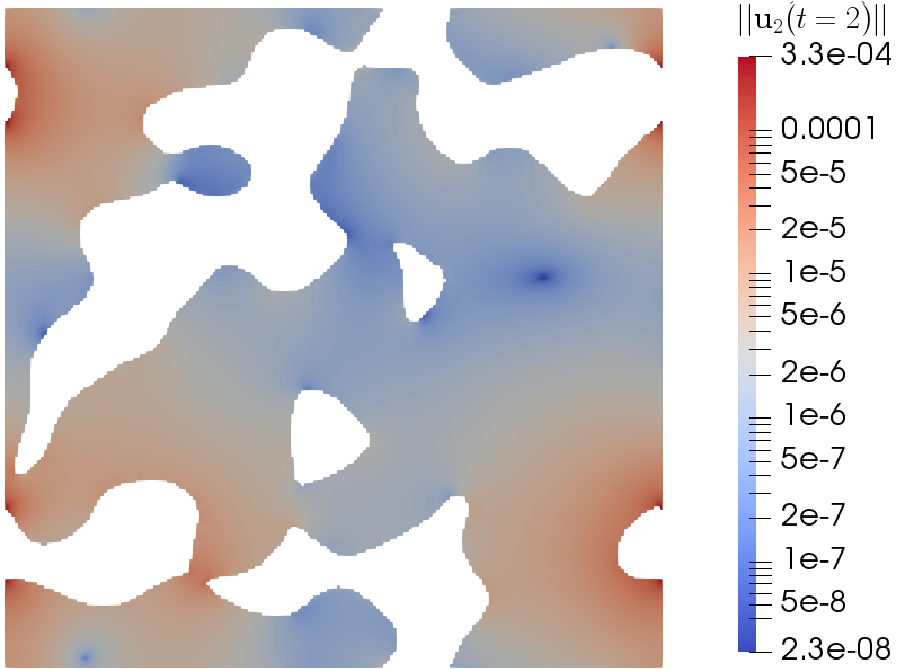}
        \end{subfigure}
        \begin{subfigure}{.32\textwidth}
            \includegraphics[width=\linewidth]{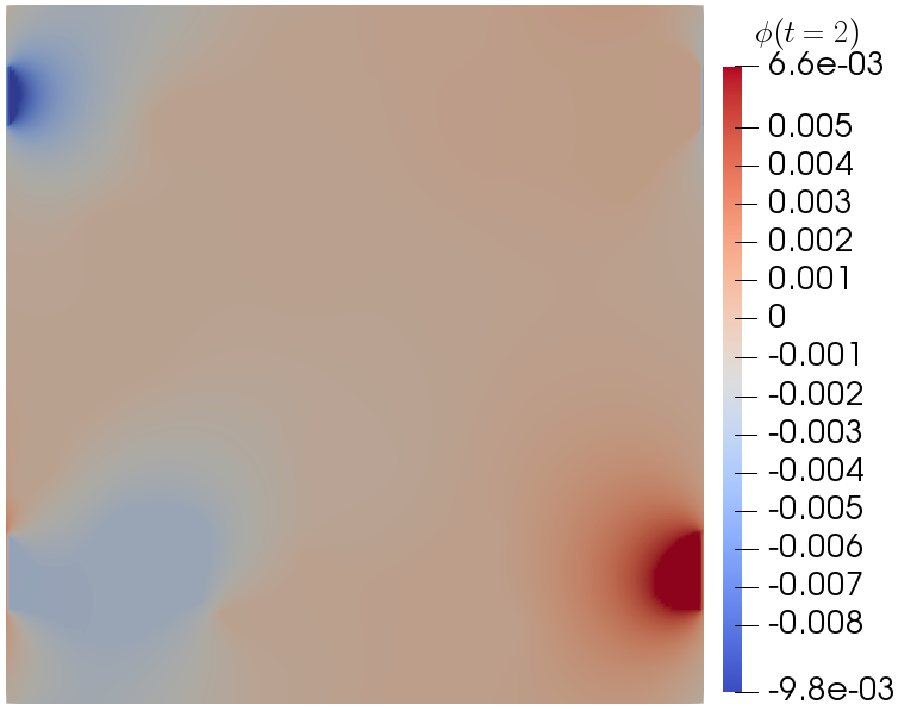}
        \end{subfigure}
        \caption{Results of the second random solid-fluid geometry, discretized using $N=4\times10^{4}$ cells. Top two rows: {Profiles of $c_1$ and $c_2$ at $t=0.3$, $1$, $2$s within $\Omega_{\text{f}}$.} Third row: {Profiles of $c_2$ at $t=0.3$, $1$, $2$s within $\Omega_{\text{s}}$.} Bottom row: {Profiles of $||\mathbf{u}_1||$, $||\mathbf{u}_2||$ and $\phi$ at $t=2$s.}}
        \label{fig::seed1_200x200Cells}
    \end{figure}

    We can compare these results with coarse mesh ($N=1\times 10 ^{4}$) results as well as with other random generation seeds. Results of the first realization with a coarse mesh, and of a coarse and fine simulation of a second random realization are reported respectively in \cref{fig::seed2_100x100Cells}, \cref{fig::seed1_100x100Cells}, and \cref{fig::seed1_200x200Cells}.
    
\section{Conclusions}\label{sec::Conclusion}

    In this paper, we present the development of open-source  solvers in \of capable of simulating microscopic electrokinetic flows of dilute ionic fluids. The underlying system, known as Stokes-Poisson-Nernst-Planck, is thoroughly reviewed where we discuss the assumptions on the fluid properties and limitations these bring. We later apply dimensional analysis to characterize the effects of dominating physical forces. This analysis is later used to give a better understanding of a common model reduction, known as the electro-neutrality assumption, as a result of such analysis.
    
    Many real-world applications of these flows involve some form of reactions between involved ionic species at the interface between fluid and solid. As such, we formulate a fully general reaction model capable of describing said heterogeneous reactions as a balance of fluxes across all reacting species.
    
    To verify flow descriptions obtained from our solvers are accurate, we compare against highly accurate solutions obtained, for simplified cases, with the  {Matlab toolbox \textit{Chebfun}} for a number of cases. With each case we find good agreement, showing spatial grid convergence of the results. We later use our solvers across a randomly generated solid-fluid porous cell, similar to a representative elementary volume (REV) used in homogenization theory, to show the solver's ability in handling complex microstructures. 

    In future works, we plan on constructing fully physically realistic cases, such as modelling the electrokinetic flow through a porous battery half-cell. Furthermore, we later plan to apply uncertainty quantification in conjunction with these randomly generated porous cells, quantifying the effects on the flow for varying geometrical properties.

\section*{Declarations}
The authors declare that they have no conflicts of interest.
MI and RB have been supported by the University of Nottingham Impact Acceleration Award. 
    
\appendix

 \section{Quasi-coupled Newton iterative effective reactive Robin conditions}\label{app:int}

\subsection{General non-linear reaction rate for restricted binary reaction}\label{sec::NonLinearReactionEffRobinBCs}
    Let there be the following reversible binary reaction between species $c_{\text{s}}$ and $c_{\text{f}}$ along the interface $\Gamma$ shared by $\Omega_{\text{s}}$ and $\Omega_{\text{f}}$,
        \begin{equation}
            \ce{\nu_{\text{s}}X_{\text{s}} <=> \nu_{\text{f}}X_{\text{f}}},
        \end{equation}
    and consider each species restricted to its respective sub-domains:
        \begin{equation}
            c_{\text{s}} = \begin{cases}
                0 & \bx\in\Omega_{\text{f}} \\
                c_{\text{s}} & \bx\in\Omega_{\text{s}}
            \end{cases}, \qquad
            c_{\text{f}} = \begin{cases}
                c_{\text{f}} & \bx\in\Omega_{\text{f}} \\
                0 & \bx\in\Omega_{\text{s}}
            \end{cases}.
        \end{equation}
    In these circumstances, we have the following condition along $\Gamma$ with reaction rate $r$:
        \begin{equation}\label{eqn::NonLinearBinaryReactionRestrictedCondition}
            \flux_{\text{s}}\big\rvert_{\Gamma_{\text{s}}}\cdot\mathbf{n} = r(c_{\text{s}}\rvert_{\Gamma_{\text{s}}},c_{\text{f}}\rvert_{\Gamma_{\text{f}}}) = -\flux_{\text{f}}\big\rvert_{\Gamma_{\text{f}}}\cdot\mathbf{n^\prime}.
        \end{equation}
    We denote here $\mathbf{n}$ to be the unit normal along $\Gamma$ facing outwards of $\Omega_{\text{s}}$ and $\mathbf{n^\prime} = -\mathbf{n}$. Consider the reaction rate $r$ is a non-linear function of $c_{\text{s}}$ and $c_{\text{f}}$. To obtain effective Robin conditions of \cref{eqn::NonLinearBinaryReactionRestrictedCondition} for $c_{\text{s}}$ and $c_{\text{f}}$ we employ the Newton-Raphson method. First, we take the vector form of \cref{eqn::NonLinearBinaryReactionRestrictedCondition} with all terms moved to the left. To simplify writing, for now we omit what side of $\Gamma$ each term is evaluated on:
        \begin{equation}\label{bcon:fluxes_Non_Linear_Reaction_Vec_F=0}
    		\mathbf{J}(c_{\text{s}},c_{\text{f}}) =
    		    \begin{pmatrix}
    		        J_1 \\
    		        J_2
    		    \end{pmatrix}
    		\equiv
    			\begin{pmatrix}
    				-D_{\text{s}}\bm{\nabla}_{\mathbf{n}}c_{\text{s}} -  r(c_{\text{s}},c_{\text{f}})  \\[6pt]
    				(\bu _{\text{f}}\cdot\mathbf{n}^\prime)c_{\text{f}} - D_{\text{f}}\bm{\nabla}_{\mathbf{n}^\prime}c_{\text{f}} - \frac{z_{\text{f}}D_{\text{f}}F}{RT}c_{\text{f}}\bm{\nabla}_{\mathbf{n}^\prime}\phi_{\text{f}} +  r(c_{\text{s}},c_{\text{f}}) 
    			\end{pmatrix}
    		= \mathbf{0}.
    	\end{equation}
    Assume we have a current approximate solution $\mathbf{c}^{N} = (c_{\text{s}}^N \hspace{.5em}  c_{\text{f}}^N)^\top$. Linearizing $\mathbf{J}$ around $\mathbf{c}^N$ then gives
    	\begin{equation}\label{bcon:Non_Linear_Reaction_Vec_F_Tay_Exp}
    		\mathbf{J}(\mathbf{c}^N) +\pardiff{\mathbf{J}}{\mathbf{c}}\bigg\rvert_{\mathbf{c}^N}\left(\mathbf{c}-\mathbf{c}^N\right) = \mathbf{0},
    	\end{equation}
    whereby setting $\mathbf{c}=\mathbf{c}^{N+1}$ is the next iterative solution. We denote here $\frac{\partial\mathbf{J}}{\partial\mathbf{c}}$ to be the Jacobian matrix of $\mathbf{J}$. To determine the Jacobian we define the directional derivative for some function $G(c)$ along the direction $\delta c = c - c^N$ as:
    	\begin{equation}\label{def:directional_derivative}
    			\pardiff{G}{c}(\delta c) = \lim_{\varepsilon \rightarrow 0}\frac{G(c + \varepsilon\delta c) - G(c)}{\varepsilon}.
    	\end{equation}
    Using this definition of the directional derivative we can then determine the Jacobian within \cref{bcon:Non_Linear_Reaction_Vec_F_Tay_Exp}. Since the method is similar for all components of $\mathbf{J}$, we will focus on the first element $J_1$. Taking the directional derivative of $J_1$ with respect to $c_{\text{s}}$ along $\delta c_{\text{s}}$ at $\mathbf{c}=\mathbf{c}^N$ gives:
    	\begin{equation}\label{eqn:deriv_J_1_wrt_c_s}
    		\pardiff{J_1}{c_{\text{s}}}\bigg\rvert_{\mathbf{c}^N}(\delta c_{\text{s}}) = \lim_{\varepsilon \rightarrow 0 }\frac{1}{\varepsilon}\left(-D_{\text{s}}\bm{\nabla}_{\mathbf{n}}(c_{\text{s}}^N + \varepsilon\delta c_{\text{s}})
    		-  r (c_{\text{s}}^N + \varepsilon\delta c_{\text{s}},c_{\text{f}}^N) + D_{\text{s}}\bm{\nabla}_{\mathbf{n}}c_{\text{s}}^N + r (c_{\text{s}}^N,c_{\text{f}}^N)\right).
    	\end{equation}
    Taking the following Taylor expansions about $\varepsilon=0$ of $r$ and the concentration gradient,
    
    	\begin{align}
    		\bm{\nabla}_{\mathbf{n}}(c_{\text{s}}^N + \varepsilon\delta c_{\text{s}}) &= \bm{\nabla}_{\mathbf{n}}c_{\text{s}}^N + \varepsilon\bm{\nabla}_{\mathbf{n}}\delta c_{\text{s}} + \mathcal{O}(\varepsilon^2), \\
    		 r (c_{\text{s}}^N + \varepsilon\delta c_{\text{s}},c_{\text{f}}^N) &=  r (c_{\text{s}}^N,c_{\text{f}}^N) + \varepsilon\delta c_{\text{s}} \pardiff{r}{c_{\text{s}}}\bigg\rvert_{\mathbf{c}^N} + \mathcal{O}(\varepsilon^2) ,\label{eqn:Tay_Exp_hat_in_epsilon}
    	\end{align}
    then substituting into \cref{eqn:deriv_J_1_wrt_c_s} we obtain:
    	\begin{equation}
    		\pardiff{J_1}{c_{\text{s}}}\bigg\rvert_{\mathbf{c}^N}(\delta c_{\text{s}}) = -D_{\text{s}}\bm{\nabla}_{\mathbf{n}}\delta c_{\text{s}}  - \delta c_{\text{s}} \pardiff{r}{c_{\text{s}}}\bigg\rvert_{\mathbf{c}^N}.
    	\end{equation}
    The same can be done to find the directional derivative of $J_1$ with respect to $c_{\text{f}}$, giving:
    	\begin{equation}
    		\pardiff{J_1}{c_{\text{f}}}\bigg\rvert_{\mathbf{c}^N}(\delta c_{\text{f}}) = -\delta c_{\text{f}} \partialderiv{ r }{c_{\text{f}}}\bigg\rvert_{\mathbf{c}^N}.
    	\end{equation}
    Repeating the same process for $J_2$ of \cref{bcon:fluxes_Non_Linear_Reaction_Vec_F=0} and collecting all elements together and evaluating at $\mathbf{c}^N$ gives
    	\begin{equation}\label{eqn::DirectionalJacobian}
    		\pardiff{\mathbf{J}}{\mathbf{c}}\bigg\rvert_{\mathbf{c}^N}\left(\mathbf{\delta c}\right) = 
    			\begin{pmatrix}
    				-D_{\text{s}}\bm{\nabla}_{\mathbf{n}}\delta c_{\text{s}}  - \delta c_{\text{s}}\partialderiv{ r }{c_{\text{s}}} - \delta c_{\text{f}}\pardiff{ r }{c_{\text{f}}} \\[6pt]
    				\delta c_{\text{s}} \partialderiv{ r }{c_{\text{s}}} + ( \bu _{\text{f}}\cdot\mathbf{n}^\prime)\delta c_{\text{f}} - \frac{z_{\text{f}}D_{\text{f}}F}{RT}\delta c_{\text{f}}\bm{\nabla}_{\mathbf{n}}\phi_{\text{f}} + \delta c_{\text{f}} \pardiff{ r }{c_{\text{f}}}
    			\end{pmatrix},
    	\end{equation}
    where derivatives in $r$ are evaluated at $\mathbf{c}^N$. We may then substitute \cref{eqn::DirectionalJacobian} into \cref{bcon:Non_Linear_Reaction_Vec_F_Tay_Exp} to get
    	\begin{equation}\label{bcon:Non_Lin_Reac_Linearized_Form_Vector}
    		\begin{pmatrix}
    			-D_{\text{s}}\bm{\nabla}_{\mathbf{n}}c_{\text{s}} -  r   - \delta c_{\text{f}} \partialderiv{ r }{c_{\text{f}}}   - \delta c_{\text{s}}\partialderiv{ r }{c_{\text{s}}}  \\[6pt]
    			( \bu_{\text{f}}\cdot\mathbf{n}^\prime)c_{\text{f}} - D_{\text{f}}\bm{\nabla}_{\mathbf{n}^\prime}c_{\text{f}} - \frac{z_{\text{f}}D_{\text{f}}F}{RT}c_{\text{f}}\bm{\nabla}_{\mathbf{n}^\prime}\phi_{\text{f}} +  r   + \delta c_{\text{s}}\partialderiv{ r }{{c}_{\text{s}}} + \delta c_{\text{f,f}} \partialderiv{ r }{c_{\text{f}}}
    		\end{pmatrix}
    		= \mathbf{0},
    	\end{equation}
    where $r$ and its derivatives are evaluated at $\mathbf{c}^N$. Re-arranging the first element of \cref{bcon:Non_Lin_Reac_Linearized_Form_Vector}, using the first order approximation $\bm{\nabla}_{\mathbf{n}}c_{\text{s}} = B_{\text{s}}(c_{\text{s}} - c_{\text{s,c}})$ where $c_{\text{s,c}}$ is the value at the cell centre of cells sharing a face with $\Gamma$, and substituting into the second element then gives:
    	\begin{equation}
    		( \bu_{\text{f}}\cdot\mathbf{n}^\prime)c_{\text{f}} - D_{\text{f}}\bm{\nabla}_{\mathbf{n}^\prime}c_{\text{f}} - \frac{z_{\text{f}}D_{\text{f}}F}{RT}c_{\text{f}}\bm{\nabla}_{\mathbf{n}^\prime}\phi_{\text{f}} =
    		 - r - \delta c_{\text{f}} \partialderiv{ r }{c_{\text{f}}}  + c_{\text{s}}^N\partialderiv{ r }{c_{\text{s}}} + P(c_{\text{f}}),
    	\end{equation}
	where
	    \begin{equation}
	        P(c_{\text{f}}) = \partialderiv{ r }{c_{\text{s}}} \left(\frac{-D_{\text{s}}B_{\text{s}}c_{\text{s,c}} +  r  + \delta c_{\text{f}}\partialderiv{ r }{c_{\text{f}}}  - c_{\text{s}}^N\partialderiv{ r }{c_{\text{s}}} }{D_{\text{s}}B_{\text{s}} + \partialderiv{ r }{c_{\text{s}}} }\right).
	    \end{equation}
    Setting $\mathbf{c} = \mathbf{c}^{N+1}$ to indicate the next iterative solution results in an iterative Robin condition of the form $-D_{\text{f}}^*\bm{\nabla}_{\mathbf{n}^\prime}c_{\text{f}}^{N+1} = -K_{\text{f}}^*c_{\text{f}}^{N+1} -F_{\text{f}}^*$, with effective coefficients dependent on the previous iteration $\mathbf{c}^N$:
    	\begin{align}
    		D_{\text{f}}^* &= D_{\text{f}} ,\label{bcon:D_f_Eff_GenReacBinary}\\
    		K_{\text{f}}^*	&=-\frac{\partialderiv{ r }{c_{\text{s}}} \partialderiv{ r }{c_{\text{f}}} }{D_{\text{s}}B_{\text{s}} +\partialderiv{ r }{c_{\text{s}}} } + \partialderiv{ r }{c_{\text{f}}}  - \frac{z_{\text{f}}D_{\text{f}}F}{RT}\bm{\nabla}_{\mathbf{n}^\prime}\phi_{\text{f}} + ( \bu _{\text{f}}\cdot\mathbf{n}^\prime), \label{bcon:K_f_Eff_GenReacBinary}\\
    		F_{\text{f}}^*	&=  r   - c_{\text{f}}^N\partialderiv{ r }{c_{\text{f}}}  - c_{\text{s}}^N\partialderiv{ r }{c_{\text{s}}}  - \partialderiv{ r }{c_{\text{s}}} \left(\frac{-D_{\text{s}}B_{\text{s}}c_{\text{s,c}} +  r   - c_{\text{f}}^N\partialderiv{ r }{c_{\text{f}}} -c_{\text{s}}^N\partialderiv{ r }{c_{\text{s}}} }{D_{\text{s}}B_{\text{s}} +\partialderiv{ r }{c_{\text{s}}} }\right).\label{bcon:F_f_Eff_GenReacBinary}
    	\end{align}
    Remember, $r$ and its derivatives are evaluated at the previous iteration $\mathbf{c}^N$. Similar steps can also be applied to the accompanying iterative Robin condition for $c_{\text{s}}$ of the form $-D_{\text{s}}^*\bm{\nabla}_{\mathbf{n}}c_{\text{s}}^{N+1} = -K_{\text{s}}^*c_{\text{s}}^{N+1} -F_{\text{s}}^*$, with coefficients:
    	\begin{align}
    		D_{\text{s}}^*	&= D_{\text{s}}, \label{bcon:D_s_Eff_General_Reac}\\
    		K_{\text{s}}^*	&= \frac{\partialderiv{ r }{c_{\text{s}}} \partialderiv{ r }{c_{\text{f}}} }{( \bu_{\text{f}}\cdot\mathbf{n}^\prime) - D_{\text{f}}B_{\text{f}} - \frac{z_{\text{f}}D_{\text{f}}F}{RT}\bm{\nabla}_{\mathbf{n}^\prime}\phi_{\text{f}} + \partialderiv{ r }{c_{\text{f}}} } - \partialderiv{ r }{c_{\text{s}}} ,  \\
    		F_{\text{s}}^*	&= - r   +c_{\text{f}}^N\partialderiv{ r }{c_{\text{f}}}  + c_{\text{s}}^N\partialderiv{ r }{c_{\text{s}}}
    		- \partialderiv{ r }{c_{\text{f}}} \left(\frac{-D_{\text{f}}B_{\text{f}}c_{\text{f,c}} -  r  + c_{\text{f}}^N\partialderiv{ r }{c_{\text{f}}}  + c_{\text{s}}^N\partialderiv{ r }{c_{\text{s}}} }{( \bu _{\text{f}}\cdot\mathbf{n}^\prime) - D_{\text{f}}B_{\text{f}} - \frac{z_{\text{f}}D_{\text{f}}F}{RT}\bm{\nabla}_{\mathbf{n}^\prime}\phi_{\text{f}} + \partialderiv{ r }{c_{\text{f}}} }\right) .\label{bcon:F_s_Eff_General_Reac}
    	\end{align}
    

\subsection{General non-linear unrestricted reaction rate}
    Here we will formulate a set of iterative Robin conditions required to solve the reactive conditions for unrestricted species and non-linear reactions. Consider a general non-linear reaction with $i=J,...,K$ reactants and $i=K+1,...,M$ products:
    	\begin{equation}
    		\sum_{i=J}^{K}\nu_{i}B_{i} + ne^{-} =\sum_{i=K+1}^{M}\nu_{i}B_{i} \,.\label{def::unresGenReacEqn}
    	\end{equation}
    We model \cref{def::unresGenReacEqn} using the following conditions with general, possibly non-linear, reaction rates $r_i$:
    	\begin{align}
    		\big[\flux_i \cdot\mathbf{n}\big]^{\Gamma_{\text{s}}}_{\Gamma_{\text{f}}} &= r_{i}(c_J,c_{J+1},\dots,c_{M}) \qquad &\bx\in\Gamma\text{ for }i=J,...,M ,\label{bcon::unresGenReac}\\
    		\big[c_{i}\big]^{\Gamma_{\text{s}}}_{\Gamma_{\text{f}}} &=0 \qquad &\bx\in\Gamma\text{ for }i=J,...,M .\label{bcon::unresGenReacContinuity}
    	\end{align}
    We begin by writing $\flux_i$ in full and use $\mathbf{n}=-\mathbf{n}^\prime$ to give $\flux_i\rvert_{\Gamma_{\text{f}}}\cdot\mathbf{n} = -\flux_{i}\rvert_{\Gamma_{\text{f}}}\cdot\mathbf{n}^\prime$, a more acceptable form for \of. Since formulating the necessary Robin conditions is the same for all species $i$, we omit the subscript from here onward unless necessary. Instead, we introduce two new subscripts $s$ and $f$ denoting evaluation on the solid or fluid side of $\Gamma$:
    	\begin{equation}
    		- D_{i}\nabla_{\mathbf{n}}c_{\text{s}} + ( \bu_{\text{f}}\cdot\mathbf{n}^\prime)c_{\text{f}} - D_i\left(\nabla_{\mathbf{n}^\prime}c_{\text{f}}+\frac{z_iF}{RT}\nabla_{\mathbf{n}^\prime}\phi_{\text{f}}\right) = r_{i}(c_J,c_{J+1},\dots,c_{M}) .\label{bcon::unresGenReacFull}
    	\end{equation}
    We next use the first order approximation $\nabla_{\mathbf{n}^\prime}c_{\text{f}} = B_{\text{f}}(c_{\text{f}} - c_{\text{f,c}})$, where the coefficient $B_{\text{f}}$ is determined by \of and $c_{\text{f,c}}$ denotes the cell centre value of cells sharing a face with $\Gamma$. Alongside the continuity condition \cref{bcon::unresGenReacContinuity} this allows us to rearrange \cref{bcon::unresGenReacFull} into a form close to Robin:
    	\begin{equation}\label{bcon::unresGenReacCsf}
    		-D_i\nabla_{\mathbf{n}}c_{s,f} = c_{s,f}\left[ - ( \bu _{\text{f,f}}\cdot\mathbf{n}^\prime) + D_iB_{\text{f}} + \frac{z_iD_iF}{RT_{\text{f,f}}}\nabla_{\mathbf{n}^\prime}\phi_{\text{f,f}}\right] - D_iB_{\text{f}}c_{\text{f,c}} + r_i(c_{j=J,...,M}).
    	\end{equation}
    Since $r_i$ may be non-linear, we linearise about an initial guess $\mathbf{c}^{N}=(c_{J}^{N},c_{J+1}^{N},...,c_i^N,...,c_{M}^{N})$ of concentrations and define $\delta c_{j} = c_{j}-c_{j}^{N}$:
    	\begin{equation}\label{eqn::LinearisedReactionRateUnrestricted}
    		r_{i}(c_{j=J,...,M}) = r_{i}(\mathbf{c^N})  + \delta c_J\partialderiv{r_i}{c_J}\bigg\rvert_{\mathbf{c}^N}  + \delta c_{J+1}\partialderiv{r_i}{c_{J+1}}\bigg\rvert_{\mathbf{c}^N} +...+ \delta c_{i}\partialderiv{r_i}{c_i}\bigg\rvert_{\mathbf{c}^N}  + ... + \delta c_{M}\partialderiv{r_i}{c_M}\bigg\rvert_{\mathbf{c}^N}  + \mathcal{O}(\delta c^2).
    	\end{equation}
    Substitute \cref{eqn::LinearisedReactionRateUnrestricted} into \cref{bcon::unresGenReacCsf} and take $c_{\text{s}}=c_{\text{s}}^{N+1}$ to be the next iterative solution. To decouple the iterative conditions, all other concentrations we set as $c_j=c_j^N$ for $i\neq j$. The result is the following iterative Robin condition for $c_{\text{s}}^{N+1}$, the concentration of species $i$ on the solid side of $\Gamma$:
    	\begin{align}
    		&-D^\star_{\text{s}} \nabla_{\mathbf{n}}c_{\text{s}}^{N+1} = -K^\star_{\text{s}} c_{\text{s}}^{N+1} - F^\star_{\text{s}}, \\
    		&D^{\star}_{\text{s}} = D_i, \\
    		&K^\star_{\text{s}} =   \bu _{\text{f}}\cdot\mathbf{n}^\prime - D_iB_{\text{f}} - \frac{z_iD_iF}{RT}\nabla_{\mathbf{n}^\prime}\phi_{\text{f}} - \partialderiv{r_i}{c_i}\bigg\rvert_{\mathbf{c}^N},  \\
    		&F^{\star}_{\text{s}} = D_iB_{\text{f}}c_{\text{f,c}} - r_i(\mathbf{c}^N) +c_{\text{s}}^N \pardiff{r_i}{c_i}\bigg\rvert_{\mathbf{c}^N}.
    	\end{align} 
    The same process can be done to determine the iterative Robin condition for $c_{\text{f}}=c_{\text{f}}^{N+1}$, the concentration of species $i$ on the fluid side of $\Gamma$:
        \begin{align}
            &-D^\star_{\text{f}} \nabla_{\mathbf{n^\prime}}c_{\text{f}}^{N+1} = -K^\star_{\text{f}} c_{\text{f}}^{N+1} - F^\star_{\text{f}}, \\
            &D^{\star}_{\text{f}} = D_i\,, \\
            &K^\star_{\text{f}} = \bu_{\text{f}}\cdot\mathbf{n}^\prime - D_iB_{\text{s}} - \frac{z_iD_iF}{RT}\nabla_{\mathbf{n}^\prime}\phi_{\text{f}} - \partialderiv{r_i}{c_i}\bigg\rvert_{\mathbf{c}^N}, \\
            &F^\star_{\text{f}} = D_iB_{\text{s}}c_{\text{s,c}} - r_i(\mathbf{c}^N) +c_{\text{f}}^N \pardiff{r_i}{c_i}\bigg\rvert_{\mathbf{c}^N} .
        \end{align}

    
	\bibliography{references}

\end{document}